\magnification=1100

\hsize 17truecm
\vsize 23truecm

\font\twelvec=msbm10 at 10pt
\font\sevenc=msbm10 at 7pt
\font\fivec=msbm10 at 5pt

\newfam\co
\textfont\co=\twelvec
\scriptfont\co=\sevenc
\scriptscriptfont\co=\fivec

\def\char{\mathop{\rm Char}\nolimits}

\def\Const{\mathop{\rm Const.}\nolimits}
\def\det{\mathop{\rm det}\nolimits}
\def\divr{\mathop{\rm div}\nolimits}
\def\exp{\mathop{\rm exp}\nolimits}
\def\Exp{\mathop{\rm Exp}\nolimits}
\def\dist{\mathop{\rm dist}\nolimits}

\def\Id{\mathop{\rm Id}\nolimits}

\def\ind{\mathop{\rm ind}\nolimits}
\def\Jac{\mathop{\rm Jac}\nolimits}

\def\Hess{\mathop{\rm Hess}\nolimits}

\def\lim{\mathop{\rm lim}\nolimits}

\def\pr{\mathop{\rm pr}\nolimits}

\def\rank{\mathop{\rm rank}\nolimits}

\def\supp{\mathop{\rm supp}\nolimits}

\def\Tr{\mathop{\rm Tr}\nolimits}

\def\sgn{\mathop{\rm sgn}\nolimits}

\def\WF{\mathop{\rm WF}\nolimits}

\def\Sum{\displaystyle\sum}
\def\e{\mathop{\rm \varepsilon}\nolimits}

\baselineskip 15pt
\centerline{\bf Semiclassical Green functions and Lagrangian intersection.}

\centerline{\bf Applications to the propagation of Bessel beams in non-homogeneous media}

\bigskip
\centerline{Michel ROULEUX}
\bigskip
\centerline{Universit\'e de Toulon, Aix Marseille Univ, CNRS, CPT, Toulon, France}

\centerline{rouleux@univ-tln.fr}
\medskip
\noindent {\it Abstract}:
We study semi-classical asymptotics for problems with localized right-hand sides by considering a Hamiltonian
$H(x,p)$ positively homogeneous of degree $m\geq1$ on $T^*{\bf R}^n\setminus0$. The energy shell is $H(x,p)=E$,
and the right-hand side $f_h$ is microlocalized: (1) on the vertical plane $\Lambda_0=\{x=x_0\}$;
(2) on the ``cylinder'' $\Lambda_0=\{(X,P)=\bigl(\varphi\omega(\psi),\omega(\psi)\bigr); \ \varphi\in {\bf R}, \omega(\psi)=(\cos\psi,\sin\psi)\}$.
when $n=2$. Most precise results are obtained in the isotropic case $H(x,p)={|p|^m\over\rho(x)}$, with $\rho$
a smooth positive function.
In case (2), $\Lambda_0$ is the frequency set of Bessel function $J_0({|x|\over h})$, 
and the solution $u_h$ of $(H(x,hD_x)-E)u_h=f_h$ when $m=1$, already provides an insight in the structure of  ``Bessel beams'',
which arise in the theory of optical fibers. We present in this work some extensions of [AnDoNaRo]. In Sect.3 we sketch
the semi-classical counterpart of the construction of parametrices for the Cauchy problem
with Lagrangian intersections, as is set up by R.Melrose and G.Uhlmann. This involves Maslov {\it bi-canonical operator}. 
\medskip
\noindent {\bf 1. Introduction}
\medskip
Let $M={\bf R}^n$, or possibly a smooth manifold. Write $\langle\theta\rangle=(1+\theta^2)^{1/2}$. For $m,k\in{\bf Z}$, 
we recall (see e.g. [Iv,Chapt.1]) the usual class of symbols  
$$S^k_m(M\times M\times {\bf R}^N)=\{a(x,y,\theta;h)\in C^\infty: |\partial_{(x,y)}^\alpha\partial_\theta^\beta a(x,y,\theta;h)|\leq
C_{\alpha,\beta}h^k \langle \theta\rangle^{m-|\beta|}\}$$
with asymptotic expansion $a(x,y,\theta;h)\sim h^k\bigl(a_0(x,y,\theta)+ha_1(x,y,\theta)+\cdots\bigr)$.
We define analogously $S^k_m(T^*M)$. Most of the time, we shall consider symbols compactly supported in $\theta$. 
Let 
$H(x,hD_x;h)$ be a $h$-PDO whose symbol $H$ belongs to $S^0_m(T^*M)$, or $S^0_m(T^*M\setminus0)$, which is for instance the case when
$H(x,p;h)$ is positively homogeneous of degree $m$ with respect to $p$.  We will denote by $z=(x,p)$ or by $z=(x,\xi)$ a point in $T^*M$.
Let $E\neq0$
be a non critical energy level for Hamiltonian $H_0$. 
Let also $f_h$ be a Lagrangian semi-classical distribution locally of the form 
$$f_h(x)=(2\pi h)^{-n/2} \int_{{\bf R}^n} e^{i\varphi(x,\theta)/h}a(x,\theta;h)\,d\theta\leqno(1.1)$$
where $\varphi$ is a non-degenerate phase function in the sense of H\"ormander
(or a generating family in the sense of Arnold) defining the Lagrangian manifold $\Lambda_0$ (see Sect.3), and 
$a\in S^k_0(M\times{\bf R}^n)$ for some $k$. We assume here $a$ to be compactly supported in $\theta$, excluding e.g. the case $a=1$,
$\varphi(x,\theta)=\langle x,\theta\rangle$, corresponding to $f_h(x)=\delta(x)$, which is too singular from our point of view. 
For brevity we will often denote such a (normalized) integral, possibly including Maslov index factor such as $e^{i\pi n/4}$, 
simply by $\int^*(\cdots)$. 

We shall assume that Hamiltonian vector field $v_{H_0}$ is transverse to $\Lambda_0$,
which we call {\it Lagrangian intersection}. In [AnDoNaRo] we have considered in the context of Maslov canonical operator
the problem of ``semi-classical Green functions'', which 
consists in solving $(H-E)u_h=f_h$, $\Lambda_0$ being the vertical plane $x=x_0$,
with $u_h$ {\it outgoing at infinity}. 
The distributions $u_h$ and $f_h$ are linearly related by
$u_h=E_+(h) f_h$, where $E_+(h)$ is the semi-classical {\it outgoing  parametrix}, that we shall compute in term of Maslov canonical operators.
\medskip
\noindent {\it 1.1 Some examples of ``localized functions''}.
\smallskip
Here are some examples of $f_h$ (expressed in a single chart):
\smallskip
\noindent ${\bf 1}$. $\Lambda_0=\{x=x_0\}=T^*_{x_0}M$, $x_0\in M$ (that we call the ``vertical plane'') is the conormal bundle to $\{x_0\}$, so that 
$$f_h(x)=\int^* e^{i(x-x_0)p/h}a(x,p;h)\,dp\leqno(1.2)$$
We say simply that $f_h$ is a ``localized function'' at $x_0$. 
This is the basic example since $\Lambda_0$ can alway be taken microlocally to such a form, and $H(x,hD_x;h)$ to $hD_{x_n}$,
where $x=(x',x_n)$, see Sect.3. Note that we can choose the amplitude $a$ in (1.2) independent of $x$, see [H\"o,Lemma 18.2.1]. 
\smallskip
\noindent ${\bf 2}$. More generally, a conormal distribution 
$f_h(x)=\int^* e^{ix'p'/h}a(x,p';h)\,dp'$, 
$x=(x',x''), p=(p', p'')$, with respect to $N=\{x'=0\}$, i.e. $\Lambda_0=T^*_N{\bf R}^n=\{(x'',p')\}$. Actually $f_h$ can be expressed with a new amplitude
$a(x'',p';h)$ not depending on $x'$, see [H\"o,Lemma 18.2.1]. 
\smallskip
\noindent ${\bf 3}$. WKB functions in Fourier representation 
$$f_h(x)=\int^* e^{i(xp+S(p))/h}a(x,p;h)\,dp$$
here $\Lambda_0=\{(-\partial_pS(p),p): p\in {\bf R}^n\}$. 
\smallskip
\noindent ${\bf 4}$. Semi-classical distributions related with Bessel functions, microlocalized on 
$$\Lambda_0=\{x=X(\varphi,\psi)=\varphi\omega(\psi), p=P(\varphi,\psi)=\omega(\psi), \varphi\in{\bf R}\}\leqno(1.3)$$ 
which is called the ``cylinder''~; here $\omega\in {\bf S}^{n-1}$ is the unit vector parametrized by $\psi$, see Sect. 5.

When $n=2$, this holds in particular for Bessel function of order 0
$$f_h(x)=\sqrt{2\pi\over h}J_0\bigl({|x|\over h}\bigr), \quad x\in{\bf R}^2$$
a radially symmetric solution of Helmholtz equation $-h^2\Delta v-v=0$, and follows from the integral representation
$$J_0\bigl({|x|\over h}\bigr)={1\over2\pi}\int_0^{2\pi}e^{i\langle\omega(\psi),x\rangle/h}\,d\psi$$
More generally, this applies to ``regular'' distributions on ${\bf R}^2$ of the form
$$f_h(x)=(2\pi h)^{-1/2}\int_0^{2\pi} e^{i\langle \omega(\psi),x\rangle/h}a(\langle\omega(\psi),x\rangle, \psi)\,d\psi\leqno(1.4)$$
with an amplitude of the special type $a(\langle\omega(\psi),x\rangle, \psi)$.
They are similar to conormal distributions in Example (2) after $x'$-variables has been removed from the amplitude.
See e.g. [DoMaNaTu], [DoNaSh], Sect.1.5.2, [DoMiNa], Sect.2.

When $n=3$ this holds also [DoMaNa] for semi-classical distributions of the form
$$f_h(x)=(2\pi h)^{-1}\int e^{i\langle \omega(\psi),x\rangle/h}\sin\theta\,d\theta\,d\widetilde\psi
=(2\pi h|x|)^{-1/2}J_{1/2}\bigl({|x|\over h}\bigr), \quad x\in{\bf R}^3$$
with $\psi=(\widetilde\psi,\theta)$ the standard angular coordinates.
One can also stick in an amplitude $a(\langle\omega(\psi),x\rangle, \psi)$
as in (1.4). 

\smallskip
\noindent ${\bf 5}$. $f_h$ identifies with a ``Bessel beam'' when
$$\Lambda_0=\{(x,p)\in T^*{\bf R}^3: \ x={\varphi\omega(\psi)\choose\phi}, \ p={\lambda(\phi)\omega(\psi)\choose\varphi\lambda'(\phi)+k}, 
\omega(\psi)={\cos\psi\choose\sin\psi}, \ \phi,\varphi\in{\bf R}\}$$
$k\in{\bf R}$ and  $\lambda$ is a smooth positive function. Such Lagrangian distributions, with the notation $x=(x',x_3)$, take the form
$$f_h(x)=(2\pi h)^{-1/2}e^{ikx_3/h}\int_0^{2\pi}e^{i\lambda(x_3)\langle x',\omega(\psi)\rangle/h} a(\lambda(x_3)\langle x',\omega(\psi)\rangle+kx_3,x_3,\psi)\,d\psi$$
When $a=a(\tau,\alpha,\phi)=\widetilde a(\alpha,\phi)$ and is even in $\alpha$, we have
$$f_h(x)={\sqrt{2\pi\over h}}\widetilde a(|x'|,x_3)e^{ikx_3/h}J_0\bigl({\lambda(x_3)|x'|\over h}\bigr)+{\cal O}(h)$$
which justifies the name ``Bessel beams'', see [DoMaNa2].  
 
\noindent ${\bf 6}$. Airy-Bessel beams, which are known also as Berry-Balasz solution 
([BeBa], [DoMaNa1]) of the paraxial approximation of wave equation in 3-D, with initial manifold
$$\Lambda_0=\{(x,p)\in T^*{\bf R}^3: \ x={\varphi\omega(\psi)\choose\phi^2/2}, \ p={\lambda\omega(\psi)\choose\phi}, \ 
\omega(\psi)={\cos\psi\choose\sin\psi}, \ \phi,\varphi\in{\bf R}\}$$

\noindent ${\bf 7}$. Lagrangian distributions with a complex phase
in the sense of Melin-Sj\"ostrand, see [MeSj]  
equivalently a complex germ in the sense of Maslov, which are superposition of 
coherent states $f_{(x_0,\xi_0)}(x;h)={1\over h^n}\exp(-\omega^2\cdot(x-x_0)^2/2h)\exp(ix\xi_0/h)$, 
and $\Lambda_0$ is a strictly positive Lagrangian manifold.  

Examples (1) and (4) will be extensively studied when $n=2$ and $H_0$ is positively homogeneous of degree $m\geq1$ ($m=1$ in Example (4)).
Since $f_h$ in Example (5) looks like a plane wave in direction $x_3$, one could expect Example (4) could be generalized to this case,
once the eikonal coordinate has been found. 
Examples (6) and (7) 
requires some special treatment and will not be considered either. 
\medskip
\noindent {\it 1.2 Global parametrices for PDO's of principal type and their semi-classical counterpart}.
\smallskip
Thus our main problem is to represent the formal asymptotic solution of 
$$(H(x,hD_x;h)-E)u_h(x)=f_h(x), \ u_h(x)=E_+(h)f_h(x)=\int_0^\infty e^{-it(H-E)/h}f_h(x)\,dt\leqno(1.5)$$
in the most explicit form. Here $e^{-itH/h}$ is the propagator, so we need consider first Cauchy problem
$$hD_tv_h+(H(x,hD_x;h)-E)v_h=0, \ v_h|_{t=0}=f_h\leqno(1.6)$$
For general $H(x,hD_x;h)$ we adapt the constructions of [MelUhl] to the semi-classical setting~; in the case of homogeneous Hamiltonians
of degree $m\geq1$ we present formulas more globally.
See [BoRo] when $\Lambda_0$ is given by (1.3), which forces $m=1$. 

Here we review some general concepts for solving the semi-classical PDE's (1.5)-(1.6) in term of
oscillating functions (i.e. semi-classical distributions) modulo ${\cal O}(h^\infty)$ terms, 
as well as their non semi-classical analogue (for standard pseudo-differential calculus)
$$P(x,D_x)u(x)=f(x), \ u(x)=\int_0^\infty e^{-itP}f(x)\,dt$$
starting from Cauchy problem
$$D_tv+P(x,hD_x)v_h=0, \ v|_{t=0}=f$$
In the latter case we look for $C^\infty$ parametrices, i.e. distributions defined modulo smooth functions.

Semi-classical approximation is also called asymptotics with respect to the small parameter $h$.
This is the most natural calculus, since it is concerned with functions oscillating rapidly with respect to a given scale $h$,
but can nevertheless be smooth in $x$. Their are called {\it semi-classical distributions}.
In the simplest case these semi-classical distributions are of WKB type. 

Global existence of an outgoing solution at infinity
provided suitable hypotheses on Hamilton vector flow, such as Lagrangian intersection, the non-trapping
and the non-return conditions, is quite involved. 
The main strategy in case of conormal distributions, has already been set up in [MelUhl]
in the context of smooth parametrices (asymptotics with respect to smoothness) for a PDE, called sometimes Melrose-Uhlmann Calculus.

The non-return condition was first formulated in [GuMel], Proposition 2.7, 
in constructing parametrices for hyperbolic PDE's with a boundary (billiard problem),
extending the case without a boundary [DuH\"o]. 

The non-return condition was then formulated in greater generality in [MelUhl], Eq.(6.5), when solving a PDE with a right hand side. 
Let $\Lambda_0$ be  a conic Lagrangian manifold, which we will call the ``initial manifold'', and $P(x,D_x)$ be a PDO
of real principal type, with real principal symbol
$p(x,\xi)$ such that the Hamilton vector field $v_p$ is never tangent to $\Lambda_0\cap\Sigma(P)$, where $\Sigma(P)=\{p(x,\xi)=0\}$. Given
a Lagrangian distribution $f$ microlocally supported on $\Lambda_0$, 
the problem is to find $u\in {\cal D}'({\bf R}^n)$ such that 
$P(x,D_x)u=f$ mod $C^\infty$. Here $\partial\Lambda_+=\Lambda_0\cap\Sigma(P)$ plays the role of the ``boundary'' (generically both are of dimension $n-1$).
Note that the results of [DuH\"o]
and [MelUhl] hold in a general pseudo-convex manifold $M$,
and that in [DuH\"o] the non-return condition is automatically satisfied.
We refer to [MelUhl] for precise
statements, see also [GreUhl], [Jos], [ForHasHil], [SteSh] and references therein. 

The natural framework of such PDE's has its counterpart in the semi-classical case described in (1.5).
In particular the non-return condition (renamed as the {\it non-refocusing} condition) has received a more systematic treatment
in relatively recent works [Ca], [Bon], [KlCa].

Namely, let $\Lambda_0=T^*_{x_0}{\bf R}^n$ be the vertical plane, and $H(x,hD_x)=-h^2\Delta+V(x)$
be a semi-classical Schr\"odinger operator
with a smooth potential $V$ with long range interaction and such that $V(x_0)\leq0$.  
In this case, the non-refocusing condition is characterized by the relation on $T^*_{x_0}M$ (the ``return set'')
$${\cal R}=\{(p,\eta):  p^2+V(x_0)=\eta^2+V(x_0)=0; \ \exists t>0: \ X(t,p)=x_0, P(t,p)=\eta\}\leqno(1.7)$$
where $\bigl(X(t,p),P(t,p)\bigr)=\exp tv_H(x_0,p)$ is the trajectory issued from $(x_0,p)\in T^*M$.

The non-refocusing condition in the restricted sense, means that ${\cal R}=\emptyset$.
But the dimension of ${\cal R}$
can also be taken into account. In [Ca], [Bo], [KlCa] it is assumed that
$\widetilde{\cal R}=\{(p,\eta,t)\in {\bf R}^{2n+1}: (p,\eta)\in{\cal R}\}$
is a submanifold of ${\bf R}^{2n+1}$ of dimension less than $n-1$. 

The {\it non trapping condition} should also be introduced in the semi-classical setting (see also [MelUhl], Eq.(6-3)-(6-4)). 
This is a condition on the set of trapped trajectories at energy $E_0$:
$${\cal K}(E_0)=\{(x,\xi): \ H_0(x,\xi)=E_0, \ (X(t),P(t))\ \hbox{does not tend to infinity as} \ |t|\to\infty\}=\emptyset $$
where $(X(t),P(t))=(X(t,x,\xi),P(t,x,\xi))=\exp tv_{H_0}(x,\xi)$
and here we can replace the non-trapping condition on the phase variables by a condition on $X(t)$ alone, e.g. 
$|X(t)|=|X(t,x,\xi)|\to \infty$ as $|t|\to\infty$. 

Outgoing solutions $u_h$,  are characterized by Sommerfeld radiation condition of the form
$${x\over|x|}\nabla_x w(x)+i\sqrt{-V(x_0)}w(x)\to0, \ |x|\to\infty, \quad n\geq2\leqno(1.8)$$
where $w=\lim_{h\to0}w_h$,  $w_h(x)=h^{n/2}u_h(hx)$, is the unique solution of $(-\Delta+V(x_0))w=f$,
and $f_h(x)=h^{-n/2}f\big({x\over h}\bigr)$.
It relates in a non trivial way the behavior of $u_h$ at infinity 
with the value of the potential $V$ at $x_0$. 
Sommerfeld radiation condition requires careful estimates on $U_h(t)=e^{-itH/h}$, 
or $U_h(t)f_h$, along with a discussion according to the relative magnitude
of $t$ and $h$. 
The proof consists, roughly speaking, in testing $U_h(t)f_h$ against some fixed $\phi\in{\cal S}({\bf R}^n)$, and then
show that $\langle u_h,\phi\rangle\to
\langle w,\phi\rangle$ as $h\to0$.
In particular one needs to know asymptotics of
$u_h$ in a $h$-dependent neighborhood of $\Lambda$. 

In this paper instead, given $H(x,hD_x;h)$ and the initial Lagrangian manifold
$\Lambda_0$, we content ourselves to present, in the sense of {\it formal asymptotics}, a ``closed form'' for
the solution of (1.5)
in term of Maslov canonical operator for bi-Lagrangian distributions. So the non-refocusing and the non-trapping conditions
can be largely ignored in case of formal asymptotics.
By formal asymptotics [Ler] we mean that, in principle, our approximate solution is only a {\it quasi-mode}, i.e.
it has no reason to be equal to $E_+(h)f_h$ mod ${\cal O}(h^2)$.
In practice however numerical simulations show that Maslov canonical operator provides an excellent agreement with the ``exact solution''.

Actually the main issue we are faced with is about existence and
uniqueness of the ``asymptotic solution''. To fix the ideas, we can already formulate the problem as follows:
how can a ``formal'' WKB solution approximate the solution of (1.6)? A possible answer relies on the construction of a normal form for $H(x,hD_x,h)$.
Assume for instance $H(x,hD_x,h)$ is self-adjoint and such that its the principal symbol $H_0(x,p)$ is real,
vanishes at $(x_0,p_0)$ and $dH_0(x_0,p_0)\neq0$ (by analogy with the terminology used for smoothing parametrices, we call $H(x,hD_x;h)$ a $h$-PDO
of real principal type near $(x_0,p_0)$). Then
$H(x,hD_x;h)$ is microlocally equivalent to $hD_{y_n}$, i.e. there is a FIO $U$ associated with the canonical transformation $\kappa$
such that $\kappa(0;0,\cdots,0,1)=(x_0,p_0)$ (Darboux theorem), $U$ is microlocally unitary, and $\|hD_{y_n}-U^*H(x,hD_x,h)U\|={\cal O}(h^\infty)$,
where $\|\cdot\|$ is a local operator norm. 
This is a version of the semi-classical Egorov Theorem, see e.g. [Iv,Sect.1.2]. So solving (1.6) amounts to 
construct a solution of
$$hD_tv'_h+hD_{y_n}v'_h=0, \ v'_h|_{t=0}=f'_h$$
where $f'_h=U^*f_h$. This PDE with constant coefficients can easily solved in the ``same class'' as $f'_h$. Taking the image of $v'_h$
by $U$ gives the suitable WKB solution of (1.6). Of course, this construction can be extended so long as $H(x,hD_x,h)$ remains of real
principal type. This is one of the basic ideas of [MelUhl], which has to be extended in the case of a (microlocal) boundary,
namely to the solution of (1.5). See also [SteSh]. Another ingredient for solving (1.5) is to construct an asymptotic solution in the
``elliptic zone'', i.e. when $H_0(x,p)-E\neq0$. There we only need to invert an elliptic $h$-PDO. Gluing together the different branches of
solutions follows from the ``compatibility condition'' (see [MelUhl] and Sect.3 below for the semi-classical case), leading here to the
notion of {\it Maslov bi-canonical operator}. 

In the framework of Maslov bi-canonical operator we can, in principle, makethese constructions global provided the non-trapping and the non-refocusing
conditions hold.

Note that making use of the non-trapping condition only, we can ensure that our construction of Maslov bi-canonical operator can be extended
to large values of $x$, but still microlocally outside the initial manifold $\Lambda_0$ (that is, in the case $\Lambda_0$ is the vertical plane
$x=x_0$, outside $x=x_0$.)

Even if we had taken care of the non-trapping and the non-return conditions, we should point out that the solution to (1.5) 
is not unique in general. Following the {\it principle of limiting absorption} we should introduce the auxiliary equation
$$(H(x,hD_x;h)-E+i\e )u_{h,\e }(x)=f_h(x), \ u_{h, \e }(x)=E_+(h,\e )f_h(x)=\int_0^\infty e^{-itH/h}e^{it(E+i\e )/h}f_h(x)\,dt$$
and take the limit $\e \to0$. In case $H(x,hD_x;h)$ is Schr\"odinger operator, this is related with the fact
that the limit of $u_{h,\e }(x)$, when $\e \to0$ satisfies  Sommerfeld radiation condition (1.8),
but in general taking the limit $\e \to 0$ can be a non trivial fact. 
So for simplicity again, we will ignore limting absorption principles in this work.

To close this general introduction, we should mention again the case of 
parametrices for hyperbolic PDE's with a boundary (billiard problem), which was given a new insight in [PeSto], [PeVo], [Vo]. 
These works make use of resolvent estimates. 

\medskip
\noindent {\it 1.3 Lagrangian intersection and microlocal Green functions}.
\smallskip
Recall the canonical 1-form on $T^*M$ takes (locally) the form $p\,dx$; on the extended phase-space $T^*(M_x\times{\bf R}_t)$,
(or locally on $T^*{\bf R}^{n+1}$) it is given by $p\,dx-E\,dt$. We consider several Lagrangian manifolds~:
\smallskip
\noindent $\bullet$ The ``initial manifold'' $\Lambda_0\subset T^*{\bf R}^n$ which contains $\WF_hf_h$. We will be concerned essentially with
the ``vertical plane'' $\Lambda_0$ (1.2) or Bessel cylinder (1.3). 

\noindent $\bullet$
The Lagrangian submanifold in extended phase-space $T^*{\bf R}^{n+1}$ (denoted instead by $\Lambda$ in [AnDoNaRo3])
$$\widetilde\Lambda_+=\{ (x,p) = g^t(z), t>0, E=H(z) : z \in \Lambda_0 \}\leqno(1.15)$$
where $g^t=\exp tv_H$ is the phase flow of the Hamiltonian vector field $v_H$ generated by $H$
$$\dot{x} = \partial_{p} H(x,p), \quad \dot{p}= - \partial_{x} H(x,p)$$
(denoting $H$ instead of $H_0$). Not to confuse $E$ as a variable in (1.15) with the given value of energy in (1.5), we shall
change the variable $E$ to $E-\tau$, so that $H=E$ for $\tau=0$.
Lagrangian manifold $\widetilde\Lambda_+$ contains
$\WF v_h(x,t)$ where $v_h$ solves cauchy problem (1.5).

\noindent $\bullet$ The Lagrangian submanifold of $T^*{\bf R}^{n}$ (denoted instead by $\Lambda_+$ in [AnDoNaRo3])
$$\Lambda_+^E=\{z \in T^*{\bf R}^n : \exists \, t\geq0, \, z \in \exp tv_{H}(\Lambda_0 \cap \Sigma_E)\}\leqno(1.16)$$
where $\Sigma_E$ is the energy surface $H=E$. Let $\partial\Lambda_+^E=\Lambda_+^E\cap\Lambda_0$ be the ``boundary'' of $\Lambda_+^E$.
We shall always assume Lagrangian intersection, i.e. $v_H$ is transverse to $\partial\Lambda_+^E$. 

For Example (4) however, there may be points on $\Lambda_+^E$ 
where transverse intersection fails (we call them {\it glancing points}) 
by analogy with the problem of diffraction by obstacles,
see Sect.5), so we miss some informations on $u_h$ nearby. Note that a glancing point $z_0\in\Lambda_0$
is a critical point of $H|_{\Lambda_0}$. This is discussed in [BoRo].

We shall assume that the amplitude $a$ defining $f_h$ is compactly supported in $p$. This hypothesis is discussed in more detail in
[AnDoNaRo3], Sect.1.6. In particular, if $f_h$ is only rapidly decreasing in $p$, then we must assume some ellipticity of $H_0$ at $|p|=\infty$.

For simplicity we restrict 
to the case where $\partial\Lambda_+^E$ is a compact, isotropic submanifold without boundary, which is certainly the case
when $H_0(x,p)$ is elliptic. This restriction however is not essential [MelUhl] and some our results carry to the case where $H(x,hD_x;h)$
is the wave operator.

We also assume that there is no finite motion on $\Sigma_E$, namely
$$|z(t)|\to\infty \ \hbox{as} \ t\to\infty\leqno(1.20)$$
for any trajectory $z(t)=(X(t),P(t))$ issued from $\partial\Lambda_+$ at $t=0$.

The non-return set ${\cal R}$ in (1.7) is irrelevant if
we content ourselves with the asymptotics of $u_h(x)$ microlocally in a compact set outside $\partial\Lambda_+^E$. 
For instance if $\Lambda_0=T_{x_0}^*M$ we shall compute $u_h(x)$,
locally uniformly in any compact set $K\subset M\setminus\{x_0\}$, as $h\to0$. A first improvement would consist in removing only a $h^\delta$-neighborhood
of $x_0$ for some $0<\delta<1$, this we have sketched in [AnDoNaRo1,Thm.2]. 

Our main goal is to represent formally (1.5) as a superposition 
$$u_h(x)={i\over h}\int_0^\infty \bigl[K_{\widetilde \Lambda_+}^hb\bigr](x)\,dt+{\cal O}(h)\leqno(1.22)$$
where $[K_{\widetilde \Lambda_+}^hb](x,t)$ is Maslov canonical operator associated with $\widetilde\Lambda_+$, and $b$ an amplitude depending linearly from
the amplitude $a$ defining $f_h$.

The relevant contributions to this integral come from $t=+\infty$, $t=0$ and the critical points $t\in]0,+\infty[$. 

The contribution of $t=+\infty$ is excluded by (1.20). If we do not have glancing
points then $\Lambda_+^E$ is smooth Lagrangian submanifold. 
We need to take into account $t=0$ in the formula (1.22). If the initial manifold $\Lambda_0$ is the vertical
then this contribution is zero outside of $x_0$. Otherwise (as is the case $\Lambda_0$ is
Bessel cylinder) the contribution of $t=0$ is not zero but maybe of
small order (as $\sqrt{h}$). But if $f_h$ has compact
support, then the contribution of
$t=0$ is zero outside $\supp (f_h)$.

The critical points $t\in]0,+\infty[$ give of course the main contributions to (1.22), and give in principle all possible types of
Lagrangian singularities in ${\bf R}^n$, microlocally outside
$\Lambda_0$. This will be investigated in detail in Sect.4 and 5 for a particular type of Hamiltonian which
allows to contruct quite explicit ``global'' asymptotic solutions. 

In case $\Lambda_0=T^*_{x_0}M$ is the vertical plane, (1.20) can be replaced by $|X(t)|\to\infty$ as $t\to+\infty$),
and the measure on $\Lambda_+$ factorizes as $d\mu_+=d\mu_0^E\wedge dt$,
where $d\mu_0^E$ is the measure on $\partial\Lambda_+^E$, see [AnDoNaRo3,Thm 2].
The accuracy ${\cal O}(h)$ in (1.22) can improve to ${\cal O}(h^\infty)$, see [AnDoNaRo3,Thm 4], and also [SteSh]. 

One of our claims is to express $\int_0^\infty dt\bigl[K_{\Lambda_+}^hb]$ as some ``bi-canonical
operator'' $[K_{\Lambda_0,\Lambda_+}^h(\sigma,\sigma^+)]$ acting on pairs of symbols $(\sigma,\sigma^+)$ depending linearly on $b$,
resp. the {\it boundary-part} and the {\it wave-part} symbol of $u_h$, which we call a {\it bi-Lagrangian} (semi-classical) distribution.
This follows from a symbolic semi-classical calculus similar to [MelUhl].

For simplicity, we shall ignore throughout Maslov indices. Corrections to the asymptotics due to Maslov indices follow for instance easily from
[Ar], [So], [Iv], [BaWe, Sect.4], [DoNaSh], [DoRo], [EstHagHedLitt]. See also [AnDoNaRo3].

\medskip
\noindent{\it 1.4 Main Results}.
\smallskip
\noindent $\bullet$ {\bf General case}. 
\smallskip
In Sect.3, we translate the setting of classical Pseudo-differential Calculus
elaborated in [MelUhl] to the semi-classical one, and make some statements more precise. 
By Proposition 3.1 below, near any $z\in \partial\Lambda_+^E$ we are reduced microlocally to the case where $\Lambda_0=T_0^*{\bf R}^n$, 
and $\Lambda_+^E={\cal L}_+^0$ is the flow out of 
the ``model Hamiltonian'' $H_0=\eta_n$ in energy surface $\eta_n=0$. In Proposition 3.4 we show that
$H(x;hD_{x};h)-E$ can be taken microlocally
near ${\cal L}_+^E$ to its normal form $hD_{y_n}$ by conjugating with a $h$-FIO. So without loss of generality, 
we may assume that, after some canonical transformation,  $\Lambda_0=T^*_{x_0}M$. 
We have~:
\smallskip
\noindent {\bf Theorem 1.2}: {\it Let $H(x,p;h)\in S^0_m(T^*M)$ (not necessarily homogeneous).
Assume the vertical plane $\Lambda_0=T_{x_0}M$, and $\Sigma_E=\{H_0(x,p)=E\}$, where $E$
is a non critical energy level, intersect transversally along the compact isotropic manifold $\partial\Lambda_+^E$. Assume
$|X(t)|\to\infty$ as $t\to\infty$ for all initial conditions $(X(0),P(0))\in\partial\Lambda_+^E$. Let also $f_h$ be a semi-classical distribution
microlocalized on $\Lambda_0$ of the form (1.2), with $x_0=0$ (to fix the ideas), and where we can assume the symbol is independent of $x$.
Then there is $\e _0>0$ such that the equation $(H(x,hD_x)-E)u_h(x)=f_h(x)+{\cal O}(h^2)$ can be solved in the form (1.22) for $0<|x-x_0|<\e _0$,
locally uniformly for $h>0$ small enough. 

Moreover there are symbols $(\sigma,\sigma^+)$, such that 
in local coordinates $\sigma(\xi;h)={1\over\xi_n}b(\xi;h)$, $\sigma^+(\xi;h)=2i\pi b(\xi',0;h)$
verify the compatibility condition 
$$\lim_{\xi_n\to 0}\xi_n\sigma(\xi',\xi_n;h)=\sigma^+(\xi',0;h)\leqno(1.23)$$ 
and such that $u_h$ can be represented by Maslov ``bi-canonical operator'', see (3.33) as
$$u_h(x)=[K_{\Lambda_0,\Lambda_+^E}^h(\sigma,\sigma^+)](x)$$
We have the commutation relation
$$H(x,hD_x;h)\bigl[K_{\Lambda_0,\Lambda_+^E}^h (\sigma,\sigma^+)\bigr](x;h)=\bigl[K_{\Lambda_0,\Lambda_+^E}^h 
((H-E)\sigma, 0)\bigr](x;h)+{\cal O}(h^2)
=f_h(x)+{\cal O}(h^2)\leqno(1.24)$$
for $0<|x-x_0|<\e _0$,
locally uniformly for $h>0$ small enough.}
\medskip

We stress that this representation holds in a punctured neighborhood of $x_0$, i.e. is not uniform in $x$ near $x_0$,
since we have neglected the non-return condition.

However, since the Hamiltonian vector field $v_H$ remains transverse to $\Lambda_t=g^t(\Lambda_0)$ for all $t$,
the condition $|x-x_0|<\e _0$ is not really a restriction. Namely, using local charts for
Maslov ``bi-canonical operator'' as in the case of Maslov canonical operator associated with the homogeneous equation (see Sect. 3.2 below),
the asymptotic solution $u_h$ can be made global. See also [AnDoNaRo3,Thm.4]. 

The remainder term ${\cal O}(h^2)$ in (1.24) can be improved by considering higher order
transport equations and using the compatibility condition (1.23).

Maslov ``bi-canonical operator'' is some operator version of the ``compatibility condition'' in [MelUhl], Formula (2.10) and Sect.4.
It encompasses the parametrix both in the elliptic region $H(x,p)\neq E$ and in the hyperbolic region $H(x,p)=E$
(see Theorems 1-2 in [AnDoNaRo3]).
\medskip
\noindent $\bullet$ {\bf Hamiltonian $H_0$ is homogeneous of degree $m$ with respect to $p$ and $\Lambda_0$ is the vertical plane}. 
\smallskip
The previous result is not very useful from the point of vue of applications, since it does not provide a ``close form'' for $u_h$.
Much more information is available when Hamiltonian $H_0$ is homogeneous of degree $m$ with respect to $p$, due to the relations
$\langle P(t,\psi),\dot X(t,\psi)\rangle=mH$ and $\langle P(t,\psi),X_\psi(t,\psi)\rangle=0$ 
(Huygens principle). We have chosen coordinates $(\lambda,\psi)$ on $\Lambda_0$,
$\lambda$ a ``radial coordinate'' such that $\lambda=1$ on $H=E$
and $\psi\in{\bf R}^{n-1}$ coordinates 
along $\partial\Lambda_+^E$, we can think of as angles parametrizing the $(n-1)$-sphere.

Lagrangian intersection always holds on $\Lambda_0$, for if $\partial_pH(x_0,p)=0$, Euler identity gives
$0=\langle \partial_pH(x_0,p),p\rangle=mH=mE$ which contradicts $E\neq0$. We will therefore assume the non-trapping condition (in $x$)
$|X(t)|\to\infty$, $t\to+\infty$. 

Note that such a symbol is not suitable
for Pseudo-differential Calculus when $m$ is not an even integer,
because of the singularity at $p=0$, but this is harmless if $(x_0,0)\notin\WF _h(f_h)$.

A particular case of these Hamiltonians is the ``conformal metric'' given by 
$$H(x,p)=|p|^m{1\over \rho(x)}\leqno(1.25)$$ 
where  $\rho$ be a smooth positive function on $M$, $m\geq1$. In case $m=2$, $H(x,p)$ identifies, through Maupertuis-Jacobi correspondence,
with Helmholtz Hamiltonian in a non-homogeneous medium with refraction index $n(x)>0$ such that $n(x)^2=\rho(x)$. 
We shall pay a great attention to this Hamiltonian,
for which computations are most explicit. 

Recall that
a {\it focal point}
for a Lagrangian embedding $\iota:\widetilde\Lambda\to T^*M$ (where $\widetilde\Lambda$ is either $\Lambda$, $\Lambda_t=\exp tv_H(\Lambda)$, 
or $\Lambda_+^E$) is a point $z\in\widetilde\Lambda$ such that
$d\pi_x:\widetilde\Lambda\to M$ has rank $<n$. The set of focal points is denoted by ${\cal F}(\widetilde\Lambda)$.   
Recall also that there exists a covering of  $\widetilde\Lambda$ by
canonical charts $U$ where $\rank d\pi_x(z)\geq k$ for all $z\in U$. These $U$ for which $k=n$ are called {\it regular} charts, and those for
which $k<n$ {\it singular} charts. 

At least for $0<|x-x_0|$ small enough i.e. $0<t$ small enough, using at most 2 canonical charts (depending if $X_\psi(t,\psi)=0$ or not)
we can construct $v_h$
solution of Cauchy problem (1.6), and hence $u_h$ by integration with respect to $t$. To this end, 
we introduce eikonal coordinates (see Sect.2),
and a generating family for  $\Lambda_+^E$ as (see (4.5)):
$$\Phi(x,t,\psi,\lambda)=mEt+\lambda\langle P(t,\psi),x-X(t,\psi)\rangle$$
The ``$\theta$-variables'' are thus $(t,\psi,\lambda)$.
This defines the 1-jet along the critical set $x=X(t,\psi),\lambda=1$ of the solution of Hamilton-Jacobi equation associated to (1.6). 
This yields (see (3.8) and (4.14)) an
(inverse) density on $\Lambda_+^E$, given by its density $\mu_+=F\bigl[\Phi,dy]=mE\det(P,P_\psi)$ ($dy$ being Lebesgue measure on the critical set)
and non vanishing precisely when $\Phi(t,x,\psi,\lambda)$
is a non-degenerate phase function, see Propositions 4.1 and 4.4. 
This holds at least for small $t$. 
For larger $t$, the non-vanishing of $mE\det(P,P_\psi)$ will only be assumed.
\medskip
\noindent {\bf Theorem 1.3}: {\it Let $n=2$, and $H_0(x,p)$ be positively homogenenous of degree $m\geq1$.   
Let $r_0>0$ such that 
$$\forall t>0, \ \forall \psi\in{\bf R}: \ \bigl[|X(t,\psi)|<r_0 \ \Longrightarrow \det(P(t,\psi),P_\psi(t,\psi))>0\bigr]$$
(this holds for $t>0$ small).
Then there is $\e _0>0$ such that the equation $(H(x,hD_x)-E)u_h(x)=f_h(x)+{\cal O}(h^{3/2})$ can be solved in the form
$$u_h(x)={i\over h}\int_0^\infty\,dt\int e^{i\Phi(x,t,\psi,\lambda)/h}b(x,t,\psi,\lambda)\,d\psi\,d\lambda\leqno(1.26)$$
for $0<|x-x_0|<\e _{0}$,
locally uniformly for $h>0$ small enough. 

Moreover we can decompose 
$$\int_0^\infty \,dt\bigl[K_{\Lambda_t}^hb](x)=\int_0^\infty \,dt\bigl[K_{\Lambda_t}^h(\chi_1b)](x)+
\int_0^\infty \,dt\bigl[K_{\Lambda_t}^h(\chi_2b)](x)\leqno(1.27)$$ 
where $\chi_1+\chi_2=1$ is a partition of unity subordinated to a chart $U_1$  where $X_\psi\neq0$ and a (singular) chart
$U_2$ where $X_\psi=0$, $\chi_2\equiv1$ near $|X_\psi|\leq{1\over3}$, and $\chi_1\equiv1$ near $|X_\psi|\geq{2\over3}$,
both contributions being discussed
in Sect. 4.4 and 4.5 below. }
\medskip
We do not attempt here to formulate the result in terms of bi-canonical Maslov operator as in Theorem 1.2.
Note the loss of accuracy (${\cal O}(h^{3/2})$) instead of ${\cal O}(h^{2})$ in Theorem 1.2.
In Sect.2.3 we discuss (somewhat informally) the case where $x$ has several pre-images under $\pi_x:\Lambda_+^E\to M$.

\smallskip

By constructing (1.5)-(1.6) we mean also determining the Lagrangian singularities of $u_h$. They are revealed when reducing the 
number of ``$\theta$-variables'' in the oscillating integral
$\int_0^\infty dt\bigl[K_{\Lambda_t}^hb](x)$.
We get virtually any kind of Lagrangian singularity, but because of homogeneity of $H$
with respect to $p$, it is convenient to introduce another classification in $T^*M$, which clarifies the construction of $u_h$.
\medskip
\noindent {\bf Definition 1.4}: 
{\it Let $H$ be positively homogeneous of degree $m$ with respect to $p$. 
We call a point $z=(x,p)$ such that $-\partial_x H(z)\neq0$
an {\it ordinary} point if $\langle-\partial_x H(z),p\rangle\neq0$, and a {\it special} point otherwise. 
If $-\partial_x H(z)=0$ we call $z$ a {\it residual} point.}
\smallskip
Here are some Examples:

(1) For Tricomi Hamiltonian, $H(x,p)=x_2p_1^2+p_2^2$, 
the residual points are those for $p_1=0$, the special points those for $p_1\neq0$
but $p_2=0$, and the ordinary points those for $p_1p_2\neq0$. Tricomi operator is used as a model for diffraction.
This model extends to the 3-D case as $H(x,p)=-p_2^2+p_1p_3+x_2p_1^2$. 

(2) For M\'etivier Hamiltonian, $H(x,p)=p_1^2+(x_1^2+x_2^2)p_2^2$, the residual points are given by $p_2=0$ or $x=0$,
the special points by $p_2\neq0$ and $x\neq0$, but $\langle x,p\rangle=0$, 
and the ordinary points by $\langle x,p\rangle p_2\neq0$. M\'etivier operator provides a counter-example for analytic
hypo-ellipticity.

(3) Let $H$ be the ``conformal metric'' given by (1.25). 
The residual points are the critical points of $\rho$~; 
at a special point, $\langle \nabla\rho,p\rangle=0$, i.e. $v_H$ is tangent to the level curves of $\rho$.
This is our main example here.

\smallskip
We denote by ${\cal S}(\Lambda_+^E)$ the set of special points on $\Lambda_+^E$.  
We shall partition points $z\in\Lambda_+$
according to the following values: (1) $z$ is a focal (or non-focal) point; (2) $z$
is a special (or ordinary, or residual) point. Thus each canonical chart splits again into ordinary, special or residual points. 
Assume $n=2$ for simplicity.

Let $x=X(t,\psi)$, and $\langle-\partial_xH(x,p),\partial_pH(x,p)\rangle\neq0$ at $(x,p)=z(t)=(X(t,\psi),P(t,\psi))$. Then
$$\partial_t\Phi(t,x,\psi,\lambda=1)=0 \ \Longrightarrow \partial^2_t\Phi(t,x,\psi,\lambda=1)\neq0\leqno(1.31)$$ 
so that we can perform asymptotic stationary phase in $t$ to simplify (1.20) at $x=X(t,\psi)$. This holds when $H$ is of the form (1.25)
and $z(t)$ is an ordinary point.
Likewise, if several values of parameters $t$ contribute, we sum over such $t$'s.   
Then $\Phi(x,t,\psi,\lambda=1)$  reduces to a phase function $\Psi(x,\psi)$ and we can further reduce the number of variables
in a standard way, according to the fact that $z(t)$ is a focal point or not. 
When $\langle-\partial_xH(x,p),\partial_pH(x,p)\rangle=0$ but $-\partial_xH(x,p)\neq0$, the situation looks like (1.10) and we have
$$\partial_t\Phi(t,x,\psi,\lambda)=\partial_\lambda\Phi(t,x,\psi,\lambda)=0 \ \Longrightarrow \det\Phi''_{(t,\lambda),(t,\lambda)}(t,x,\psi,\lambda)\neq0
\leqno(1.32)$$ 
where $\lambda$ is constrained to be equal to 1 on the critical set $C_\Phi$, 
so that we can perform asymptotic stationary phase with respect to $(t,\lambda)$. Then $\pi_x:\Lambda_+^E\to T^*{\bf R}^n$ 
has rank 1 if $z(t)$ is a special point or rank 2 otherwise, see Proposition 4.3(ii). Note that $\Lambda_+^E$ never turns vertical, 
since $\partial_pH(x,p)\neq0$ on $H=E$.

In this generality, we only succeed (see Proposition 4.3) to describe the contribution
to (1.20) (by asymptotic stationary phase) of short times $t$  (``near field''), i.e. so long as $F\bigl[\Phi,dy\bigr]=mE\det(P,P_\psi)$
$\neq0$,
But this is actually sufficient to compute $u_h$ microlocally near $\Lambda_0$, when $x\neq x_0$. 
\medskip
\noindent $\bullet$ {\bf $H_0$ is the ``conformal metric'' $H_0(x,p)={|p|^m\over\rho(x)}$ and $\Lambda_0$ is the ``vertical plane''}. 
\smallskip
Using that $P(t,\psi)$ is parallel to $\partial_p H(x,p)$ we get more complete results in this case.  
If $\rho$ is bounded, a sufficient condition for (1.22) is
that energy $E$ is non trapping.  
The following stronger condition 
excludes natural potentials having a limit as $|x|\to\infty$, as shows the example 
$\rho(x)=\rho_0+\langle x\rangle ^{-\e }$. However it turns out to be convenient
from the point of vue of Definition 1.4.  
\medskip
\noindent {\bf Definition 1.5} We say that $\rho$ has the {\it defocussing condition} iff
$$G(\rho)(x,p)=\langle\nabla^2\rho(x)\cdot p,p\rangle
+{|\nabla\rho(x)|^2\over m\rho(x)}|p|^2>0, \ \forall (x,p)\in\{H_0=E\}\leqno(1.35)$$
\smallskip
In particular if $\rho$ has a  critical point, this is a non degenerate minimum. 
Under defocussing condition (1.35), if $z(s)$ is a 
special point along some bicharacteristic $\gamma$ issued from $\Lambda_0$, then for all $t>s$, $z(t)\in\gamma$ is an ordinary point.
We expect that (1.35) is related to the non-trapping condition (1.22), and provides an information
on the ``return set'' in (1.7), e.g. ${\cal R}=\emptyset$. 
It also implies that special and residual points are ``exceptional'' compared to ordinary points, 
in the same way singular points are ``exceptional'' with respect
to regular points. This allows a natural subdivision of canonical charts into ordinary and special (residual) points, so we can
speak of a regular-ordinary chart, or regular-special chart and so on. 
The main results are summarized in Proposition 4.8.
\smallskip
\noindent {\it Remark 1.1}: In case (1.25) with $m=1$ and $f_h=0$ (scattering problem), 
the asymptotic solution of $H(x,hD_x)u_h=Eu_h$ has been constructed in 
[DoMaNaTu1], Example 6, in term of Bessel functions.

\medskip
\noindent $\bullet$ {\bf $\Lambda_0$ is ``Bessel cylinder'' $(1.3)$, and $m=1$}. 
\smallskip
To the former ``$\theta$-variables'' $(t,\psi,\lambda)$,
one has now to add $\varphi$ as a parameter. Note that $\Lambda_0$ has a Lagrangian singularity at $\varphi=0$. 
This is the most technical part of the paper, and the results are only partial, because we are ignoring glancing points.
It is necessary here to assume $m=1$. 
For $m=1$, Euler identity shows that the 1-form $p\,dx-E\,dt$ vanishes on $\Lambda_+$. We are led to assume $m=1$.
For simplicity we shall also assume essentially that $H$ is of the form (1.25) with $\rho$ radially symmetric. 

Due to possible glancing points, we cannot formulate a global result in term of Maslov canonical operator as in 
Theorem 1.3. For an Hamiltonian positively homogeneous of degree 1 in the $p$ variables,
we content ourselves with computing the phase and density on $\Lambda_+^E$, so our results are most complete 
in case of the conformal metric. This simplifies further in case of a radially symmetric conformal metric. 
We refer to Sect.5 for detailed statements

\medskip
\noindent{\it 1.5 Outline of the paper}.
\smallskip
In Sect.2 we first construct eikonal coordinates on $\Lambda_+^E$ when $H_0(x,p)$ is positively homogeneous of degree $m$, 
and $\Lambda_0$ is either the vertical plane, or Bessel cylinder. Then we discuss some well-known facts about the extension of the solution of 
Cauchy problem (1.6) for large $t$. In particular we examine the case where
there are several branches of $\Lambda_+^E$ lying over
$x$, i.e. $\pi_x(X(t,\psi),P(t,\psi))=x$, leading to Van-Vleck formula. 
Following [CdV], [GuSt] we then focus to the case when $H_0$ defines a metric, according to $m=1$ (Finsler metric, or Randers symbol),
or $m>1$. 

In Sect.3 we prove Theorem 1.2. We start to recall some basic facts on Maslov theory. Then we sketch its generalisation to bi-Lagrangian distributions,
following mainly [MelUhl], where $H$ is taken microlocally to its normal form $hD_{x_n}$ on a non critical energy surface.
We make also the results of [MelUhl] more precise and adapted to asymptotics with respect to the small parameter $h$.
Thus we can construct $u_h(x)$ microlocally outside $\partial\Lambda_+^E$, and locally uniformly with respect to $h$. 
We end up by computing explicitely  
$u_h$ when $n=2$, $H=-h^2\Delta$ and $f_h$ is compactly supported, and verify that $u_h$ can be written as the sum of 2 terms, 
microlocally supported
on $\Lambda_0$ and $\Lambda_+^E$ respectively.

We start in Sect.4 to recall from [DoMaNa2], [DoNaSh] the matrix $\widetilde n\times \widetilde n$ matrix 
${\cal M}(\widetilde\phi,\widetilde\psi)$ defined on a
local chart of a Lagrangian manifold $\widetilde\Lambda$, whose determinant turns out to be the (inverse) density on $\widetilde\Lambda$.
It will be most useful in Sect.5. Then we define the phase function from which compute directly the (inverse) density $F[\Phi,dy]|_{C_\Phi}$
on $\Lambda_+^E$. Lateron we restrict to the 2-D case. In Sect.4.3, 
assuming this density is non zero, or equivalently, that $\Phi$ is a non-degenerate phase function in the sense of 
H\"ormander, we investigate some configurations of $\Lambda_+^E$ in $T^*M$ (according to $X_\psi=0$ or $X_\psi\neq0$)
and describe more closely the corresponding Lagrangian singularities (focal points) in the chart where $X_\psi\neq0$.
We relate focal points with ordinary, special or residual points.
In Sect.4.4 we complete the 1:st order asymptotics by considering the transport equations, and prove Theorem 1.3.
In Sect.4.5 we specialize further to the case of the ``conformal metric'', using also the defocussing condition (1.25), which allows 
a more complete description of focal points.

Sect.5 is the most technical and sketchy part, since we do not take glancing point into account. 
In Sect.5.1 we give necessary and sufficient for a point of the ``cylinder'' (1.3) be glancing with respect to $v_{H_0}$. In particular
we show that a glancing point at $t=0$ is also a special point. Then we describe the 
parametrization of $\partial\Lambda_+(\tau)$ provided this is a closed manifold without boundary. We compute the matrix
${\cal M}(t,\varphi,\psi,\tau)$ we introduced already in Sect.4, and show that we should take $m=1$ for its determinant 
identifies with the density on $\Lambda_+$. All computations should be carried in the extended phase-space $T^*(M\times{\bf R}_t)$. 

In Appendix we prove the density is non vanishing near focal points in the case of the ``conformal metric''.  
\medskip
\noindent{\it 1.6 Some open problems}.
\smallskip
\noindent $\bullet$ Semi-classical structure of the Green function outside a $h^\delta$-neighborhood of $\Lambda$. 

\noindent $\bullet$ Other types of initial Lagrangian manifolds, e.g. more general Bessel or Airy-Bessel beams, for which the initial manifold
$\Lambda_0$ is similar to Bessel cylinder. 

\noindent $\bullet$ Structure of the Green function near residual points, 
in particular glancing points, where Lagrangian intersection fails to be transverse, see [BoRo]. 

\noindent $\bullet$ Hyperbolic equations ($\partial\Lambda_+$ non compact)

\noindent $\bullet$ Case of multiple characteristics, involving Lagrangian manifolds with boundary $\partial\Lambda_+$
and corner $c\Lambda_+$ [MelUhl].

\noindent $\bullet$ Complex phases as in Example (7).

\noindent $\bullet$ Non linear PDE's: Melrose-Uhlmann calculus has been used in
the study of propagation of singularities for nonlinear wave equations, see [UhlZha], and it
would be natural to try the semiclassical analogue in the analysis of oscillatory solutions,
as in the work of Joly-M\'etivier-Rauch [JoMeRa,1,2] 
\medskip
\noindent {\bf Acknowledgements}: Special thanks are due to Anatoly Anikin, Serguei Dobrokhotov and Vladimir Nazaikinski for their collaboration
at an early stage of this work, especially concerning eikonal coordinates, the construction of the phase functions in Sect.4-5,
and the notion of Maslov ``bi-canonical'' operator. 
I also thank Ilya Bogaevsky for collaborating in a subsequent work on glancing intersection [BoRo], which greatly helped in editing
Sect.5, Vesselin Petkov for his valuable remarks, and a Referee for his useful comments.
Last, but no least, cheerful thanks to Johannes Sj\"ostrand and Michael Hitrik
for organizing a Workshop
on WKB theory, held on-line during the Covid time. 
This work was initially supported by the grant
PRC No 1556 CNRS-RFBR 2017-2019.
\bigskip
\noindent{\bf 2. Hamiltonians and phase functions}
\smallskip
In this Sect. we consider integral manifolds for positively homogeneous Hamiltonians on $T^*M\setminus0$, which is the first
step in constructing semi-classical Green kernels. We discuss first general facts (eikonal coordinates, Hamilton-Jacobi equation),
but more specific points will be 
discussed in Sect. 4 and 5. Basic references are [Ar], [H\"o], [GuSt].  
\smallskip
\noindent {\it 2.1 Eikonal coordinates}
\smallskip
First we recall some general facts about canonical coordinates on Lagrangian manifolds, see [DoMaNaTu], [DoNaSh].
Let $\iota:L\to T^*M$ be a smooth embedded Lagrangian manifold. We write $(x,p)=\iota(\alpha)=(X(\alpha),P(\alpha))$, where $\alpha$
are local coordinates on $L$. 
The 1-form $p\,dx$ is closed on $L$, so is locally exact, and $p\,dx=dS$ on any simply connected domain $U$ (so called canonical chart).
Such a $S$ is called an {\it eikonal} (or action) and is defined up to a constant.
If $p\,dx=dS\neq0$ on $U$, $\phi=S$ 
can thus be chosen as a coordinate on $U$, i.e. a local coordinate on $L$, and completed by coordinates $\psi\in{\bf R}^{n-1}$ such that
$$\langle P,\partial_\phi X\rangle=1, \ \langle P,\partial_\psi X\rangle=0, \ \langle \partial_\psi P,\partial_\phi X\rangle=
\langle \partial_\phi P,\partial_\psi X\rangle\leqno(2.2)$$
This holds true if $L$ is projectable on $U$, i.e. $U$ of rank $n$,
but this condition
is not necessary. Namely, consider a local chart $U$ of rank $k=1$, and let $\alpha=(\phi,\psi)$, $\phi\in{\bf R}$, $\psi\in{\bf R}^{n-1}$,
so that $L$ is defined by $x=X(\phi,\psi)$, $p=P(\phi,\psi)$, then
$$d\phi=\langle P(\phi,\psi),dX(\phi,\psi)\rangle\leqno(2.3)$$
and $\phi$ is an eikonal on $L$. See [DoNaSh], Sect. 1.4.3.

We shall deal with elliptic positively
homogeneous Hamiltonians of degree $m\geq1$ with respect to $p$ on the cotangent bundle $T^*M\setminus0$ ($M={\bf R}^n$ for simplicity),
and eventually restrict to ``conformal metrics'' to get most explicit results.  
In this Sect. we write $H$ for $H_0$.
\smallskip
\noindent {\it Examples 2.1}:

(1) $m=2$ if $H(x,p)$ is a geodesic flow associated with 
a Riemannian metric $ds^2=g_{ij}(x)\,dx^i\otimes dx^j$.
In the Riemannian case, when $E=1$, geodesics are parametrized by arc-length.  

(2) $m=1$ if $H(x,p)$ is a ``Randers symbol'', associated with a Finsler metric [Tay], [DoRo2] and reference therein. 

(3) $H$ is of the form (1.25) with $m\geq1$. Hamilton equations $(\dot x,\dot p)=v_H(x,p)$ then read
$$\dot x=\partial_pH=m|p|^{m-1}{1\over\rho(x)}{p\over|p|}, \ \dot p=-\partial_xH=|p|^m{\nabla\rho(x)\over\rho(x)^2}\leqno(2.4)$$
Our most complete results hold for such Hamiltonians with $n=2$. 
\smallskip

\smallskip
\noindent $\bullet$ {\it Case of the ``vertical plane''}.
When $H$ is positively homogeneous of degree $m$ with respect to $p$, and $\Lambda_0$ is the vertical plane, $v_H$ is always transverse to
$\Lambda_0$ when $E\neq0$, i.e. there are no glancing points.

Let $L=\Lambda_+^E$ be the flow-out of $H$ with initial data on $\Lambda_0=T_{x_0}^*M$. Thus $L$ is the union of maximally extended bicharacteristics
starting at $\Lambda_0$, and $\iota:\Lambda_+^E\to T^*M$ a Lagrangian immersion.

Let also $\psi\in{\bf R}^{n-1}$ be smooth coordinates on $\partial\Lambda_+^E$,
which we complete by $-\tau$,
the dual coordinate of $t$, so that $\partial\Lambda_+^E$ is given in $\Lambda_0$ by $\tau=0$,
and in $\Lambda_+^E$ by $t=0$. 

In the special case $H(x,p)={|p|^m\over\rho(x)}$, we have $P(\psi,\tau)=|P|_\tau\omega(\psi)$, with $\omega(\psi)\in {\bf S}^{n-1}$, 
$$|P|_\tau=(H\rho(0))^{1/m}=\bigl((E-\tau)\rho(0)\bigr)^{1/m}\leqno(2.5)$$ 
Sections are defined as follows: for small $\tau$, let $\Lambda_+(\tau)=\Lambda_+^{E-\tau}$ be the Lagrangian manifold in the energy shell $\tau+H(x,p)=E$
issued from $\Lambda_0$ at $t=0$. We consider the isotropic manifold $\partial\Lambda_+(\tau)=\Lambda_0\cap\Lambda_+(\tau)$, viewing $\Lambda_+(\tau)$
as a manifold with boundary. When $\tau=0$, we simply write $\Lambda_+(0)=\Lambda_+^E$. 

For $t\geq0$, let $\Lambda_t=\exp tv_{H}(\Lambda)$, we have the group property $\Lambda_{t+t'}=\exp tv_{H}(\Lambda_{t'})$ for all $t,t'\geq0$.

We define in a similar way the family of isotropic manifolds $\partial\Lambda_t(\tau)=\exp t\partial\Lambda_+(\tau)=\Lambda_t\cap\Lambda_+(\tau)$.

We compute the eikonal $S$ on $\Lambda_{+}(\tau)$ by integrating along a piecewise $C^1$ path
connecting the base point, say $(x_0,0)\in \Lambda_0$ (where $dx=0$) to $(x_0,P(\psi))\in \partial\Lambda_+(\tau)$, $\psi\in{\bf R}^{n-1}$, 
followed by the integral curve $x=X(t,\psi), p=P(t,\psi)$ of $v_H$ starting at $(x_0,P(\psi))$, where 
$dS=p\,dx|_{\Lambda_+(\tau)}=\langle P(t,\psi),dX(t,\psi)\rangle$. Because $\Lambda_+(\tau)$ is Lagrangian, $S$ doesn't depend 
on the choice of the base point. 
Since  $S(0,\psi)=\Const =S_0$ on $\Lambda_0$, we get
$$\eqalign{
&S(t,\psi)=S(0,\psi)+\int_{(0,\psi)}^{(t,\psi)}p\,dx|_{\Lambda_+}=S(x_0,\psi)+\int_{(0,\psi)}^{(t,\psi)}P(s,\psi,\tau)\,dX(s,\psi,\tau)\cr
&=S_0+\int_0^t\langle P(s,\psi,\tau),\dot X(s,\psi,\tau)\rangle\,ds\cr
}\leqno(2.6)$$
By Hamilton equations and Euler identity 
$\langle P(s,\psi,\tau),\dot X(s,\psi,\tau)\rangle=\langle P,\partial_pH\rangle=mH=m(E-\tau)$. 
Now $S(t,\psi)=S_0+m(E-\tau)t=S_0+mHt$ is the action on $\Lambda_+(\tau)$, and the eikonal coordinate is just $S(t,\psi)=mHt$ up to a constant $S_0$,
and we may set $S_0=0$.
Restricting to $\Lambda_+(\tau)$ we get
$dS=d\bigl(m(E-\tau)t\bigr)=m\,(E-\tau)dt$, or
$$m\,(E-\tau)dt=\langle P(t,\psi,\tau),dX(t,\psi,\tau)\rangle=\langle P(t,\psi,\tau),\dot X(t,\psi,\tau)\rangle\,dt+
\langle P(t,\psi,\tau),X_\psi(t,\psi,\tau)\rangle\,d\psi\leqno(2.7)$$
and it follows (Huygens' principle) that 
$$\langle P(t,\psi,\tau),\dot X(t,\psi,\tau)\rangle=mH, \quad \langle P(t,\psi,\tau),\partial_\psi X(t,\psi,\tau)\rangle=0
\leqno(2.8)$$ 
We will denote for short
$(X,P)=\bigl(X(t,\psi,\tau),P(t,\psi,\tau)\bigr)$ or $X(t,\psi),P(t,\psi)$ when $\tau=0$. This is called the {\it leading front}. 
\smallskip
\noindent $\bullet$ {\it Case of ``Bessel cylinder''}. Even for positively homogeneous $H$ of degree $m$ with respect to $p$, there may be
glancing points, but $\Lambda_+^E$ is still an immersed Lagrangian manifold away from the trajectories starting at
the glancing points. In this work we restrict to local charts where $\Lambda_+^E$ is immersed, glancing intersection being considered in [BoRo]. 

Let $L=\Lambda_0$, then $p\,dx|_{\Lambda_0}=d\varphi$, which shows that $\varphi$ is the eikonal on $\Lambda_0$. 

Let $H$ be homogeneous of degree $m=1$ with respect to $p$, and 
$L=\Lambda_+$ be the flow out of $\Lambda_0$ by $v_H$ in the extended phase-space.
According to Euler identity, $\dot S=p\dot x-H=p\partial_pH-H=0$ along the trajectories,
the eikonal on $\Lambda_+$ is again $S=\varphi$, and  $(\varphi,\psi,t)=\varphi$ are eikonal coordinates. 
This results also from the fact that $v_H$ preserves $p\,dx$ on $\Lambda_+$ (see also (2.29) below).

Let $L=\Lambda_+(\tau)$ be the flow-out of $H$ at energy $E-\tau$ with initial data on (1.3).
Along $\Lambda_+(\tau)$ we have as in (2.6)
$$\eqalign{
S(t,\varphi,\psi,\tau)&=S(0,\varphi,\psi,\tau)+\int_{(0,\varphi,\psi)}^{(t,\varphi,\psi)}p\,dx|_{\Lambda_+}=
S(0,\varphi,\psi,\tau)+\int_{(0,\varphi,\psi)}^{(t,\varphi,\psi)}P(s,\varphi,\psi,\tau)\,dX(s,\varphi,\psi,\tau)=\cr
&=S_0+\varphi+\int_0^t\langle P(s,\varphi,\psi,\tau),\dot X(s,\varphi,\psi,\tau)\rangle\,ds
}$$
As before, $\langle P,\dot X\rangle=mH$, so on $\Lambda_+(\tau)$
$$S(t,\varphi,\psi,\tau)=\varphi+m(E-\tau)t+S_0\leqno(2.10)$$
and we may set again $S_0=0$. This is the eikonal on $\Lambda_+(\tau)$.
Identifying the differential of $S$ we get (omitting again variables)
$$\langle P,\dot X\rangle=mH,\ \langle P,\partial_\psi X\rangle=0, \ 
\langle P,\partial_\varphi X\rangle=1
\leqno(2.11)$$
In case of (1.25) with $m=1$ and radially symmetric $\rho$, with $\nabla\rho(x)\neq0$ if $x\neq0$. Let $E>0$ be such that
$\rho(0)\neq1/E$. Then we have $\varphi=\varphi(\tau)\neq0$ on $\Lambda_+(\tau)$ for small $\tau$, 
and $(t,\psi)$ parametrize $\Lambda_+(\tau)$, at least for small $t$. See Example 5.2 below. 

In Sect.5 we develop a slightly different point of vue, extending the phase-space to $T^*(M\times{\bf R})$, which amount
to introduce variable $\tau$ as the dual  variable of $t$, and change 
$\langle P(t,\psi,\tau),dX(t,\psi,\tau)\rangle$ to $\langle P(t,\psi,\tau),dX(t,\psi,\tau)\rangle+\tau\,dt$ accordingly.
But for simplicity we restrict to radially symmetric $\rho$ as above. The general case is complicated by the occurence of glancing points,
i.e. where $v_H$ fails to be transverse to the initial manifold $\Lambda_0$. Nevertheless, (2.10) and (2.11)
are sufficiently general to be considered as a ``starting point''
in order to investigate the case of glancing intersection. 

\medskip
\noindent {\it 2.2 Hamilton-Jacobi equation for small $t$ and phase functions}
\smallskip
Because of focal (or glancing) points we cannot in general find a phase-function $\Psi^E(x)$ such that $H(x,\partial_x\Psi^E(x))=E$, 
so we obtain it as a critical value
of a phase $\Psi(x,t)$. To do this we
solve Hamilton-Jacobi equation (HJ) in the extended phase space $T^*(M\times{\bf R})$, which is 
the suitable framework to vary $(t,\tau)$ as well ($\tau$ being set eventually to 0). So we look for a phase function 
$\Psi(x,t)$ satisfying
$$\partial_t\Psi+H(x,\partial_x\Psi)=E, \  \Psi|_{t=0}=\phi\leqno(2.12)$$
with given $\phi$ (to be chosen lateron), and prescribed $\partial_t\Psi(x_0,0)=\tau_0,\partial_x\Psi(x_0,0)=\eta_0$
satifying $\tau_0+H(x_0,\eta_0)=E$. By Hamilton Eq., $\tau=\tau_0$ is a constant of the motion. 
It is well-known (Hamilton-Jacobi theory, see [H\"o,Thm6.4.5]) that (2.12) as a unique solution for small $t$. 
This is the generating  function of the Lagrangian manifold the extended phase-space
$$\Lambda_+=\{p=\partial_x\Psi(x,t), \ \tau=\partial_t\Psi(x,t), \ x,t\in M\times{\bf R}_+ \}\subset T^*(M\times{\bf R})\leqno(2.13)$$ 
constructed along the integral curves of $v_H$ starting at $t=0$ from the Lagrangian manifold
$\Lambda_\phi$ in $T^*M$ given by $p=\partial_x\phi$. 
Its section at fixed $t,\tau$ is the Lagrangian manifold 
$$\Lambda_{\Psi,t,\tau}=\{p=\partial_x\Psi(x,t), \ x\in M\times{\bf R} \}\subset T^*M$$
which is simply the flow out $\Lambda_t(\tau)$ of $\{p=\phi'(x)\}$ in $H(x,p)=E-\tau$ at time $t$. 
We choose the initial condition to be the standard pseudo-differential phase function
of the form $\phi(x)=x\eta$. Here $\eta$ is a parameter, we choose so that $\phi'|_\Lambda=P(\widetilde\psi)$
where $\widetilde\psi$ are coordinates on $\Lambda_0$, that could be taken of the form $\widetilde\psi=(\psi,\tau)$, $\psi$ being
coordinates on $\partial\Lambda_+(\tau)$.

Actually Sect.4.2 or 5.2 will not directly rely on Hamilton-Jacobi theory, but rather on the
construction of a generating family (in the sense of Arnold) or a non degenerate phase-function (in the sense of H\"ormander)
$\Phi(x,t,\theta)$ for $\Lambda_+$. This is the 1-jet of $\Psi$ along $\Lambda_+$. On the other hand, since $\Phi(x,t,\theta)$
is smooth in a neighborhood of $\Lambda_+$, it can be differentiated and thus provides the jet of $\Psi$ at infinite order along $\Lambda_+$. 
\smallskip
The phase $\Psi(x,t)$ has the property (for a general Hamiltonian) that along each of these curves
$$\Psi(x(t),t)=\phi(x)+\int_0^t\bigl[\langle{\partial H\over\partial p}(x(s),p(s)),p(s)\rangle+\tau(s))\bigr]\,ds\leqno(2.17)$$
which is $\phi(x)+(mH+\tau)t$ for $H$ positively homogeneous of degree $m$, and simplifies to $\phi(x)+Et$ when $m=1$, see Sect.5. This is also the action 
$\int _{(0,x_0)}^{(t,x)}L(q(s),\dot q(s))\,ds=\int _{(0,x_0)}^{(t,x)}\langle p,dq\rangle-H\,dt$, where the integral is computed 
along an integral curve of $v_H$ from $x_0$ to $x$, see e.g. [Ar2,Sect.46]. 
\medskip
\noindent {\it 2.3 The phase functions ``in the large''  and the semi-classical Cauchy problem}
\smallskip
We discuss the case of the vertical plane $\Lambda_0=T^*_{x_0}M$, which reduces to standard
variational problems in the space variable. 

Let $\theta$ parametrize the initial condition $\phi$ in (2.12). The relevant case is $\phi=x\eta$, we can take $\theta=(\psi,\lambda)$
as local coordinates on $\Lambda_0$ near $\partial\Lambda_+^E$.
So for $\tau=0$ the initial surface, 
$\partial\Lambda_+^E$ is compact and of the form
$(x,p)=(X(\psi)=0, P(\psi)=\eta))$, $\psi\in {\bf S}^{n-1}$.
We will add $t$ to the ``$\theta$''-variables since we
require $\tau=\partial_t\Psi=0$.

We consider the map $(t,\psi)\mapsto x=\pi_x\exp tv_H(x_0,P(\psi))$, or which is the same, $(t,\eta)\mapsto x=\pi_x\exp tv_H(x_0,\eta))=\Exp_{x_0}^t\eta$. 

So far we have described the phase function when ``moving along'' $\Lambda_+$ for small $t$. 
Thus the critical point of $(t,\theta)\mapsto\Psi(x,t,\theta)$, 
is such that $x=X(t,\theta)$.   
A ``dual'' point of vue is to fix $x$ and find the set of $(t,\theta)$ with $x=X(t,\theta)$.
In the Riemannian case ($m=2$) this is related to the problem of geodesic completeness, which holds locally. 
Namely if $|x|$ is small enough,
there is a unique $(t,\eta)$ such that $x=\Exp_{x_0}^t\eta$. This holds globally if the Riemaniann manifold ($m=2$)
is geodesically convex.
Otherwise, the ``inverse map'' $x\mapsto(t,\eta)$ may be multivalued. 

It is well known [CdV, pdf p.132] that the global geodesic convexity can be relaxed (locally) to a non-degeneracy condition. Namely, let
$H$ be associated with a Lagrangian ${\cal L}(x,\dot x)$ strictly convex with respect to $\dot x$,
in particular if $H$ is positively
homogeneous of degree $m>1$ with respect to $p$.
 
Let ${\cal K}\subset{\bf R}^n_p$ be a compact set, which will be identified with the support of $a(p;h)$ in (1.2) (after eliminating $x$).
For fixed $(x,t)$ we make the generic assumption~:
\smallskip
\noindent (2.18) {\it  For all $\eta\in {\cal K}$ such that $x=\Exp_{x_0}^t\eta$, the map ${\bf R}^n\to{\bf R}^n$, 
$\xi\mapsto \Exp^t_{x_0}\xi$ is a local diffeomorphism near $\eta$~: in other terms,
$x_0$ and $x$ are not conjugated along any trajectory that links them together within time $t$,
with initial momentum $\xi$. }
\smallskip
The set of such $(x,t)$ is an open set $\Omega_{x,t}\subset{\bf R}^{n+1}$, and by Sard's Theorem, its complement has Lebesgue measure 0. 

Fixing $(x,t)$, (2.18) implies by Morse theory that
$\eta\mapsto x=\Exp_{x_0}^t\eta$ has a discrete set of pre-images $\eta\in{\cal K}$. 

Fixing $x$, consider now the pre-images of $(t,\eta)\mapsto x=\Exp_{x_0}^t\eta$. It can happen that the integral manifold of $v_H$
has several sheets over $x$, so several values of $t$ can contribute to the same $x=X(t,\psi)$. However, under the non-trapping condition 
$|X(t,\psi)|\to\infty$ as $t\to\infty$, there is again, generically, a finite number of such $t_j$. 
Namely, it suffices that (2.18) holds with a time $T$ such that for $t\geq T$, $X(t,\psi)$ will never coincide again with $x$.
Moreover these $t_j$ are 
non-degenerate critical points of $t\mapsto\Phi(x,t,\eta)$.  
This holds in particular when $x_0$ and $x$ are connected by (possibly several) minimal geodesics for the Riemannian
metric associated to $H$, each indexed by some $\eta_\alpha$. 
\smallskip
At least for small $t$ we seek for a solution of (1.6)  of the form
$$v_h(t,x)=\int ^* e^{i\Psi(t,x,p)/h} b(t,x,p;h)\,dp \leqno(2.20)$$ 
where the phase function $\Psi$ 
as in (2.12) with initial condition $\Psi|_{t=0}=xp$~; the symbol $b(t,x,p;h)$ verifies
$b(0,x,p;h)=a(x,p;h)$, and solves some 
transport equations along the integral curves of $v_{H}$.
One may address the problem of a semi-classical ``close form'' of (2.20), i.e. of performing the integration with respect to $p$,
so that the final expression is of WKB type. Under Assumption (2.18) the answer to this problem is given by
Van Vleck formula [CdV, pdf p.132] which gives $v_h$
as a finite sum
$$\Sum_\alpha {A(\eta_\alpha)\over\sqrt{\Jac _{x_0}(\eta_\alpha)}}e^{i\Psi(t,x,\eta_\alpha)/h}e^{-i\pi\ind (\gamma_\alpha)/2}\leqno(2.21)$$
(at leading order in $h$) indexed by all $\eta_\alpha\in\supp A$ such that  $\Exp_{x_0}^t\eta=x$. Here $A$ is the principal
part of $a$ (where we have eliminated the dependence in $x$), $\Jac _{x_0}(\eta_\alpha)$ is the Jacobian of 
$\Exp_{x_0}^t$ at $\eta_\alpha$ and $\ind (\gamma_\alpha)$ Morse index of the integral curve $s\mapsto\exp sv_{H}(x_0,\eta_\alpha)$, $s\in[0,t]$.
In other words, under Assumption (2.18) it suffices to use only non singular charts on $\Lambda_+$ over $x$, and the solution is expressed
in term of finitely many oscillating functions.

When Assumption (2.18) is not met, i.e. $x=x_*$ is conjugated to $x_0$, then there is at least one focal point $(x_*,p_*)$
over $x_*$ in $\Lambda_+$.
The construction of the canonical operator (see Sect.3) necessarily uses a singular chart in a neighborhood of 
$(x_*,p_*)$, and the solution in a neighborhood of $x_*$ involves not only simple oscillating functions corresponding to nonsingular 
charts as in (2.21) (if any) but also an  integral of an oscillating function over some of the momenta (or ``$\theta$-variables''). 
The total number of singular and nonsingular charts over $x_*$ however remains finite, and so (generically) only a finite sum 
of integrals and simple oscillating functions contribute (one summand per each chart).
So Maslov canonical operator encompasses Van Vleck formula. 
See Sect.3 for more details.
\smallskip
\noindent {\it 2.4 Distances and generating functions}
\smallskip
When $\Lambda_0$ is the vertical plane, the phase function $\Phi$ is related to the ``distance'' to $x_0$ for the
``metric'' implied by $H_0$, which is of special interest.  We make here some general remarks, mainly following [CdV], [GuSt].

In Sect.4-5 we shall discuss how to parametrize, by a non degenerate phase function, 
the flow of $v_H$ out of some Lagrangian plane, when $H$ is positively homegeneous of degree $m$.
It includes the case $m=1$ which plays an important role because of Finsler metrics. 
So we begin with a general discussion on
corresponding symplectic maps. 

Let $H$ be a positively
homogeneous Hamiltonian of degree $m$ with respect to $p$, defined on $T^*M\setminus0$, 
and $\Gamma\subset T^*M\setminus0\times T^*M\setminus0$ be the graph of $\exp v_H$ (time-1 flow). 
Recall from [GuSt,Chap.5, formula (5.6) and Thm 5.4.1] that
$$(\exp v_H)^* (p\,dx)-p\,dx=(m-1) dH\leqno(2.29)$$
Integrating over a path $\gamma_E\subset\{H=E\}$, we recover the fact that 
$\int_{\gamma_E}(\exp v_H)^* (p\,dx)=\int_{\gamma_E}p\,dx$. 

So when $m=1$, not only the 2-form, but also the 1-form $p\,dx$ are preserved under the flow of $v_H$. 
In this case, $v_H$ is actually the lift of a vector field on $M$. We are forced to assume $m=1$ when $\Lambda_0$ is Bessel cylinder, see Sect.5.

When $m>1$,  formula (2.29) gives a generating function for $\Gamma$
under the following assumption (which however might never be verified when $m\neq2$)~: 
\smallskip
\noindent (2.30) {\it Let $\pi_{M\times M}:T^*(M\times M)\to M\times M$ be the natural projection, and assume 
$\pi_{M\times M}:\Gamma\to M\times M$ is a
diffeomorphism, i.e. for all $(x,y)\in M\times M$, there is a unique $\xi\in T_x^*M$ such that $y=\pi_M\exp v_H(x,\xi)$. }
\smallskip
Then we say that $\Gamma$ is {\it horizontal}. In case of a geodesic flow ($m=2$) (2.30) holds true when $M$ is geodesically convex. 
Provided (2.30), $\Gamma$ has a generating function $\chi$, i.e. $d\chi=\pr _2^*(p\,dx)-\pr _1^*(p\,dx)$, where 
$\pr _i:T^*(M\times M)=T^*M\times T^*M\to T^*M$ is the projection onto the $i$:th factor, and
$\pr _1\circ(\pi_{M\times M}|_\Gamma)^{-1}$ is a diffeomorphism $M\times M\to T^*M$. 
Moreover we can then represent $\chi$ as
$$\chi=\bigl(\pr _1\circ(\pi_{M\times M}|_\Gamma)^{-1}\bigr)^*(m-1)H\leqno(2.31)$$
In case of the geodesic flow ($m=2$)
$\chi(x,y)={1\over2}\dist (x,y)^2$. 
Formula (2.29) is related to {\it exact symplectic twist maps} as follows. An exact symplectic twist map [Ar], [Ka], [GuSt], [CdV7]
$F:T^*M\to T^*M$ is a symplectic map with a generating function $S_1:M\times M\to {\bf R}$, $(x,X)\mapsto S_1(x,X)$ which satisfies 
$$F^*(p\,dx)-p\,dx=P\,dX-p\,dx=dS_1(x,X)\leqno(2.32)$$
$(p,x)$ and $(P,X)$ are related by $p=-\partial_xS_1, P=\partial_X S_1$. In notation $S_1$ the subscript 1 refers to time-1 flow.
In case $H(x,p)=p^2$ (flat metric on ${\bf R}^n$), comparing (2.32) with (2.29), i.e. $dS_1=dH$, we get 
$S_1(x,X)={1\over4}(x-X)^2$, and more generally, if $H(x,p)=|p|^m$,
with $m>1$, $S_1(x,X)=\bigl({m-1\over m}\bigr)^{m/(m-1)}|x-X|^{m/(m-1)}$.

Again, $S_1$ is not well defined when $m=1$. More generally $F(x,y)$ coincides with
$\chi(x,y)$ above for the geodesic flow. 

For Hamilton-Jacobi (HJ) equation we have the following Proposition, extending (2.15) for large $t$. 
Assume $H$ is associated with a Lagrangian convex with respect to $\dot x$. 
Let $x_0,y_0\in M$ be non conjugate points 
along an extremal curve $\gamma_0(t)$ such that $x_0=\gamma_0(0)$ and 
$y_0=\gamma(t_0)$, and $(x_0,\xi_0),(y_0,\eta_0)$
the corresponding points in $T^*M$. 
\medskip
\noindent {\bf Proposition 2.2} [CdV,Thm 14, pdf p.45]: {\it Let $(t_0,x_0,y_0)$ be as above. 
Then for any $(x,y)$ close to $(x_0,y_0)$, and $t$ close to $t_0$, there is a unique extremal curve $\gamma$  
such that $x=\gamma(0)$ and $y=\gamma(t)$. Let $\widetilde S(t,x,y)$ be the action along these curves (minimizing the Lagrangian action)
This is a generating function for the Hamiltonian flow near $(x_0,\xi_0)$, verifying HJ equation}
$$\partial_t\widetilde S+H(y,\partial_y\widetilde S)=0\leqno(2.33)$$
\smallskip
This is verified 
in the Riemannian case
$\widetilde S(t,x,y)={F(x,y)\over2(1+t)}={\dist^2(x,y)\over2(1+t)}$ where $F$ is the exact symplectic twist map considered above, 
and can be identified with the phase in the Heat kernel. We can check (2.33) trivially when $H={1\over2}p^2$, $\widetilde S(0,x,y)={1\over2}(x-y)^2$,
see [Ar,p.255].
This holds also under assumption (2.30).
Clearly under Hypothesis (2.18), (2.33) extends (2.12) for large times.

So far we have assumed some convexity of $H$ with respect to $p$. The case $m=1$ (Finsler metric and Randers symbols)
is investigated in [Ta]~: it turns out that similar results hold when the {\it square} of Finsler metric or Randers symbol enjoys some 
convexity property, so for a ``conformal metric'' the case $m=1$ makes no difference. In Sect.5 we shall require $m=1$, but $\Phi$ is no 
longer associated with a distance on $M$. 
\bigskip
\noindent{\bf 3. Maslov canonical operators and bi-Lagrangian distributions}
\smallskip
Here we prove Theorem 1.2. Our purpose is to describe the solution globally, including unfolding of Lagrangian singularities.
Among many references to the subject we make use in particular of 
[M], [H\"o], [Du], [Iv], [BaWe], [CdV], [GuSt], [GrSj], [DoZh], [DoNaSh], [DoRo]. 

First we recall the asymptotic stationary phase formula
for a quadratic phase function [H\"o,Lemma 7.7.3]. Let $A$ be a symmetric non-degenerate 
matrix, then
$$\int e^{i\langle Ax,x\rangle/2h}u(x)\,dx=\bigl(\det\bigl(A/(2i\pi h\bigr)\bigr)^{-1/2}\Sum_0^{k-1}(h/(2i))^j\langle A^{-1}D,D\rangle^ju(0)/j!+{\cal O}(h^k)\leqno(3.1)$$
Since we shall ignore for simplicity Maslov indices, this formula has the advantage of hiding phase factors like $e^{-i\pi n/4}$,
which we could restore by choosing an appropriate branch of the square root in the complex plane. A similar formula [H\"o,Theorem 7.7.5] holds
for $\langle Ax,x\rangle/2$ replaced by $f$ with a non-degenerate critical point at $x_0$ and Hessian matrix $A$.

\medskip
\noindent {\it 3.1. Lagrange immersions and non-degenerate phase functions}:
\smallskip
A smooth function $\Phi:(x,\theta)\mapsto\Phi(x,\theta)$, $\theta\in{\bf R}^N$, 
defined near $(x_0,\theta_0)$ with $\xi_0=\partial_x\Phi(x_0,\theta_0)$ 
is called a non-degenerate phase function in the sense of H\"ormander iff the $(n+N)\times N$ matrix 
$(\Phi''_{\theta x},\Phi''_{\theta\theta})$ has rank $N$ on the critical set 
$$C_\Phi=\{(x,\theta)\in M\times{\bf R}^N: {\partial \Phi\over\partial\theta}(x,\theta)=0\}\leqno(3.3)$$  
Then 
$$\iota_\Phi:C_\Phi\to\Lambda_\Phi=\{(x,\Phi'_x(x,\theta): (x,\theta)\in C_\Phi\}\leqno(3.4)$$ 
is a local Lagrangian embedding (diffeomorphism). 

It is easy to prove [Iv,(1.2.7)] that
$$N-\rank\Phi''_{\theta\theta}=n-\rank d\pi_\Lambda(\iota_\Phi(x,\theta))\leqno(3.5)$$

Let also $\pi_{\Lambda_\Phi}:\Lambda_\Phi\to M$ (or simply $\pi_x$) be the natural projection. 
If $k=\rank d\pi_\Lambda(\iota_\Phi(x_0,\theta_0))$, we say that $\Lambda_\Phi$ has rank $k$ in a neighbhd $U$ of
$(x_0,\xi_0)$, and call $U$ a {\it local chart} of rank $(\geq)k$ near $(x_0,\xi_0)$. 
If $k=n$, $U$ is called a ``regular'' chart, and
$\Lambda_\Phi$ is called ``projectable'' or ``horizontal'' on $U$. On the other extreme, if $k=0$, $U$ is called a 
``maximally singular'' chart, and $\Lambda_\Phi$ is called ``vertical'' on $U$. 

If at some $z=(x,\xi)\in T^*M$, $T_z\Lambda_\Phi$ is transverse to the vertical plane $V_z=\{(0,\delta\xi)\}$ (i.e. $z$ is a regular
point) then (3.5) shows that
$\Phi''_{\theta\theta}$ is of maximal rank $N$. 

When $k<n$ we start to add some extra variables: namely there exists a partition of variables $x=(x',x'')$ such that the $(N+n-k)\times(N+n-k)$ matrix 
$$\Hess_{(x'',\theta)}(\Phi)=\pmatrix{\Phi''_{x''x''}&\Phi''_{x''\theta}\cr\Phi''_{\theta x''}&\Phi''_{\theta\theta}}\leqno(3.6)$$ 
is non degenerate. So the map $(x'',\theta)\mapsto\Phi(x,\theta)-x''\xi''$ has a non-degenerate critical point 
$\theta_c=\theta(x',\xi''), x''_c=x''(x',\xi'')$ with the critical value 
$S(x',\xi'')=\Phi(x',x''_c,\theta)-x'\xi''_c$. The projection 
$\widetilde\pi:\Lambda_\Phi\to T^*_x{\bf R}^n, (x,\xi)\mapsto\widetilde\pi(x,\xi)=(x',\xi'')$ becomes of maximal rank $n$. Hence $\Lambda_\Phi$ near $x$ 
is parametrized by $S(x',\xi'')$. 
\medskip
\noindent {\it Remark 3.1}: The above non-degeneracy condition on $\Phi$ is equivalent to (non-degeneracy
in the sense of H\"ormander) $d_{x,\theta}\Phi(x_0,\theta_0)\neq0$, 
and $d_{(x,\theta)}\partial_{\theta_1}\Phi,\cdots,d_{(x,\theta)}\partial_{\theta_N}\Phi$ are linearly independent on 
the critical set $C_\Phi$. 
The property stated above means that, if $\Phi$ is non-degenerate
in the sense of H\"ormander, then it is always possible to find coordinates such that $\pi_\xi:(x,\xi)\mapsto\xi$ has rank $n$.
Actually, there are coordinates such that 
$(x,\theta)\mapsto\Phi(x,\theta)-x\xi$ has a non-degenerate critical point, so that 
$$\Hess_{x,\theta}(\Phi)=\pmatrix{\Phi''_{xx}&\Phi''_{x\theta}\cr \Phi''_{\theta x}&\Phi''_{\theta\theta}}\leqno(3.7)$$ 
is non degenerate and $\Lambda_\Phi$ is of the form
$\Lambda_\Phi=\{(-\phi'(\xi),\xi)\}$, see [H\"o,Proposition 25.1.5].
This follows from the fact that while $\pi_x:(x,\xi)\mapsto x$ is invariantly defined under diffeomorphisms
in $M$, this is not the case for the horizontal projection
$\pi_\xi:(x,\xi)\mapsto \xi$. 
For the generating function $\Phi$ of $\Lambda_+^E$ constructed in Sect.4 below, 
$\Hess_{x,\theta}(\Phi)$ is actually degenerate in the ``natural'' coordinates of the problem, while $\Hess_{x'',\theta}(\Phi)$ is not, see Remark 4.2.
\medskip
Let us recall the expression for the (inverse) density on $C_\Phi$. Let $y=(y_1,\cdots,y_n)$ be some local coordinates on $C_\Phi$ and $|dy|$
corresponding Lebesgue measure. Then the non vanishing, real function 
$$F[\Phi,|dy|]={dy\wedge\,d\Phi'_\theta\over dx\wedge\,d\theta}={dy\wedge\,d\Phi'_{\theta_1}\wedge \cdots\wedge\,d\Phi'_{\theta_N}\over
dx\wedge\,d\theta_1\wedge\cdots\wedge\,d\theta_N}\leqno(3.8)$$
is well-defined near $C_\Phi$ as the quotient of two volume forms, see [H\"oIV,p.14], [NaStSh], [DoNaSh, (2.8)], [AnDoNaRo3].
Restricted to $C_\Phi$, its absolute value $\mu_\Phi$ defines the (inverse) density on $C_\Phi$. Computed on the complexified tangent space to $C_\Phi$,
the variations of the argument of $F[\Phi]$ can define also the variations of Maslov index (see [DoNaSh], [AnDoNaRo3]),  
which we shall ignore in this paper.
We can also write the absolute value of (3.8) on $C_\Phi$ in the form, see [GrSj, Sect.11]
$$\mu_\Phi=|F[\Phi,|dy|]|_{C_\Phi}=\Big|\det\pmatrix{{\partial y\over\partial x}&{\partial y\over\partial\theta}\cr
{\partial^2\Phi\over\partial x\partial\theta}&{\partial^2 \Phi\over\partial\theta^2}}\Big|_{C_\Phi}$$
It is actually independent of the choice of coordinates on $C_\Phi$ but it does depend on the choice of local coordinates $x$. 
\medskip
\noindent {\it 3.2 Maslov canonical operator acting on Lagrangian distributions}.
\smallskip
Let $u_h$ be a semi-classical Lagrangian distribution (or oscillatory integral) i.e. locally
$$u_h(x)=\int ^* e^{i\Phi(x,\theta)/h}b(x,\theta;h)\,d\theta\leqno(3.9)$$
where $\Phi(x,\theta)$ is a non-degenerate phase function in the sense above, and 
$b(x,\theta;h)=b_0(x,\theta)+hb_1(x,\theta)+\cdots$ an amplitude.
With $u_h$ we associate as in (3.4) the Lagrangian submanifold
$\Lambda_\Phi$. 

It is proved in [H\"o,Proposition 25.1.5] that, using that (3.7) is non degenerate, we can choose local coordinates $\xi\in{\bf R}^n$ on $\Lambda_\Phi$,
take $h$-Fourier transform
${\cal F}_{h,x\to\xi}u_h(\xi)=(2\pi h)^{-n/2}\int e^{-ix\xi/h}u(x)\,dx$ and expand by stationary phase.
The half-density in the local chart $(C_\Phi,\iota_\Phi)$ is the given by
$\sqrt{d\mu_\Phi}=|\det \Phi''|^{-1/2}|d\xi|^{1/2}$ (denoting $\Phi''=\Hess_{x,\theta}(\Phi)$ for short),
and the (oscillating) principal symbol of $u_h$ in $\Lambda_\Phi$ by
$$e^{i\phi(\xi)/h}A_0(\xi)=e^{i\phi(\xi)/h}e^{i\pi\sgn\Phi''/4}b_0(x(\xi),\theta(\xi))\sqrt{d\mu_\Phi}\leqno(3.11)$$
Here $\phi$ is the ``reduced phase function'' such that $\Lambda_\Phi=\{(-\phi'(\xi),\xi)\}$.
We emphasize that the factor $e^{i\phi(\xi)/h}$ in (3.11) is due to the fact that the Lagrangian manifold $\Lambda_\Phi$ is not necessarily conic
in $\xi$ (see e.g. [Dui1,p.222]). 

Alternatively when (3.6) is non degenerate, one can express the (oscillating) principal symbol of $u_h$ 
taking partial Fourier transform
$$v_h(x',\xi'';h)={\cal F}_{h,x''\to\xi''}u_h(x',\xi'')=(2\pi h)^{(k-n)/2}\int e^{-ix''\xi''/h}\,dx''\int e^{i\Phi(x',x'',\theta)/h}b(x',x'',\theta;h)\,d\theta$$
leading again to an expression of WKB type as in (3.11), and locally
$$\Lambda_\Phi=\Lambda_S=\{(x,\xi): x''=-\partial_{\xi''}S(x',\xi''), \xi'=\partial_{x'}S(x',\xi'')\}$$
Thus we obtained a {\it reduced phase functions}, with least possible number of variables $\theta$, i.e. at most $n$. 
When $k=n$, then $u_h$ assumes simply a WKB form in $x$ variables.

Conversely, let $\iota:L\to T^*M$ be a Lagrangian immersion. We know ([H\"o, Theorem 21.2.16]) that 
there exists a covering of $L$ by canonical charts $U$, such that $L$
is parametrized in each $U$ by a non-degenerate phase function $\Phi$. 
The Lagrangian immersions $\iota$ and (3.4) 
have the same image on $U$ and $C_\Phi$ is a submanifold
of dimension $n$. In particular, $\iota_\Phi:C_\Phi\to L$ is a diffeomorphism onto its image.
These phases can be chosen coherently, and define a class of ``reduced phase functions'' $\phi$, parametrizing $\iota$ locally.
This gives the fibre bundle of phases ${\cal L}_h$, 
including Maslov indices, equipped with transition functions.
We are also given local
smooth half-densities $|d\mu_\Phi|^{1/2}$ on $L$, defining the fibre bundle of half-densities $\Omega^{1/2}$,
equipped with transition functions. 
The collection of these objects make a fibre bundle  $\Omega^{1/2}\otimes{\cal L}_h$  over $L$. 
A section of $\Omega^{1/2}\otimes {\cal L}_h$ will be written as 
$\bigl[K_{(L,\mu)}^h a\bigr](x;h)$, 
where $K_{(L,\mu)}^h$ is called Maslov canonical operator. At leading order $\bigl[K_{(L,\mu)}^h a\bigr]$
reduces to its oscillating symbol (3.11). 
The set of such Lagrangian distributions microlocally supported on $L$
will be denoted by $I(M;L)$.
\smallskip
We apply Maslov canonical operator for constructing
solutions to homogeneous equation 
$$(H(x,hD_x;h)-E)u_h=0$$ 
microlocally supported on $L=\Lambda_+^E$ in the characteristic foliation of $\Sigma_E$.
Here $H(x,hD_x;h)$ is a $h$-PDO with principal symbol $H(x,p)$, and we assume for short $E=0$, so that we denote $\Lambda_+^E$ by $\Lambda_+^0$.
If $u_h$ is a Lagrangian distribution
locally of the form (3.9), then the same holds for $Hu_h$.

The phase $\Psi$ (with time $t$ in Hamilton equations as one of the $\theta$-parameters)
is determined by HJ equation (2.12), with initial data on $t=0$, which gives (locally) the Lagrangian embedding
(3.4) with image $\Lambda_+$. Alternatively we can use the 1-jet $\Phi$ of $\Psi$ along $\Lambda_+^0$ as above. 
In particular, $C_\Phi\subset\{\partial_t\Phi=0\}$.
We prescribe the amplitude such that $b|_{t=0}=a$. The construction of $\bigl[K_{(\Lambda_+^0,\mu)}^h a\bigr]$ goes as follows.
The amplitude of $H(x,hD_x;h)u_h(x)$ has leading term 
$$H_0\bigl(x,\partial_x\Phi(x,\theta)\bigr)b_0(x,\theta)$$ 
Moreover if $H(x,p;h)$ has sub-principal symbol $H_1(x,p)$, and $H_0(x,p)=0$ on $\Lambda_+^0$,
$H(x,hD_x;h)u_h(x)$ has principal symbol 
$$h\bigl({1\over i}{\cal L}_{v_{H_0}}|_{\Lambda_+}+H_1\bigr)b_0|d\mu_\Phi|^{1/2}\leqno(3.12)$$
where ${\cal L}_{v_{H_0}}$ denotes Lie derivative along $v_{H_0}$ acting on half-densities $u|dy|^{1/2}$ as
$${\cal L}_v(u|dy|^{1/2})=\bigl(v_j(y){\partial u\over\partial y_j}+{1\over2}(\divr v(y))u(y)\bigr)|dy|^{1/2}\leqno(3.13)$$
in local coordinates. 
For Schr\"odinger operator $H(x,hD_x)=-h^2\Delta+V(x)$, $\Lambda_+$ is ``horizontal'', and
${\cal L}_{v_{H_0}}|_{\Lambda_+}\bigl(b_0|d\mu_\Phi|^{1/2}\bigr)$ takes the form 
$\bigl(\Sum_j\Phi'_j(x){\partial b_0\over\partial x_j}+{1\over2}\Delta\Phi(x)b_0(x)\bigr)|dx|^{1/2}$,
and a similar expression when $\Lambda_+$ is ``vertical'', see e.g. [DoRo,(b.14)]. 
Then (3.9) solves $Hu_h=0$ mod ${\cal O}(h^2)$. 

Note that on $\Lambda_+$, ${\cal L}_{v_{H_0}}=d/dt$.
Provided $\Phi$ is a non-degenerate phase-function, (3.12) admits a global solution, computed on each canonical chart.
For instance on a regular chart, this is just WKB construction. In a totally singular chart instead, we solve (3.12) in Fourier representation, 
and more generally in the mixed representation.
 
The function $b_0(x(\xi),\theta(\xi))$ is smooth in $\xi$, but of course when expressed in $x$-variable, singularities may occur
do to singular Jacobians at focal points. 

The fact that $u_h$ solves $Hu_h=0$ mod ${\cal O}(h^2)$ is also expressed by the {\it commutation relation} 
$$H(x,hD_x;h)\bigl[K_{(\Lambda_+,\mu)}^h a\bigr](x;h)=\bigl[K_{(\Lambda_+,\mu)}^h h\bigl({db_0\over dt}+iH_1b_0\bigr)\bigr](x;h)
+{\cal O}(h^2)={\cal O}(h^2)\leqno(3.14)$$
\medskip
\noindent {\it 3.3  Bi-Lagrangian distributions}.
\smallskip

As in [MelUhl] our constructions make use of symbolic calculus adapted to Lagrangian intersection. So we need first to translate 
some notions
relative to asymptotics with respect to smoothness (or ``standard Pseudo-differential Calculus''), to the framework of asymptotics
with respect to small parameter $h$ (or ``$h$-Pseudo-differential Calculus''), 
in particular to allow for general phase functions, without homogeneity in the momentum variable. We need also to keep track of the energy parameter.

Let $\iota_0:\Lambda_0\to T^*M$ be a smooth embedded Lagrangian manifold, and $\iota_1:\Lambda_1\to T^*M$ be a smooth embedded Lagrangian manifold
with smooth boundary $\partial\Lambda_1$ (isotropic manifold). Following [MelUhl] we say that $(\Lambda_0,\Lambda_1)$ 
is an {\it intersecting pair} of Lagrangian manifolds iff $\Lambda_0\cap\Lambda_1=\partial\Lambda_1$ and the intersection is {\it clean}, i.e.
$$\forall z\in\partial\Lambda_1 \quad T_z\Lambda_0\cap T_z\Lambda_1=T_z\partial\Lambda_1$$
On the set of intersecting pairs of Lagrangian manifolds 
we define an equivalence relation by saying that $(\Lambda_0,\Lambda_1)\sim({\cal L}_0,{\cal L}_1)$ 
iff near any $z\in\partial\Lambda_1$,
$z'\in\partial{\cal L}_1$, there is a symplectic map $\kappa$ such that $\kappa(z)=z'$, 
and a neighborhood $V\subset T^*M$ of $z$ such that $\kappa(\Lambda_0\cap V)\subset{\cal L}_0$,
$\kappa(\Lambda_1\cap V)\subset{\cal L}_1$. We will call the equivalence class a 
{\it Lagrangian pair}.

All intersecting pairs of manifolds in $T^*M$ are locally equivalent. 
This results from Darboux theorem (see e.g. [H\"o,Vol.3], or [GrSj]) suitably adapted to a pair of Lagrangian manifolds
with non-glancing intersection by [MelUhl,Prop.1.3]. The following Proposition extends [MelUhl, Prop.1.3] in the case of homogeneous
Lagrangian manifolds~:
\medskip
\noindent {\bf  Proposition 3.1}: {\it Let $(\Lambda_0,\Lambda_1)$ be an intersecting pair, and $z\in \Lambda_0\cap\Lambda_1$.
Then there exists a neighbhd $V\subset T^*M$ of $z$, and a canonical map $\kappa: T^*M\to T^*{\bf R}^n$ such that $\kappa(z)=(0,0)$, 
$\kappa(\Lambda_0\cap V)\subset T_0^*{\bf R}^n$ (vertical fiber at 0),
and $\kappa(\Lambda_1\cap V)\subset{\cal L}_+^0$ , ${\cal L}_+^0$ as in (1.15) being the flow-out of $T_0^*{\bf R}^n$ 
by the Hamilton vector field
$v_{\xi_n}=((0,\cdots,0,1),0)$ of $\xi_n$, passing through some $(0;\xi)=(0;(\xi',0))$ i.e. 
$${\cal L}_+^0=\{(x,\xi)\in T^*{\bf R}^n: x=(0,x_n), \xi=(\xi',0), x_n>0\}\leqno(3.15)$$
(superscript 0 is relative to the energy level).}
\smallskip
So this pair of Lagrangian manifolds is in the class of intersecting pairs,
    i.e. in some local canonical charts $\Lambda_0 = {\cal L}_+^0$ and $H(x,p)= E+p_n$.
In particular
    $(\Lambda_0,\Lambda_+^E)$ are mapped onto $({\cal L}_0,{\cal L}_+^0)$ by a canonical transformation
    sending
    $\partial\Lambda_+^E=\Lambda_0\cap\Lambda_+^E$ onto $\partial{\cal L}_+^0={\cal L}_0\cap{\cal L}_+^0$.
\smallskip
As a warm-up,
let us construct $u_h=E_+(h)f_h$ mod ${\cal O}(h^\infty)$ for $hD_{x_n}$, and $f_h$ as in (1.2)
with $a$ compactly supported. 
By a gauge transformation and a shift of the support of $\xi_n\mapsto a(x,\xi)$ in (1.1) 
we may assume $E=0$. So we just need to compute a primitive of $f_h(x)$. Hamilton equations give in this case $X(t)=(X'(t),X_n(t))=(0,t)$ and
$z(t)=(X(t),(\xi',0))$, so in particular $|z(t)|\to\infty$ as $t\to+\infty$,
and the non trapping condition holds, as well as the non return condition. 
Let $T>0$ and $\theta_T\in C^\infty({\bf R}_+)$ vanishing near $+\infty$ and $\theta_T(t)=1$ for $t\leq T$. We consider
$$u(x,h)={i\over h}\int_0^\infty\theta_T(t)\,dt\int^* e^{i(x'\xi'+(x_n-t)\xi_n)/h}a(x',x_n-t,\xi;h)\,d\xi\leqno(3.16)$$
Provided $x_n\leq T/2$ (say) integration by parts and a non-stationary phase argument in variables
$(t,\xi_n)$ show that
$$hD_{x_n}u_h(x)=f_h(x)+{\cal O}(h^\infty)\leqno(3.17)$$
Note that by [H\"o,Lemma 18.2.1] we could already assume $a=a(\xi;h)$. This will be crucial for the {\it compatibility condition}
(see below).

We need to adapt to the semi-classical case
the symbolic calculus considered in [MelUhl].

To this end, we need to generalize (3.16), including local semi-classical distributions of the form
$$u_h(x)={i\over h}\int_0^\infty\theta_T(t)\,dt\int ^* e^{i(x'\xi'+(x_n-t)\xi_n)/h}b(t,x,\xi;h)\,d\xi\leqno(3.18)$$
where $b$ is an amplitude, and $\theta_T(t)$ a cut-off as in (3.16), which we omit for simplicity, since our analysis is local in $x$.

Let us first compute the semi-classical wave-front set $\WF _hu_h$ in (3.18). 
Fix $\overline z=(\overline x,\overline\xi)\in T^*M$, $\overline\xi=(\overline\xi',\overline\xi_n)$. 
It is well known that $\WF _h(u_h)$ is characterized by 
the following property: $\overline z\notin\WF _h(u_h)$ iff there exists
$\chi\in C_0^\infty(T^*M)$ equal to 1 near $\overline z$, such that 
$$\chi(x,hD_x) u_h(x)=(2\pi h)^{-1}
\int\int e^{i(x-y)\eta/h}\chi(y,\eta)u_h(y)\,dy\,d\eta={\cal O}(h^\infty)\leqno(3.20)$$ 
We have the standard
\medskip
\noindent {\bf Proposition 3.2}: {\it Let $u_h$ be as in (3.18). Then 
$$\WF _h(u_h)\subset\{x_n\geq0\}\cap\bigl(\{x_n=0\}\cup\{x'=0\}\cup\{\xi_n=0\}\bigr)\leqno(3.19)$$
If moreover $b(t,x,\xi;h)$ verifies the transport equation $\partial_tb+\partial_{x_n}b=0$, we get the
sharper estimate
$\WF _h(u_h)\subset \Lambda_0\cup\Lambda^0_+$, $\Lambda_0=T_0^*{\bf R}^n$,
which is the conormal bundle of the manifold with boundary $x'=0,x_n\geq0$. Here we recall from Proposition 3.1
that $\Lambda^0_+\subset\char hD_{x_n}$
is the flow out of $\Lambda_0$ by $v_{\xi_n}$ in $\xi_n=0$.} 
\smallskip
\noindent {\it Proof}: Since we work locally in $x$, we can safely omit the cut-off $\theta_T$ in (3.18). 
Let $\Phi_t(x,\xi,y,\eta)=y'\xi'+(y_n-t)\xi_n+(x-y)\eta$ be the phase-function in (3.18). 

\noindent ${\bf (i)}$ If $\overline x_n<0$, we choose $\chi$ such that $\chi(y,\eta)=\widetilde\chi(y',\eta)\chi_{n}(y_n)$, with 
$\chi_{n}(\overline x_n)\neq0$. It follows that $\eta_n\mapsto\Phi_t$ is non stationary in $\supp\chi_{n}$, so $\WF _hu_h\subset\{x_n\geq0\}$. 
\smallskip
\noindent ${\bf (ii)}$ Assume $\overline\xi_n\neq0$ and choose $\chi$ such that $\chi=\widetilde\chi(y,\eta')\chi_{2n}(\eta_n)$, with 
$\chi_{2n}(\overline\xi_n)\neq0$. 
Let $\e $ be so small that $|\xi_n-\eta_n|>\delta>0$ on $|\xi_n|<\e$ and $\eta_n\in \supp\chi_{2n}$,
we split $\chi$ according to $\chi^{\e }$ and $\widehat\chi^{\e }=\chi-\chi^{\e }$, where 
$\chi^{\e }(x,hD_x)u_h(x)=(2\pi h)^{-n}\int dy\int_{|\xi_n|\geq\e }\int e^{i(x-y)\eta/h}\chi(y,\eta)u_h(y)\,d\eta$.

We have  $hD_te^{i\Phi_t/h}=-\xi_n e^{i\Phi_t/h}$, so that integrating by parts $N$ times with respect to $t$
we get $\chi^{\e }(x,hD_x)u_h(x)=\bigl({h\over i}\bigr)^{N}A_N^{\e }(x;h)+B_N^{\e }(x;h)$ where
$$\eqalign{
&A_N^{\e }(x;h)=(2\pi h)^{-n}\int\chi(y,\eta)\,dy\,d\eta\int_0^\infty dt
\int^*_{|\xi_n|\geq\e }e^{i\Phi_t/h}{1\over\xi_n^{n+1}}\partial_t^{N+1}b(t,y,\xi;h)\,d\xi\cr
&B_N^{\e }(x;h)=(2\pi h)^{-n}\int\chi(y,\eta)\,dy\,d\eta\int^*_{|\xi_n|\geq\e }
e^{i\Phi_0/h}{1\over\xi_n}\Sum_{j=0}^N\big({h\over i\xi_n}\bigr)^j \partial_t^jb(0,y,\xi;h)\,d\xi
}\leqno(3.21)$$
Now $(y,\xi)\mapsto\Phi_0(x,\xi,y,\eta)=y\xi+(x-y)\eta$ 
has a non-degenerate critical point at $\xi=\eta,y=0$.
Assume $\overline x\neq0$, i.e. we choose $\chi$ such that $0\notin\pi_x(\supp\chi)$; so $\Phi_0$ is not stationary and 
$B_N^{\e }(x;h)={\cal O}(h^\infty)$. Hence 
$\chi^{\e }(x,hD_x) u_h(x)={\cal O}(h^N)$ for any $N$. 

Consider next the contribution $\widehat\chi^{\e }(x,hD_x)u_h(x)$ of $|\xi_n|<\e $ to $\chi(x,hD_x) u_h(x)$. The map
$$y_n\mapsto\Phi_t(x,\xi,y,\eta)$$ 
has a critical point at $\xi_n=\eta_n$. Since $|\xi_n-\eta_n|>\delta>0$ when $\eta_n\in\supp\chi_{2n}$,
$y_n\mapsto\Phi_t(x,\xi,y,\eta)$ is non stationary and $\widehat\chi^{\e }u_h(x)={\cal O}(h^\infty)$. Altogether,
$\chi(x,hD_x) u_h(x)=\chi^{\e }u_h(x)+\widehat\chi^{\e }u_h(x)={\cal O}(h^\infty)$ so $\overline z\notin\WF _hu_h$ if $\overline\xi_n\neq0$.
In particular
$\WF _hu_h\subset\{x_n\geq0\}\cap\bigl(\{\xi_n=0\}\cup\{x=0\}\bigr)$. 
\smallskip
\noindent ${\bf (iii)}$ Assume next $\overline x_n>0$ and take as above
$\chi(y,\eta)=\widetilde\chi(y',\eta)\chi_{n}(y_n)$, with 
$\chi_{n}(\overline x_n)\neq0$. 
Then $(t,\xi_n,y,\eta)\mapsto\Phi_t$ is critical
at $t=y_n=x_n,\xi_n=0,y'=x',\eta'=x'$, and this is a non-degenerate critical point. So when $x_n>0$, asymptotic stationary phase (3.1)
shows that 
$$\chi(x,hD_x)u_h(x)=2i\pi\int ^* e^{ix'\xi'/h}\widetilde\chi(x',\xi')\bigl[b(x_n,x,(\xi',0);h)+{h\over i}
{\partial^2b_0\over\partial t\partial\xi_n}\bigr)(x_n,x,(\xi',0);h)\bigr]\,d\xi'\leqno(3.22)$$
In particular, using (i), we see that
$\WF _h(u_h)\subset\{x_n\geq0\}\cap\bigl(\{x_n=0\}\cup\{x'=0\}\bigr)$, which altogether proves (3.20). 

For the last statement of Proposition 3.1, apply
$hD_{x_n}$ to $u_h$ and integrate (3.17) by parts once with respect to $t$. We find
$$hD_{x_n}u_h(x)=\int^*e^{ix\xi/h}b(0,x,\xi;h)\,d\xi+
\int_0^\infty dt\int^*e^{i\Phi_t/h}(\partial_t+\partial_{x_n})b(t,x,\xi;h)\,d\xi\leqno(3.23)$$
so if $b(t,x,\xi;h)$ verifies the transport equation, $hD_{x_n}u_h(x)=g_h(x)=\int^*e^{ix\xi/h}b(0,x,\xi;h)\,d\xi$
and since $\WF _hg_h\subset T_0^*{\bf R}^n$, the last (sharper) estimate on $\WF _hu_h$
follows from the well-known property $\WF _h u_h\subset \WF _h hD_{x_n}u_h\cup\char hD_{x_n}$. $\clubsuit$
\smallskip

First we must show that a $h$-FIO $A$ ``quantizing'' the canonical transformation in Proposition 3.1,
i.e. whose canonical map preserves the 
Lagrangian intersection, preserves also bi-Lagrangian distributions of the form (3.18).
Namely, let $({\cal L}, {\cal L}_+)\sim(\Lambda,\Lambda_+)$
be Lagrangian pairs
in the sense of Proposition 3.1, $(y_0,\eta_0)\in\partial{\cal L}_+$, and $u_h$ be of the form (3.18).
Proposition 3.2 of [UhlMel] readily extends as follows:
\medskip
\noindent {\bf Proposition 3.3}: {\it Let $({\cal L},{\cal L}_+)\sim(\Lambda,\Lambda_+)$ be Lagrangian pairs 
in the sense of Proposition 3.1. Let $A$ be a $h$-FIO of the form 
$$Av(x;h)=(2\pi h)^{-(n+N)/2}\int e^{i\phi(x,y,\theta)/h}c(x,y,\theta)v(y;h)\,dy\,d\theta$$
associated with the canonical transformation $\kappa_A$, $\kappa_A(y_0,\eta_0)=(x_0,\xi_0)$ with graph
$$\Lambda_A=\{(x,\phi'_x(x,y,\theta),y;\phi'_y(x,y,\theta):\phi'_\theta(x,y,\theta)=0\}$$
such that (locally) $\Lambda_A\circ{\cal L}=\Lambda$, $\Lambda_A\circ{\cal L}_+=\Lambda_+$
and the compositions are transversal [H\"oIV,p.19\&44].
Here have denoted as usual $\Lambda'_A=\{(x,\xi;y,-\eta): (x,\xi;y,\eta)\in\Lambda_A\}$.  
Let $u_h$ be defined near $(y_0,\eta_0)$ on the Lagrangian pair $({\cal L},{\cal L}_+)$
by (3.18). Then $Au_h$ defined near $(x_0,\xi_0)$ on the Lagrangian pair $(\Lambda,\Lambda_+)$ is again of the form (3.18).}
\smallskip
The proof essentially reduces to show that $Au_h$, after a change of $(t,\theta)$ variables, can be rewritten as an integral of the form (3.18), i.e. 
with the same phase $\Phi_t$, and a new amplitude $b'(t,x,\xi;h)$.

So we can define the class 
$I(M,\Lambda,\Lambda_+)$ of bi-Lagrangian distributions supported on the Lagrangian pair $(\Lambda,\Lambda_+)\sim({\cal L},{\cal L}_+)$,
all of which take locally the form (3.18).

We say that $u_h\in I(M,\Lambda,\Lambda_+)$ is a {\it bi-Lagrangian} (semi-classical) {\it distribution}
on the intersecting pair $(\Lambda,\Lambda_+)$. 

Consider now the inhomogeneous equation $H(x,hD_x;h)u_h=f_h$, where 
$$f_h(x)=\int^* e^{ix\xi/h}a(x,\xi;h)\,d\xi$$ 
is conormal to $\Lambda=\Lambda_0=T^*_0M$.

When $H(x,hD_x;h)=hD_{x_n}$, $u_h$ solves $H(x,hD_x;h)u_h=f_h$ whenever $b$ solves the transport equation, i.e.
$b(t,x,\xi;h)=a(x',x_n-t,\xi;h)$, so that (3.23) simplifies to $hD_{x_n}u_h(x)=f_h(x)$ mod
${\cal O}(h^\infty)$.

In the general case, we can show 
that we can take $H(x,hD_x;)$ to its normal form $hD_{y_n}$ by conjugating with a suitable $h$-FIO as in Proposition 3.3. Namely we have
\medskip
\noindent {\bf Proposition 3.4}: {\it Let the energy surface $H_0=E$ be non critical, and $v_{H_0}$ be transverse to
$\Lambda$ at $(x_0,\xi_0)\in\partial\Lambda_+$, where $\Lambda_+$ is the flow out of $v_H$ from $\Lambda$.
Then there is a $h$-FIO $A$, as in Proposition 3.3, defined microlocally near $\bigl((x_0,\xi_0),(0,0)\bigr)$,
quantizing the canonical transformation of Proposition 1.1 such that $A^{-1}H(x,hD_{x})A=hD_{y_n}$.}
\smallskip
We could presumably take $A$ (microlocally) unitary as in the case without a boundary, but this will not be needed.
\medskip
\noindent{\it 3.4 Compatibility condition and symbolic calculus. Maslov canonical operator for bi-Lagrangian disributions}.
\smallskip
Now we prove Theorem 1.2.
By Propositions 3.3 and 3.4 it suffices to consider the Lagrangian pair $(\Lambda^0,\Lambda_+^0)$, with $\Lambda^0=T^*_0{\bf R}^n$
(which we denote again $\Lambda$ for short). 
When $u_h$ solves $\bigl(H(x,hD_x;h)-E)u_h(x)=f_h$ mod ${\cal O}(h^N)$ we want to
define the ``boundary-part'' $\sigma$ and ``wave-part'' $\sigma^+$ symbols of $u_h$ satisfying the 
{\it compatibility condition}. See also [MelUhl,Sect.4] for a more intrinsic discussion on $(\sigma,\sigma^+)$. 

Our argument is similar to 
[MelUhl,Prop.2.3], and relies on the fact that  if $\chi(x,hD_x)$ is a 
$h$-PDO with $\WF _h(\chi(x,hD_x))\cap\Lambda=\WF '_h(\chi(x,hD_x))\circ\Lambda=\emptyset$, then $\chi(x,hD_x)u_h$ 
is a Lagrangian distribution supported on $\Lambda_+^0$,
while if $\chi(x,hD_x)$ is a 
$h$-PDO with $\WF _h(\chi(x,hD_x)\cap\Lambda_+^0=\emptyset$, then $\chi(x,hD_x)u_h$ is a Lagrangian distribution supported on $\Lambda$.

We proceed in two steps. In Step 1, we construct a family $(\chi_\delta(x,hD_x))_{\delta>0}$ such that $\chi_{\delta }(y,\eta)\to1$ pointwise and 
$\chi_\delta(x,hD_x)u_h$ 
is supported on $\Lambda^0$. This will give the {\it boundary-part symbol}.
In Step 2, we construct a family $(\chi_\delta(x,hD_x))_{\delta>0}$ such that $\chi^+_{\delta }(y,\eta)\to1$ pointwise and 
$\chi^+_\delta(x,hD_x)u_h$ 
is supported on $\Lambda_+^0$. This will give the {\it wave-part symbol}. The {\it compatibility condition} results in comparing these 2 symbols.
\smallskip
\noindent $\bullet$ Step 1: The boundary-part symbol.

We check first the compatibility condition for $N=1$. 
Consider $\chi_{\delta }\in C^\infty({\bf R}^{2n})$ of the form
$\chi_{\delta }(y,\eta)=\chi'_{\delta }(y,\eta')\chi_{2n}^\delta(\eta_n)$ with $|\eta_n|\geq\delta$ on $\supp\chi_{2n}^\delta$,
and $\chi_n^{\delta }(\eta_n)\to1$ pointwise for $\eta_n\neq0$, as $\delta\to0$. Let
$$q(y,\xi;h)={1\over\xi_n}\bigl[b(0,y,\xi;h)+{h\over i\xi_n}\partial_t b_0(0,y,\xi)+{\cal O}(h^2)\bigr]$$
From the part (ii) of the proof of Proposition 3.1 for $N=1$, 
taking $\e <\delta$ so that the contribution of $|\xi_n|<\e $ to  $\chi_\delta(x,hD_h)u_h(x)$ is ${\cal O}(h^\infty)$,
we know that
$$\chi_\delta(x,hD_x)u_h(x)=B_1^{\delta,\e }(x;h)+{\cal O}(h^2)\leqno(3.26)$$
We may assume $\chi_\delta$ is such that  ${\partial\chi_\delta\over\partial y}(0,\eta)=0$, hence computing
$$B_1^{\delta, \e }(x;h)=\int e^{ix\eta/h}\,d\eta\int^*_{|\xi_n|>\e}e^{iy(\xi-\eta)/h}q(y,\xi;h)\chi_\delta(y,\eta)\,dy\,d\xi$$
by asymptotic stationary phase in $(y,\xi)$ with (3.1), where $A=A_\delta={\partial^2\Phi_0\over\partial(y,\xi)^2}$,
$\det A_\delta=(-1)^n$, gives
$$\eqalign{
B_1^{\delta,\e }&(x;h)=(\det A_\delta)^{-1/2}(2i\pi h)^{n}\int d\eta\,e^{ix\eta/h}\chi_\delta(0,\eta)\cr
&\bigl[q(0,\eta;h)-
{h\over i\eta_n}{\partial^2 q_0\over\partial y'\partial\xi'}(0,0,\eta)-{h\over i\eta_n}
{\partial^2 b_0\over\partial y_n\partial\xi_n}(0,0,\eta)+{h\over i \eta_n^2}{\partial b_0\over\partial t}(0,0,\eta)+{\cal O}(h^2)\bigr]\cr
&=(\det A_\delta)^{-1/2}(2i\pi h)^{n}\int d\eta\,e^{ix\eta/h}\chi_\delta(0,\eta)\cr
&\bigl[{1\over\eta_n}b(0,0,\eta;h)+{h\over i \eta_n^2}\bigl(\partial_t b_0(0,0,\eta)+\partial_{y_n}b_0(0,0,\eta)\bigr)-
{h\over i\eta_n}{\partial^2 b_0\over\partial y\partial\xi}(0,0,\eta)+{\cal O}(h^2)]\cr
}\leqno(3.27)$$
The oscillating integral $u_h$ solves $hD_{x_n}u_h=f_h$ iff $b$ satisfies the transport equation 
$$\partial_t b(0,x,\xi;h)+\partial_{x_n} b(0,x,\xi;h)=0$$
i.e. $b(t,x,\xi;h)=b(x',x_n-t,\xi;h)$. Then (3.27) reduces to
$$B_1^{\delta,\e }(x;h)=\int e^{ix\eta/h}\chi_\delta(0,\eta){1\over\eta_n}\bigl[
b(0,0,\eta;h)+ih{\partial^2 b_0\over\partial y\partial\xi}(0,0,\eta)\bigr]\,d\eta+{\cal O}(h^2)\leqno(3.28)$$
which makes sense since $|\eta_n|\geq\delta$ on $\supp\chi_{2n}^\delta$.
So for $\eta_n\neq0$ we define the {\it boundary symbol} of $u_h$ as
$$\sigma(\eta;h)={1\over\eta_n}\bigl[b(0,0,\eta;h)+ih{\partial^2 b_0\over\partial y\partial\xi}(0,0,\eta)\bigr]\leqno(3.29)$$
\smallskip
\noindent $\bullet$ Step 2: The wave-part symbol. 

Consider next $\chi^+_{\delta }\in C^\infty({\bf R}^{2n})$ of the form
$\chi^+_{\delta }(y,\eta)=\chi'^+_{\delta }(y',\eta)\chi_{n}^\delta(y_n)$ with $y_n\geq\delta$ on $\supp\chi_n^\delta$,
and $\chi_n^{\delta }(\eta_n)\to1$ pointwise for $\eta_n>0$, as $\delta\to0$. 
As in the part (iii) in Proposition 3.1, we perform the integration 
$$\chi^+_\delta(x,hD_x)u_h(x)={i\over h}\int_0^\infty dt\int\chi^+_\delta(y,\eta)\,dy\,d\eta\int^*e^{i\Phi_t(x,\xi,y,\eta)/h}b(t,y,\xi;h)\,d\xi$$
by asymptotic stationary phase (3.1) with respect to $(t,\xi_n,y,\eta)$. Here the Hessian is 
$A=A^+_\delta={\partial^2\Phi_t\over\partial(t,\xi_n,y,\eta)^2}$,
$\det A^+_\delta=(-1)^{n+1}$, 
and the critical value of the phase is
$\Phi_c=x'\xi'$. We find
$$\eqalign{
\chi^+_\delta&(x,hD_h)u_h(x)={i\over h}(\det A^+_\delta)^{-1/2}(2i\pi h)^{n+1}\int^*e^{ix'\xi'/h}\chi^+_\delta(x;\xi',0)\cr
&\bigl[b(x_n,x,\xi',0;h)-{2h\over i}
{\partial^2 b_0\over\partial t\partial\xi_n}(x_n,x,\xi',0)\bigr]\,d\xi'+{\cal O}(h^2)
}\leqno(3.30)$$
So we define the {\it wave-part symbol} of $u_h$ by letting $\delta\to0$ as
$$\sigma^+(\xi;h)={i\over h}(\det A^+_\delta)^{-1/2}(2i\pi h)^{n+1} \bigl[b(x_n,x,\xi',0;h)-{2h\over i}
{\partial^2 b_0\over\partial t\partial\xi_n}(x_n,x,\xi',0)\bigr]
\leqno(3.31)$$
and we are to compare (3.29) with (3.31).

From [H\"o,Lemma 18.2.1] and its proof, we know that
if $v_h(x)=\int e^{ix\xi/h} a(x,\xi;h)\,d\xi$, then we also have $v_h(x)=\int e^{ix\xi/h} \widetilde a(\xi;h)\,d\xi$,
with a symbol $\widetilde a(\xi;h)\sim\Sum_j h^j\langle-iD_x,D_\xi\rangle a(x,\xi)/j!$ independent of $x$. 
Applying this to amplitude $b(x',x_n-t,\xi)$, 
the terms ${\partial^2 b_0\over\partial y\partial\xi}(0,0,\eta)$ and ${\partial^2 b_0\over\partial t\partial\xi_n}(x_n,x,\xi',0)$
resp. in (3.28) and (3.30) disappear, $\sigma^+(\xi;h)$ is continuous up to $\partial\Lambda_+$, 
and we end up with the compatibility condition (1.23) mod ${\cal O}(h^2)$ between the 
wave-part and boundary-part symbols on $\partial\Lambda_+$.

It is clear, following again the proof of Proposition 3.1, that (3.26) carries by induction mod ${\cal O}(h^N)$, all $N$,
when replacing $B_1^{\delta,\e }(x;h)$ by $B_N^{\delta,\e }(x;h)$. 

This allow to define coherently the (bi-)symbol $(\sigma,\sigma^+)=(\sigma(u_h), \sigma^+(u_h))$ 
computed as above in local coordinates, and thus by analogy with Sect.3.2, an ``effective''
Maslov canonical operator 
$K^h_{\Lambda,\Lambda^0_+}(\sigma,\sigma^+)$. 
The commutation formula for bi-Lagrangian distributions takes the form
$$\eqalign{
H&(x,hD_x;h)\bigl[K_{\Lambda,\Lambda_+}^h (\sigma,\sigma^+)\bigr](x;h)=\bigl[K_{\Lambda,\Lambda_+}^h 
\bigl((H-E)\sigma, h\bigl({1\over i}{\cal L}_{v_{H_0}}|_{\Lambda_+}+H_1\bigr)\sigma^+\bigr](x;h)+{\cal O}(h^{2})\cr
&=f_h(x)+{\cal O}(h^{2})
}\leqno(3.33)$$
which gives (1.24) once the transport equation has been solved as in (3.12). This brings the proof of Theorem 1.2 to an end. $\clubsuit$
\smallskip
\noindent {\it 3.5 The constant coefficient case}. 
\smallskip
In general it is difficult to obtain a decomposition of $u_h$ adapted to the splitting $u_h=u_h^0+u_h^1$ where 
$\WF _h(u_h^0)\subset\Lambda_0$ and $\WF _h(u_h^1)\subset\Lambda_+^E$.  

Here we compute $u_h$ explicitely in the 2-D case for Helmholtz operator $-h^2\Delta-E$, 
but $f$ with compact support. Let $f$ also be radially symmetric;
its Fourier transform $g={\cal F}_1f$ is again of the form 
$g(p)=g(|p|)=g(r)$ and extends holomorphically to ${\bf C}^2$. For $E=k^2$, $k>0$,
we rewrite 
$$u_h(x)=(2\pi h)^{-n}\int e^{ix\xi/h}{{\cal F}_1f(\xi)\over\xi^2-E-i0} \,d\xi$$ 
as $u_h(x)=u(x)=u_0(x)+u_1(x)$ with
$$\eqalign{
&u_0(x)={k+i\e \over(2\pi h)^2}\int_0^{2\pi}\,d\theta\int_0^\infty\exp[i|x|r\cos\theta/h]{g(r)\over r^2-(k+i\e )^2}\,dr\cr
&u_1(x)={1\over(2\pi h)^2}\int_0^{2\pi}\,d\theta\int_0^\infty\exp[i|x|r\cos\theta/h]{g(r)\over r+k+i\e }\,dr\cr
}$$
To compute $u_0$ we use contour integrals. When $\theta\in]-{\pi\over2},{\pi\over2}[$, we shift the contour of integration
to the positive imaginary axis and get by the residues formula
$$\eqalign{
\int_0^\infty&\exp[i|x|r\cos\theta/h]{g(r)\over r^2-(k+i\e )^2}\,dr+\int_0^\infty\exp[-|x|r\cos\theta/h]{g(ir)\over r^2+(k+i\e )^2}\,idr=\cr
&2i\pi{g(k+i\e )\over2(k+i\e )}\exp[i|x|(k+i\e )\cos\theta/h]
}\leqno(3.35)$$
while for $\theta\in]{\pi\over2},{3\pi\over2}[$,
$$\int_0^\infty\exp[i|x|r\cos\theta/h]{g(r)\over r^2-(k+i\e )^2}\,dr-\int_0^\infty\exp[|x|r\cos\theta/h]{g(-ir)\over r^2+(k+i\e )^2}\,idr=0
\leqno(3.36)$$
Summing up (3.35) and (3.36), integrating over $\theta\in]0,2\pi[$ and letting $\e \to0$, we obtain
$$\eqalign{
u_0(x)&={i\pi g(k)\over (2\pi h)^2}\int_{-\pi/2}^{\pi/2}\exp[i|x|k\cos\theta/h]\,d\theta+\cr
&\int_0^\infty{dr\over r^2+k^2}\bigl[\int_{-\pi/2}^{\pi/2}g(ir)-\int_{\pi/2}^{3\pi/2}g(-ir)\bigr]\exp[-|x||\cos\theta|/h]\,d\theta
}$$
Since $g(ir)=g(-ir)$, the latter integral vanishes, so we end up with
$$u_0(x)={i\pi g(k)\over (2\pi h)^2}\int_{-\pi/2}^{\pi/2}\exp[i|x|k\cos\theta/h]\,d\theta$$
It is readily seen that 
$$\WF _h u_0\subset\{x=0\}\cup\{(x,k{x\over|x|}), x\neq 0\}=\Lambda_0\cup\Lambda_+^E$$
Consider now $u_1$. We let $\e \to 0$ and set $\widetilde g(r)={g(r)\over r(r+k)}$. Since
$\widetilde g(r)\sqrt r\in L^1({\bf R}_+)$, we have  
$u_1(x)=H_0(\widetilde g)({|x|\over h})$, where $H_0$ denotes Hankel transform of order 0. 
Let $\chi\in C_0^\infty({\bf R}^2)$ be
radially symmetric, and equal to 1 near 0, since $\WF _h f_h=\{x=0\}$, we have 
$$g={\cal F}_h(\chi f_h)+{\cal O}(h^\infty)=(2\pi h)^{-2}{\cal F}_h(\chi)*g+{\cal O}(h^\infty)$$
so in the expression for $u_1$ we may replace mod ${\cal O}(h^\infty)$, $\widetilde g(r)$ by a constant times
$\widehat g(r)={({\cal F}_h(\chi)*g)(r)\over r(r+m)}$ (see [Bad] for 2-D convolution and Fourier transform in polar coordinates). 
To estimate $\WF _h u_1$, we compute again the Fourier transform of $(1-\widetilde\chi)\widehat g$ where $\widetilde \chi$
is a cut-off equal to 1 near 0, and we find it is again ${\cal O}(h^\infty)$ if $\chi\equiv1$ on $\supp\widetilde\chi$.
This shows that  
$\WF _hu_1\subset \{x=0\}$. 

Note that this Example makes use of Bessel function $J_0({|x|\over h})$, we shall return to such ``localized functions'' in Sect.5. 
\bigskip

\noindent {\bf 4. $f_h$ is supported microlocally on the ``vertical plane''}
\smallskip
Here we shall construct objects written globally in a suitable coordinate system, using the phase functions in Sect.2.
This system consists only of coordinates on $\Lambda$ and of time parameter $t$ in Hamilton equations.
This Sect. relies in a strong way on [DoNaSh]. 

Consider the case where $H_0$ is positively homogeneous of degree $m\geq1$ with respect to $p$ and
$f_h$ is microlocally concentrated on the vertical plane $\Lambda=\{x=0\}$,
e.g. $f_h(x)=h^{-n}f({x\over h})$, with $f\in{\cal S}({\bf R}^n)$ (Schwartz space). 
\medskip
\noindent {\it 4.1 Some non-degeneracy condition}.
\smallskip
Recall first from [DoNaSh, Lemma 6] the following result. Let $\widetilde \iota:\widetilde\Lambda\to T^*\widetilde M$ 
be a Lagrangian embedding of dimension $\widetilde n$, $U\subset\widetilde\Lambda$ 
a connected simply connected open set, 
$$(\widetilde\phi,\widetilde\psi)=(\phi_1,\cdots,\phi_k,\psi_1,\cdots,\psi_{\widetilde n-k})$$
local coordinates on $U$. Here the additional assumption of [DoNaSh, Lemma 6] that 
$k$ is the rank of $\pi_x:\widetilde\Lambda\to\widetilde M$, is not required. Thus $\widetilde\Lambda$ is defined by $x=X(\widetilde\phi,\widetilde\psi), 
p=P(\widetilde\phi,\widetilde\psi)$ in the chart $U$. 
Let $\Pi(\widetilde\phi,\widetilde\psi)$ be a smooth $\widetilde n\times k$ matrix function defined in $U$ such that: 
$$\eqalign{ 
&\Pi^*(\widetilde\phi,\widetilde\psi)X_{\widetilde\phi}(\widetilde\phi,\widetilde\psi)=\Id_{k\times k}\cr
&\kappa:(\widetilde\phi,\widetilde\psi)\mapsto 
(X(\widetilde\phi,\widetilde\psi),\widetilde\psi) \ \hbox{is an embedding}
}\leqno(4.1)$$ 
Then there is a neighbhd $V$ of $\kappa(U)$ such that the system
$$\Pi^*(\widetilde\phi,\widetilde\psi)\bigl(x-X(\widetilde\phi,\widetilde\psi)\bigr)=0, \quad (x,\widetilde\psi)\in V$$
has a unique smooth solution $\widetilde\phi=\widetilde\phi(x,\widetilde\psi)$ satisfying the condition 
$X(\widetilde\phi(x,\widetilde\psi),\widetilde\psi)=x$, when $(x,\widetilde\psi)\in \kappa(U)$.

For $(\widetilde\phi,\widetilde\psi)\in U$, define the $\widetilde n\times \widetilde n$ matrix 
$${\cal M}(\widetilde\phi,\widetilde\psi)=\bigl(\Pi(\widetilde\phi,\widetilde\psi);
P_{\widetilde\psi}(\widetilde\phi,\widetilde\psi)-P_{\widetilde\phi}(\widetilde\phi,\widetilde\psi)\Pi^*(\widetilde\phi,\widetilde\psi)
X_{\widetilde\psi}(\widetilde\phi,\widetilde\psi)\bigr)$$
As we shall see, invertibility of ${\cal M}(\widetilde\phi,\widetilde\psi)$ plays an important role [DoMaNa2], [DoNaSh]. 

Consider now our special setting where $H$ is positively homogeneous of degree $m$, $\Lambda_0$ is the vertical plane, and recall 
$\langle P(t,\psi,\tau),\dot X(t,\psi,\tau)\rangle=mH$ from (2.8). 
Here $\tau$ small enough is taken as a parameter, everything depends smoothly on $\tau$ and $\Lambda_+(0)=\Lambda_+^E$. So
$\Pi(t,\psi,\tau)={1\over mH}P(t,\psi,\tau)$ is a left inverse of $\dot X$: $\Pi^*\dot X={1\over mH}\langle P,\dot X\rangle=1$.
Further, the map $\Lambda_+^E\to{\bf R}^n$, $(t,\psi)\mapsto(X(t,\psi,\tau),\psi)$ is clearly an embedding.  
This fulfills conditions (4.1) above for $\widetilde\Lambda=\Lambda_+(\tau)$, with $\widetilde n=n$, $k=1$
and $\widetilde\phi=t$, $\widetilde\psi=\psi$. 
So the system
$$\Pi^*(t,\psi,\tau)(x-X(t,\psi,\tau))=\langle P,x-X(t,\psi,\tau)\rangle=0\leqno(4.2)$$
has a unique solution $t=t_0(x,\psi,\tau)$ satisfying the condition 
$X(t_0(x,\psi,\tau),\psi,\tau)=x$,
and this solution is a smooth function. When $\tau=0$ we omit it from the notations, and write for instance $X(t,\psi)$ for $X(t,\psi,0)$.

Moreover by (2.8) again, the matrix
$${\cal M}(t,\psi,\tau)=\bigl(\Pi(t,\psi,\tau),P_\psi(t,\psi,\tau)-\dot P{1\over mH}\,{}^tP\,X_\psi\bigr)=({1\over mH}P,P_\psi)\leqno(4.3)$$ 
has determinant ${1\over mH}\det(P,P_\psi)$. As we shall see, it turns out that $\det{\cal M}$ gives the invariant (inverse) density on $\Lambda_+$.
\smallskip
\noindent {\it Example 4.1}: Let us compute $\det (P,P_\psi)$ at $t=0$ for a geodesic flow $H(x,p)$, on the energy shell $E=1$ when $n=2$ or $n=3$. 
When $n=2$, up to a change of $x$ coordinates such that at $x=0$, the metric $H(0,p)$
takes the diagonal form $H(x,p)={p_1^2\over a_1^2}+{p_2^2\over a_2^2}$ (elliptic polarization), and 
$P=(a_1\cos\psi,a_2\sin\psi)$ for some $a_1,a_2>0$. Hence $\det (P,P_\psi)=a_1a_2$ at $x=0$, and for  
small $|t|$ we have $\det (P,P_\psi)>0$.
For $n=3$, $H(0,p)={p_1^2\over a_1^2}+{p_2^2\over a_2^2}+{p_3^2\over a_3^2}$, for some $a_1,a_2,a_>0$ and in spherical coordinates $(\psi_1,\psi_2)$
where $0<\psi_1<\pi$,
we find $\det (P,P_\psi)=a_1a_2a_3\sin\psi_1>0$ for small $t$ and away from the poles $(0,0,\pm1)$.
\smallskip
\noindent {\it Example 4.2}: When $H(x,p)={|p|^m\over\rho(x)}$, recall from (2.5) that $P(\psi,\tau)=|P|_\tau\omega(\psi)$ at $t=0$. 
Since $\det(\omega(\psi),\omega^\perp(\psi))=1$, again we have
$\det(P,P_\psi)\neq0$ for small $t$. 
\medskip
\noindent {\it 4.2 Construction of the phase function and half-density, general case}.
\smallskip
We first construct by HJ theory a generating function $\Phi$ of $\Lambda_+$ that verifies the initial condition
$\Phi|_{t=0}=\langle x,\omega(\psi)\rangle$.
Our approach consists in looking for a parametric form of the phase, depending on the initial data through the ``front variables'' 
$(X(t,\psi,\tau),P(t,\psi,\tau))$ only.

The most natural Ansatz (recall $\tau+H=E$), would be 
$$\Phi_0(x,t,\psi,\tau)=mEt+\langle P(t,\psi,\tau),x-X(t,\psi,\tau)\rangle\leqno(4.4)$$
with initial condition $\Phi|_{t=0}=\langle p,x\rangle$, $p=P(\psi,\tau)$ arbitrary. 
The ``$\theta$ variables'' in H\"ormander's definition are then $(\psi,\tau)$. 

In the simplest example where $n=1$, $\tau+H(x,p)=\tau+p=E$,
$\Phi_0=Et+p(x-t)$ (there are no variable $\psi$, and $X(t)=t$ is independent of $\tau$). 
This is actually a parametrization of $\exp tv_H(z)$, $z\in T^*M$, for $t\in{\bf R}$. 
But the drawback of $\Phi_0$ is to depend on $\tau$ (that has eventually to bet set to 0)
in a complicated way, when taking variations with respect to parameters. 

The second one [DoNa2] consists  in choosing a new ``radial'' coordinate $\lambda=\lambda(\tau)$, $\lambda(0)=1$, on $\Lambda_0$ 
completing the $\psi$ variables,
such that $\partial\Lambda_+^E$ is given by $\lambda=1$. We could think of $\lambda$ as a Lagrange multiplier. We define
$$\Phi(x,t,\psi,\lambda)=mEt+\lambda\langle P(t,\psi),x-X(t,\psi)\rangle\leqno(4.5)$$
where now $(X,P)$ are evaluated on $\Lambda_+^E$ (and not on $\Lambda_+(\tau)$). 
The ``$\theta$ variables'' in H\"ormander's definition are then $(\psi,\lambda)$. In the Example above, $\Phi=Et+\lambda p(x-t)$.
The critical value
of $\Phi$ with respect to $\theta$, is viewed either as a function on the critical set 
$\widetilde C_\Phi$, with the Lagrangian embedding 
$$\widetilde C_\Phi=\{(x,t,\theta): \partial_\theta\Phi=0\}\to \widetilde\Lambda_+\leqno(4.6)$$ 
where we recall $\widetilde\Lambda_+$ from (2.13), or on
$$C_\Phi=\{(x,t,\theta): \partial_\theta\Phi=\partial_t\Phi=0\}\to \Lambda_+^E\leqno(4.7)$$ 
In both cases, (2.12) holds precisely on the critical set.  

Eikonal equation (2.12) verified at second order on $C_\Phi$ reads
$$\partial_t\Phi+H(x,\partial_x\Phi)-E={\cal O}(|x-X(t,\psi),\lambda-1|^2)\leqno(4.8)$$
Variables $\tau$ and $\lambda$ are diffeomorphically mapped onto each other. In case (1.25) this goes as follows~: comparing
(4.6) with (4.7) at $t=0$, we get $P(\psi,\tau)=\lambda P(\psi)$, so by (2.3)
$$\lambda=\bigl(1-{\tau\over E}\bigr)^{1/m}\leqno(4.9)$$ 
A similar correspondence holds in Example 4.1.
\smallskip
\noindent {\bf Proposition 4.1}: {\it Let $H(x,p)$ be positively
homogeneous of degree $m\geq1$ with respect to $p$ on $T^*M\setminus 0$, and $\det(P,P_\psi)\neq0$ at some point $(t',\psi')$. Then 
$\Phi(x,t,\psi,\lambda)$ 
given in (4.5) is a non-degenerate phase function defining $\Lambda_+^E$ near $(t',\psi')$, 
with initial condition $\Phi|_{t=0}=\langle x,p\rangle$, thus is the 1-jet on $\Lambda_+$ of the solution of HJ equation (2.12). 
The positive invariant (inverse) density on $\Lambda_+^E$ we recall from (3.8) is given by 
$$F[\Phi,|dy|]|_{C_\Phi}=mE\det (P,P_\psi)\neq0\leqno(4.10)$$
The critical set $C_\Phi$ is then determined by $x=X(t,\psi)$ (which can be inverted as 
$t=t_1(x,\psi)$) and $\lambda=1$. It coincides with the set $\kappa(U)$ defined after (4.1).}
\smallskip
Recall that in case of Hamiltonian (1.25), the condition $\det(P,P_\psi)\neq0$ holds at $t=0$, since $\det (P,P_\psi)=|P|^2$ there, see Example 4.2.
Thus $\Lambda_+^E$ is parametrized by $\Phi(x,t,\psi,\lambda)$ for small $t$. Recall $|dy|$ is Lebesgue measuer on $C_\Phi$. 
\smallskip
\noindent {\it Proof}: 
We have using (2.6), evaluated at $\tau=0$
$$\partial_t\Phi=\dot\Phi=mE+\lambda\langle\dot P,x-X(t,\psi)\rangle-\lambda\langle P,\dot X\rangle=mE(1-\lambda)+\lambda
\langle\dot P(t,\psi),x-X(t,\psi)\rangle\leqno(4.11)$$
so $\partial_t\Phi=0$ and 
$\partial_x\Phi=P(t,\psi)$ along $x=X(t,\psi)$ when $\lambda=1$.
We are left to show that $\Phi$ is non degenerate phase function, with $(\psi,\lambda)$ as ``$\theta$-parameters''. From (4.5)
$$\eqalign{
&\partial_\psi\Phi=\lambda\langle \partial_\psi P(t,\psi),x-X(t,\psi)\rangle\cr
&\partial_\lambda\Phi=\langle P(t,\psi),x-X(t,\psi)\rangle\cr
}\leqno(4.12)$$
Let us add $t$ to the ``$\theta$-variables'', and consider
the variational system $\partial_t\Phi=\partial_\psi\Phi=\partial_\lambda\Phi=0$, which determines the critical set $C_\Phi$. 
Last 2 equations $\partial_\psi\Phi=0,\partial_\lambda\Phi=0$ give an homogeneous linear system
in $x-X(t,\psi)$ with determinant $\det (P,P_\psi)$. 

So for $(t,\psi)$ near $(t',\psi')$, the phase is critical with respect to $(\psi,\lambda)$ precisely for $\lambda=1$ and
$x=X(t,\psi)$, in particular it is critical along $\partial\Lambda_+^E$ when $\lambda=1$. Recall $t_0$ from (4.2). By the discussion above 
and [DoNaSh, Lemma 6], we find that  
$\langle P(t,\psi),x-X(t,\psi,\tau)\rangle=0$ when $\lambda=1$ has a unique solution $t=t_1(x,\psi)=t_0(x,\psi,\tau=0)$ satisfying the condition:
if $(x; \psi,\lambda=1)\in C_\Phi$, then $X\bigl(t_1(x,\psi),\psi\bigr)=x$. 
Moreover $t_1$ is the critical point of $t\mapsto\Phi$ when $\lambda=1$.

Condition $\det (P,P_\psi)\neq0$ actually ensures that $\Phi$ is a non degenerate phase function, i.e. the vectors
$\bigl(d\partial_{t}\Phi, d\partial_{\psi}\Phi, d\partial_{\lambda}\Phi\bigr)$ are linearly independent on the set $x=X(t,\psi)$.
Namely, look at the variational system and use (4.11) and (4.12) to compute on $C_\Phi$ the differentials
$$\eqalign{
&d\dot\Phi=-mE\,d\lambda+\lambda\langle\dot P(t,\psi),dx-dX(t,\psi)\rangle\cr
&d(\partial_\psi\Phi)=\lambda\langle P_\psi(t,\psi),dx-dX(t,\psi)\rangle\cr
&d(\partial_\lambda\Phi)=\langle P,dx-dX(t,\psi)\rangle
}\leqno(4.13)
$$
Introduce the Jacobian (3.8), quotient of two $2n+1$ forms. 
$$F[\Phi,d\mu_+]|_{C_\Phi}={dt\wedge\,d\psi\wedge\,d\dot\Phi\wedge\,d(\partial_\lambda\Phi)\wedge\,d(\partial_\psi\Phi)\over
dx\wedge\,dt\wedge\,d\psi\wedge\,d\lambda}\leqno(4.14)$$
Here $dx$ is the volume form. 
Substituting (4.13) into $\omega=dt\wedge\,d\psi\wedge\,d\dot\Phi\wedge\,d(\partial_\lambda\Phi)\wedge\,d(\partial_\psi\Phi)$ we get 
$$\eqalign{
&\omega=-mE\,d\psi\wedge\,dt\wedge\,d\lambda\wedge\langle P,dx-dX(t,\psi)\rangle\wedge\langle P_\psi,dx-dX(t,\psi)\rangle+\cr
&d\psi\wedge\,dt\wedge\langle \dot P,dx-dX(t,\psi)\rangle\wedge \langle P,dx-dX(t,\psi)\rangle\wedge \langle P_\psi,dx-dX(t,\psi)\rangle
}$$
Writing $dX=\dot X\,dt+X_\psi\,d\psi$, we check that the second term vanishes, so we are left with
$$\omega=-mE\,\det(P,P_\psi)\,dt\wedge\,d\psi\wedge\,d\lambda\wedge\,dx$$
which gives (4.10). 
So if $\det (P,P_\psi)\neq0$,
$\Phi$ is a non-degenerate phase function, and (4.10) the invariant (inverse) density on $\Lambda_+^E$.
$\clubsuit$
\medskip
\noindent {\it 4.3 New parametrizations, general case in 2-D}.
\smallskip
We investigate some configurations of $\Lambda_+^E$,
and describe the corresponding Lagrangian singularities.
Consider first the critical points of the phase.  Let (see Proposition A.1) 
$$a=\langle\dot P,X_\psi\rangle, \ c=\langle P_\psi,X_\psi\rangle, \ d=\langle\dot P,\dot X\rangle, \ \alpha=\det(P,P_\psi), \ 
\beta=\det(P,\dot P), \ \gamma=\det(\dot P,P_\psi) \leqno(4.16)$$
At the critical point
$$-\Hess_{(t,\psi,\lambda)}\Phi=\pmatrix{\lambda\langle\dot P,\dot X\rangle&\lambda\langle\dot P,X_\psi\rangle&\langle P,\dot X\rangle\cr
\lambda\langle \dot P,X_\psi\rangle&\lambda\langle P_\psi,X_\psi\rangle&0\cr
\langle P,\dot X\rangle&0&0}\leqno(4.17)$$
thus $\det\Hess_{(t,\psi,\lambda)}\Phi=(mE)^2\langle P_\psi,X_\psi\rangle$. 

When $T_z\Lambda_+$ is not transverse to the vertical plane
$V_z=\{(0,\delta p)\}\subset T_z^*M$, we know from Sect.3.2 that we need to
express $u_h$ in Fourier representation.
This will be needed in Sect.4.5 to derive the commutation formula at some $z\in\Lambda_+$.

Assume 
$\alpha=\det(P,P_\psi)\neq0$ (which holds near $t=0$). Since $\langle P,X_\psi\rangle=0$, the relation $\langle P_\psi,X_\psi\rangle=0$ implies
$X_\psi=0$, so at such a point, we need to change some of the ``$\theta$-variables'' $(t,\psi,\lambda)$. 
If $z\in\Lambda_+^E$ is such that $X_\psi=0$, $T_z\Lambda_+^E$ is not transverse to the vertical plane: indeed
$\dim(V_z\cap T_z\Lambda_+^E)=1$, for $\langle P,\dot X\rangle=1$ and $\dot X\neq0$.

We proceed as in Sect.3. 
Consider the embedding $C_\Phi\to\Lambda_+^E$ as in (4.7).
Let $x',x''$ be a partition of $x$, we introduce a partial
Legendre transformation as in Sect.3, implement the latter equations for the critical point by $\dot\Phi=0$, and compute the Hessian
$${\cal H}(x'',\xi')=\Hess _{(x',t,\psi,\lambda)}\bigl(\Phi(x,t,\psi,\lambda)-x'\xi'\bigr)$$ 
First we try $x'=x_2,x''=x_1$ so that 
$${\cal H}(x_1,\xi_2)=\pmatrix{0&\lambda\dot P_2&\lambda\partial_\psi P_2&P_2\cr
\lambda\dot P_2&-\lambda\langle \dot P,\dot X\rangle&-\lambda\langle\dot P,X_\psi\rangle&-mE\cr
\lambda\partial_\psi P_2&-\lambda\langle\dot P,X_\psi\rangle&-\lambda\langle P_\psi,X_\psi\rangle&0\cr
P_2&-mE&0&0}$$
and we find, with notations (4.16)
$$\lambda^{-2}\det{\cal H}(x_1,\xi_2)=P_2^2(a^2-cd)+2mEcP_2\dot P_2-mEaP_2^2-mEaP_2\partial_\psi P_2+mE(\partial_\psi P_2)^2\leqno(4.20)$$
Similarly, choosing $x'=x_1, x''=x_2$, we get the same expression with $P_2$ replaced by $P_1$. Now if $X_\psi=0$, then $a=c=0$,
and since $P_\psi\neq0$ (we assume here $\det(P,P_\psi)\neq0$, there is a partition of variables $x',x''$ such that 
$\det{\cal H}(x',\xi'')\neq0$. Actually variables $(x',\xi'')$ are implicit in the expression of ${\cal H}(x',\xi'')$, but fixing $(x',x'')$
on the critical set determines the front variables $(X(t,\psi),P(t,\psi)$.~)
\medskip
\noindent {\it Remark 4.2}: Compute instead ${\cal H}(\xi)=\Hess_{x,t,\psi,\lambda}\bigl(\Phi(x,t,\psi,\lambda)-x\xi\bigr)$ 
at the critical point. We have
$${\cal H}(\xi)=\pmatrix{0&0&\lambda\dot P_1&\lambda\partial_\psi P_1&P_1\cr
0&0&\lambda\dot P_2&\lambda\partial_\psi P_2&P_2\cr
\lambda \dot P_1&\lambda\dot P_2&-\lambda\langle \dot P,\dot X\rangle&-\lambda\langle\dot P,X_\psi\rangle&-mE\cr
\lambda\partial_\psi P_1&\lambda\partial_\psi P_2&-\lambda\langle\dot P,X_\psi\rangle&-\lambda\langle P_\psi,X_\psi\rangle&0\cr
P_1&P_2&-mE&0&0}\leqno(4.21)$$
and $-\lambda^{-3}\det{\cal H}(\xi)=\beta^2c+\alpha^2d-2\alpha\beta a-2mE\alpha\gamma$.
When $a=c=0$, vanishing of the determinant (4.21) reduces to
$$-\lambda^{-3}\det{\cal H}(\xi)=\det^2(P,P_\psi)\langle-\partial_xH,\partial_pH\rangle-2mE\det(P,P_\psi)\langle-\partial_xH,P\rangle$$
In the particular case of Hamiltonian (1.25),  at $t=0$, using $\det(P,P_\psi)=|P|^2$, we find
$$\lambda^{-3}|P|^{-1}\det{\cal H}(\xi)={mE^2\over\rho(0)}\langle P,\nabla\rho(0)\rangle$$
which vanishes at a special or residual point. See Remark 3.1.   
\medskip
\noindent {\it Remark 4.3}: Consider instead the embedding $\widetilde C_\Phi\to\widetilde\Lambda_+$ as in (4.6), and 
compute the critical points of $(x,\psi,\lambda)\mapsto\Phi(x,t,\psi,\lambda)-x\xi$.
We add $\lambda P(t,\psi)-\xi=0$ to the previous equations $\partial_\psi\Phi=0, \partial_\lambda\Phi=0$. 
and compute the Hessian 
$${\cal H}_0(t,\xi)=\Hess _{(x,\psi,\lambda)}\bigl(\Phi(x,t,\psi,\lambda)-x\xi\bigr)$$ 
at $x=X(t,\psi)$, namely
$${\cal H}_0(t,\xi)=\lambda\pmatrix{0&0&\partial_\psi P_1&P_1\cr
0&0&\partial_\psi P_2&P_2\cr
\partial_\psi P_1&\partial_\psi P_2&-\langle P_\psi,X_\psi\rangle&0\cr
P_1&P_2&0&0}$$
so that 
${\cal H}_0(t,\xi)=\lambda^4\det(P,P_\psi)$
is non-degenerate, but as a function 
on $\widetilde C_\Phi$
(extended phase-space) instead of $C_\Phi$. This is related to the general fact that $\widetilde\Lambda_+$
is always projectable (for small $t$) on ${\bf R}^{n+1}$,
which is not the case for $\Lambda_+^E$ on ${\bf R}^{n}$.
\medskip
On the other hand, in order to investigate Lagrangian singularities of $\Lambda_+^E$, we need to eliminate some of the ``$\theta$-variables''.
For short we will do it only in the case where $\Lambda_+^E$ is tranverse to the vertical fiber of $T^*M$.

So let $z(t)=\bigl(X(t,\psi),P(t,\psi)\bigr)$ be such that $T_{z(t)}\Lambda_+$ is transverse to the vertical plane 
$V_{z(t)}$ (namely
$X_\psi\neq0$), or $t=0$ but $\det{\cal H}(\xi)\neq0$, see (4.21); we parametrize $\Lambda_+^E$ with $\Phi(x,t,\psi,\lambda)$, 
and $\Phi_{\theta\theta}\neq0$. When $z(t)$ is a focal point, we discuss according to the case $z(t)$ is a special point (in the sense
of Definition 1.4) or not.
\medskip
\noindent {\bf Proposition 4.3}: {\it Let $n=2$ for simplicity. Let $z\in\Lambda_+^E$ (possibly on $\partial\Lambda_+^E$) and assume $\Phi$
is a non-degenerate phase function near $z$ (which holds true when $z\in\partial\Lambda_+^E$ except for exceptional points where $X_\psi=0$). We have:

\noindent (i) Let $z\in\Lambda_+^E$ such that $\langle-\partial_xH(z),\partial_pH(z)\rangle=\langle\dot X,\dot P\rangle\neq0$. Then near $z$
the rank of $d\pi|_{\Lambda_+^E}$ is 1 when $c=0$ (i.e. $X_\psi=0$), and 2 when $c\neq0$.

\noindent (ii) Let $z_0=\bigl(X(t_0,\psi_0),P(t_0,\psi_0)\bigr)\in\Lambda_+^E$ be a special point for some $(t_0,\psi_0)$. 
Then the rank of $d\pi|_{\Lambda_+}(z)$ is 1 or 2.
When the rank is 1, the tangent space of the caustics at $X(t_0,\psi_0)$ takes the form
$${\partial x_1\over\partial \psi}{\partial\lambda\over\partial x_1}+{\partial x_2\over\partial \psi}{\partial\lambda\over\partial x_2}=0\leqno(4.22)$$
where ${\partial\lambda\over\partial x}(x,\psi)\neq0$. 

\noindent (iii) Let $z_0=\bigl(X(t_0,\psi_0),P(t_0,\psi_0)\bigr)\in\Lambda_+$ be a residual point for some $(t_0,\psi_0)$, i.e. $\dot P(t_0,\psi_0)=0$.
Then the eikonal is $mE\,dt=\langle P,dx\rangle\neq0$ at $z_0$.   } 
\smallskip 
\noindent {\it Proof}: {\bf (i)} On $C_\Phi$ we have $\partial^2_t\Phi=-\langle \dot P,\dot X\rangle\neq0$, so implicit function theorem shows that 
(for small $t$)
$\dot\Phi=0$ is equivalent to $t=t(\psi,\lambda)$. 
Since we have eliminated $t$, the ``$\theta$-parameters'' are now $(\psi,\lambda)$, and we set 
$$\Psi(x,\psi,\lambda)=\Phi(x,t(x,\psi,\lambda),\psi,\lambda)$$
Differentiating the relation $\partial_t\Phi=0$, we get that on $C_\Phi$ and for $\lambda=1$
$${\partial t\over\partial\psi}=-{\langle P_\psi,\dot X\rangle\over\langle \dot P,\dot X\rangle}, \quad 
{\partial t\over\partial\lambda}=-{mE\over\langle \dot P,\dot X\rangle}$$
and a  straightforward computation using (2.6) yields
$$\Psi''_{\theta\theta}=\Psi''_{(\psi,\lambda)}=\pmatrix{
{\langle \dot P,X_\psi\rangle^2\over\langle\dot P,\dot X\rangle}-\langle P_\psi,X_\psi\rangle& mE
{\langle \dot P,X_\psi\rangle\over\langle\dot P,\dot X\rangle}\cr 
*&{(mE)^2\over\langle\dot P,\dot X\rangle}}\leqno(4.23)$$
Applying (3.5) to the non-degenerate phase function $\Phi$ with $N=n=2$, we find that the rank of $d\pi|_{\Lambda_+}(z)$ is 1 
when $c=0$ or 2 when $c\neq0$. 
\smallskip
\noindent {\bf (ii)} We could attempt to solve
$\partial_t\Phi(x,t,\psi,\lambda)=\partial_\psi\Phi(x,t,\psi,\lambda)=0$ 
but already for $t=0$, the determinant of the Hessian of $\Phi$
with respect to $(t,\psi)$ vanishes on $\Lambda_+$. We can solve instead (locally)
$\Phi'_{(t,\lambda)}(x,t,\psi,\lambda)=(\partial_t\Phi,\partial_\lambda\Phi)=0$. 
Namely since ${\partial^2\Phi\over\partial t\partial\lambda}=-mE\neq0$,
the implicit function theorem gives $(t,\lambda)=\bigl(t(x,\psi),\lambda(x,\psi)\bigr)$.
 
We want to keep $\lambda(x,\psi)=1$. Differentiating $\Phi'_{(t,\lambda)}=0$ along $\Lambda_+$ 
with respect to $x$ and $\psi$ we find,
using (2.6) and Hamilton equations
$$\eqalign{
&\langle \dot P,\dot X\rangle{\partial t\over\partial x}+mE{\partial\lambda\over\partial x}={}^t \dot P, \quad 
mE{\partial t\over\partial x}={}^t P\cr
&\langle \dot P,\dot X\rangle{\partial t\over\partial\psi}+mE{\partial\lambda\over\partial\psi}=-\langle\dot P,X_\psi\rangle=-a, \quad
mE{\partial t\over\partial\psi}=0
}\leqno(4.24)$$
Assume ${\partial\lambda\over\partial x}=0$ at $z_0$. This implies $\dot P={\langle\dot P,\dot X\rangle\over mE}P$, 
i.e. $\partial_xH+{1\over mE}\langle-\partial_xH,\dot X\rangle P=0$. Taking scalar product with $P\neq0$ we find
$\langle \partial_xH,P\rangle+{|P|^2\over mE}\langle-\partial_xH,\dot X\rangle P=0$, and since $z$ is a special point,
$\langle-\partial_xH,\dot X\rangle P=0$. It follows that $\partial_xH=0$ which is a contradiction ($z$ is not a residual point). 

Now we need $\lambda=\lambda(x,\psi)=1$~; since $\partial_x\lambda\neq0$, the implicit functions theorem shows that (possibly
after renumbering the coordinates) that $x_2=x_2(x_1,\psi)$. 
By second line (4.24) we have
${\partial t\over\partial\psi}=0$, and 
$-mE{\partial\lambda\over\partial\psi}=a=\langle\dot P,X_\psi\rangle$.

\noindent $\bullet$ Assume $a=0$. 

Since we have eliminated $t,\lambda$, the ``$\theta$-parameter'' is simply $\psi$, and we set 
$$\Psi(x_1,\psi)=\Phi(x_1,x_2(x_1,\psi),t_2(x_1,x_2(x_1,\psi)),\psi,\lambda(x_1,x_2(x_1,\psi),\psi)$$
By (3.5) with $N=1, n=2$, it follows that $\rank d\pi_x=2$ if $\partial^2_\psi\Psi(x_1,\psi)|_{x_1=X_1}\neq0$, 
and $\rank d\pi_x=1$ if $\partial^2_\psi\Psi(x_1,\psi)|_{x_1=X_1}=0$ ($X_1$ being evaluated at $(t,\psi)=(t_0,\psi_0)$. 
In the latter case, differentiating $\lambda=\lambda(x,\psi)=1$ gives
${\partial\lambda\over\partial \psi}+{\partial\lambda\over\partial x}{\partial x\over\partial \psi}=0$. Since ${\partial\lambda\over\partial \psi}=0$
at point $z$, (4.22) easily follows.

\noindent $\bullet$ Assume $a\neq0$. From $\lambda(x,\psi)=1$, we get $\psi=\psi(x)$ by implicit function theorem, so we have eliminated all
``$\theta$-variables'' and $\rank d\pi_x=2$. 
\smallskip
\noindent {\bf (iii)} We consider a residual point as a limit of special points, for which ${\partial \lambda\over\partial x}=0$. Since $a=0$, we have 
${\partial \lambda\over\partial \psi}={\partial \lambda\over\partial x}=0$ at $z_0$, and $\lambda(x,\psi)=1+{\cal O}(|x-X(t_0,\psi_0),\psi-\psi_0|^2)$.
Then (4.24) reduces to $mE{\partial t\over\partial x}={}^t P$ at $z_0$, which can be cast in the form $dt=\langle P,dx\rangle\neq0$. 
$\clubsuit$. 
\smallskip
For residual points Proposition 4.3 tells nothing however about $\rank d\pi_x(z_0)$.  
For instance, if $\dot P=0$, hence $\langle \dot P,\dot X\rangle=0$ and $\partial_t\Phi=\partial^2_t\Phi=0$, we could have $\partial^3_t\Phi\neq0$,
and we have a cusp described by Pearcy functions (see e.g. [DoMaNaTu1,App.2]). Alternatively we could think of Hamiltonian $p^2$
for which $\rank d\pi_x(z_0)=n$ is maximal, or of Hamiltonian $H(x,p)=p_1$ for which $\rank d\pi_x(z_0)=1$.
It tells nothing either about ordinary points, see however Lemma 4.6 below when $H$ is of the form (1.25). 
\medskip
\noindent {\it 4.4 Construction of $E_+(h)f_h$ and the commutation formula}.
\smallskip
Here we prove Theorem 1.3. First we look for a solution $v(t,x;h)$ to the Cauchy problem (1.6)
that can be expressed as an oscillatory integral $\int e^{i\Phi(x,t,\psi,\lambda)/h}b(x,t,\psi,\lambda)\,d\psi\,d\lambda$, see (1.26).

Assume for simplicity $H(x,hD_x;h)$ has no sub-principal symbol: $H_1(x,p)=0$. Then it is well known that the principal term $b=b_0$
of the amplitude restricted to $C_\Phi$, since $mE\det(P,P_\psi)$ is the (inverse) density, is of the form
$$b_0=\bigl(mE\det(P,P_\psi)\bigr)^{-1/2}a(\psi,\lambda)\leqno(4.25)$$
with $a$ independent of $t$. Since $P,P_\psi$ are linearly independent, we look for 
$$\eqalign{
b(x,&t,\psi,\lambda)=\bigl(mE\det(P,P_\psi)\bigr)^{-1/2}a(\psi,\lambda)+\langle FP(t,\psi)+GP_\psi(t,\psi),x-X(t,\psi)\rangle\cr
&+{\cal O}(|x-X(t,\psi)|^2+|\lambda-1|^2)
}$$
where we can determine functions $F=F(t,\psi,\lambda),G=G(t,\psi,\lambda)$ from second derivatives of $H$
by taking variations.
Set $\widetilde b=\langle FP+GP_\psi,x-X(t,\psi)\rangle$, it is readily seen that
$$\partial_t\widetilde b+\langle\partial_pH,\partial_x\widetilde b\rangle=0 \ \hbox{on} \ C_\Phi\leqno(4.26)$$
Let $v_h$ solves Cauchy problem (1.6), and $u_h(x)={i\over h}\int_0^\infty v(t,x;h)\,dt$ (after sticking in a cut-off $\theta_T(t)$ as in (3.16)).
We start with computing $(H-E)u_h$ and assume the
general case of $H$ homogeneous of degree $m$ in the $p$ variables, and $\det(P,P_\psi)\neq0$.
For simplicity we present the calculations as if $H$ were a differential operator, see [Du].  
We assume also the sub-principal symbol of $H$ (as a $h$-PDO) vanishes.
\smallskip
\noindent $\bullet$ Let first $x\in M$ be such that $x=X(t,\psi)$ with $X_\psi\neq0$, so that $T_z\Lambda_+$ is transverse to the vertical plane. 
We use representation (1.26).
Applying $H-E$ to (1.26), since $\Phi$ is just linear in $x$, we get first
$$e^{-i\Phi/h}(H-E)e^{i\Phi/h}b=(H(x,\partial_x\Phi)-E)b+{h\over i}\langle\partial_pH(x,\partial_x\Phi),\partial_xb\rangle+{\cal O}(h^2)$$
By (4.8) we have, integrating by parts
$$\eqalign{
(H&-E)u_h(x)=\int e^{i\Phi/h}b\,d\psi\,d\lambda|_{t=0}+\int_0^\infty dt\int e^{i\Phi/h}\bigl(\partial_tb+\langle\partial_pH(x,\partial_x\Phi),
\partial_xb\rangle\bigr)\,d\psi\,d\lambda+\cr
&{i\over h}\int_0^\infty dt\int e^{i\Phi/h}{\cal O}(|x-X(t,\psi),\lambda-1|^2)b\,d\psi\,d\lambda+{\cal O}(h^2)
}\leqno(4.27)$$
To the second integral we apply asymptotic stationary phase [H\"o,Theorem 7.7.5]; denote
$c(x,t,\psi,\lambda)={\cal O}(|x-X(t,\psi),\lambda-1|^2)b(x,t,\psi,\lambda)$,
and by $\Phi_c$ the critical value of $(t,\psi,\lambda)\mapsto\Phi(x,t,\psi,\lambda)$ with non critical degenerate point 
$\eta(x)=(t=t(x),\psi=\psi(x),\lambda=1)$ (see (4.17)) we have
$$\eqalign{
{i\over h}&\int_0^\infty dt\int e^{i\Phi/h}
e^{i\Phi/h}c(x,t,\psi,\lambda;h)\,d\psi\,d\lambda=\cr
&e^{i\Phi_c/h}\bigl(\det(\Phi''/2i\pi h)\bigr)^{-1/2}\bigl(c(x,\eta(x);h)+{h\over i}\langle \Phi''(x,\eta(x))^{-1}D,D\rangle
c(x,\eta(x))+{\cal O}(h^2)
}$$
where $D$ denotes the gradient with respect to the 3 variables $t,\psi,\lambda$ (of course we still assumed $n=2$). We have  
$c(x,\eta(x);h)=0$, but the next term $\langle \Phi''(x,\eta(x))^{-1}D,D\rangle
c(x,\eta(x))$ may not vanish because of the partial derivative ${\partial^2c\over\partial\lambda^2}$, as shows (4.17). So
$${i\over h}\int_0^\infty dt\int e^{i\Phi/h}
e^{i\Phi/h}c(x,t,\psi,\lambda;h)\,d\psi\,d\lambda={\cal O}(h^{3/2})$$
We consider next the first integral in (4.27). Because of (4.25) and (4.26) which implies 
$$\partial_tb+\langle\partial_pH(x,\partial_x\Phi),
\partial_xb\rangle={\cal O}(|x-X(t,\psi)|)$$
we apply asymptotic stationary phase as before and obtain
$$\int_0^\infty dt\int e^{i\Phi/h}\bigl(\partial_tb+\langle\partial_pH(x,\partial_x\Phi),
\partial_xb\rangle\bigr)\,d\psi\,d\lambda={\cal O}(h^{5/2})$$
Collecting these estimates in (4.27) yields 
$$(H-E)u_h(x)=\int e^{i\Phi/h}b\,d\psi\,d\lambda|_{t=0}+{\cal O}(h^{3/2})$$
Since we can choose $b_0$ in (4.25) to be equal to the amplitude defining $f_h$, the RHS is just $f_h$ mod ${\cal O}(h^{3/2})$.
\smallskip
\noindent {\it Remark 4.4}: Note the loss of $h^{1/2}$ with respect to the remainder term ${\cal O}(h^2)$ when solving the 
homogeneous equation $(H-E)u_h=0$, see the discussion after (3.13). 
\medskip
\noindent $\bullet$ Take next $x\in M$ near $X(t,\psi)$ with $X_\psi=0$, by the discussion after (4.21), 
up to a permutation of $x_1$ and $x_2$, we may consider in the mixed representation $H(-hD_{p_1},x_2,\xi_1,hD_{x_2};h)$. We try as new phase function
$$\Phi_1(x_2,p_1,t,\psi,\lambda)=mEt-\lambda X_1(t,\psi)p_1+\lambda P_2(t,\psi)\bigl(x_2-X_2(t,\psi)\bigr)$$
so that the eikonal equation reads
$$\partial_t\Phi_1+H(-\partial_{p_1}\Phi_1,x_2,p_1,\partial_{x_2}\Phi_1))={\cal O}(|x_2-X_2(t,\psi),p_1-P_1(t,\psi),\lambda-1|^2)$$
Transport equations are derived similarly. Using again (1.26) we can present
$$(H(-hD_{p_1},x_2,\xi_1,hD_{x_2};h)-E){\cal F}^h_{x_1\to p_1}u_h(p_1,x_2)$$
in the form
$$(H(-hD_{p_1},x_2,\xi_1,hD_{x_2};h)-E){i\over h}\int_0^\infty\,dt\int e^{i\Phi(x_2,p_1,t,\psi,\lambda)/h}b(x_2,p_1,t,\psi,\lambda)\,d\psi\,d\lambda$$
which we compute as before by asymptotic stationary phase. Theorem 1.3 easily follows $\clubsuit$
\medskip
\noindent {\it 4.5 Reduced parametrizations of $\Lambda_+$ in case of the ``conformal metric'', $n=2$}.
\smallskip
In case of the conformal metric we can make the results more precise (at least for $n=2$), due to fact that $\dot X$ is parallel to $P$.
First information is related with the density.
By Proposition 4.1, $\Phi$ is a non-degenerate phase function parametrizing $\Lambda_+$ iff $\det (P,P_\psi)\neq0$, see (4.10). 
This certainly holds for small $t$.
We want to allow for larger values of $t$ (the far field). We have no direct proof that (4.10)
is valid everywhere on $\Lambda_+$. See however [DoMaNaTu1], Example 6, in case case $m=1$, and $\rho$ is radially symmetric. 
In general, this property is related with parametrization of Lagrangian submanifolds, see [H\"o,Thm 21.2.16]. In case of the conformal metric, 
Lemma A.2 readily implies~:
\smallskip
\noindent {\bf Proposition 4.4}: {\it Let $H(x,p)$ be as in (1.25), $n=2$. Then at least near focal points,
representation (4.5) defines a non degenerate phase function parametrizing $\Lambda_+^E$, and
the (inverse) density on $\Lambda_+^E$ is $mE\det (P,P_\psi)\neq0$. }
\smallskip
This holds also when $\Lambda_0$ is the ``cylinder'', see Proposition 5.4 below. 
\medskip
Next information is related to eliminating extra ``$\theta$-variables'' in the phase function and determining the rank of $\pi_x:\Lambda_+^E\to M$.
For simplicity we consider only the position representation of $u_h$, i.e. the case $X_\psi\neq0$. 
As in Proposition 4.3, we proceed to find the critical value of $t\mapsto \Phi(x,t,\psi,\lambda)$ when $z(t)=\bigl(X(t,\psi),P(t,\psi)\bigr)$ 
is an ordinary point (which is equivalent to $\langle\dot X,\dot P\rangle\neq0$ in case of Hamiltonian (1.25). 
In a simple scenario there would be at most one special or residual point on each bicharacteristic.  
Definition 1.5 provides such a scenario. Recall
$\partial_t\Phi(x,t,\psi,\lambda=1)=\langle\dot P(t,\psi),x-X(t,\psi)\rangle=0$
This holds on $\Lambda_+^E$, i.e. for $x=X(t,\psi)$. Taking second derivative at critical point gives
$$\partial^2_t\Phi(X(t,\psi),t,\psi,\lambda=1)=-{m|P(t,\psi)|^{2m-2}\over\rho(X(t,\psi))^3}\langle\nabla\rho(X(t,\psi)),P(t,\psi)\rangle
\leqno(4.30)$$
so we have to take into account the set of $\psi$ such that $\langle\nabla\rho(X(t,\psi)),P(t,\psi)\rangle=0$, i.e.
of the special or residual points. Consider
$f(t,\psi)=\langle\nabla\rho(X(t,\psi)),P(t,\psi)\rangle$, so that $f=0$ iff $(X(t,\psi)),P(t,\psi)$ is special or residual. 
Using Hamilton equations, we find 
$$\partial_tf(t,\psi)={m|P(t,\psi)|^{m-2}\over\rho(X(t,\psi))}\bigl[ \langle\nabla^2\rho(X(t,\psi)\cdot P(t,\psi),P(t,\psi)\rangle
+{|\nabla\rho(X(t,\psi))|^2\over m\rho(X(t,\psi))}|P(t,\psi)|^2\bigr]\leqno(4.31)$$
Let $z(s,\psi)$ be a special (or residual) point for some $s\geq0$, then whenever 
$\partial_t\Phi(t,x,\theta)=0$ at some $t>s$ (this occurs when the bicharacteristic $t\mapsto z(t)$ projects again on $x$),
$z(t)$ is no longer special (or residual). This holds under
assumption (1.35), namely $\partial_tf(t,\psi)>0$, and $t\mapsto\partial^2_t\Phi(X(t,\psi),t,\psi,\lambda=1)$ is strictly decaying on $\Lambda_+^E$. 
\smallskip
\noindent $\bullet$ {\it Ordinary critical points}. They correspond to non degenerate critical points of $t\mapsto\Phi(x,t,\psi,\lambda=1)$.
\medskip
\noindent {\bf Lemma 4.6}: {\it Assume (1.35), $n=2$ (no condition on $\nabla\rho$ is required here). Let
$I_{\psi}=\{t: z(t,\psi)\notin{\cal S}(\Lambda_+^E)\}$ (we have already evaluated $\lambda=1$)
Then $I_\psi$ is an interval, and
$$\forall t\in I_\psi, \ x=X(t,\psi) \ \Longleftrightarrow \ t=t_1(x,\psi) \quad \hbox{on} \quad C_\Phi\leqno(4.33)$$ 
where $t_1$ is a smooth function. Moreover $\pi_x:{\Lambda_+^E}\to M$ at every ordinary critical point has same rank as the symmetric matrix
(4.23) i.e. $d\pi_x$ has rank 1 ($\langle P_\psi,X_\psi\rangle=0$) or 2 ($\langle P_\psi,X_\psi\rangle\neq0$). }
\smallskip
\noindent {\it Proof}: 
Note that when $t=0$, $f(0,\psi)=\langle\nabla\rho(0),P(\psi)\rangle$. So when $f(0,\psi)>0$,  
0 is non-degenerate critical point of $t\mapsto\Phi(x,t,\psi,\lambda)$, and the implicit function
theorem shows that (4.33) holds. Since $t\mapsto f(t,\psi)$ is increasing,
this holds for all $t$ in the maximal interval of definition of the integral curve starting at $(0,P(\psi))$. 
When $f(0,\psi)<0$ instead, (4.33) holds on an interval ending at some $s$ such that $f(s,\psi)=0$. 
Let us compute the rank of $\pi_x:\Lambda_+^E\to M$ at an ordinary point. Let 
$U=\{(t,\psi): t\in I_\psi, \ z(t)=(X(t,\psi),P(t,\psi))\notin{\cal S}(\Lambda_t)\}$, then the same computation as in Proposition 4.3 
shows that $U$ is a canonical chart rank 1 or 2, which gives the Lemma. $\clubsuit$
\medskip
\noindent $\bullet$ {\it Special and residual critical points}. Near the end point $s$ of $I_\psi$ 
we can solve (locally) as in Proposition 4.3,
$\partial_t\Phi(x,t,\psi,\lambda)=\partial_\lambda\Phi(x,t,\psi,\lambda)=0$ 
which gives $t=t(x,\psi)$ and $\lambda=\lambda(x,\psi)$. Namely,  
the Hessian of $\Phi$ with respect to $(t,\lambda)$ at $(s,1)$  has determinant
$-\bigl({\partial^2\Phi\over\partial t\partial\lambda}\bigr)^2=-(mE)^2<0$ on $\Lambda_+^E$. 
So if $(X(s,\psi),P(s,\psi)\bigr)$ is a special point 
then $(s,1)$ is a non degenerate point of $(t,\lambda)\mapsto \Phi(x,t,\psi,\lambda)$.  
Integrating Hamilton equations also for $t<0$ gives the Lagrangian manifold $\Lambda_-^E\cup\Lambda_+^E$.
So there is no loss of generality in assuming the special point is at $s=0$. The following Lemma strengthens Proposition 4.3 in case $X_\psi\neq0$. 
\medskip
\noindent {\bf Lemma 4.7}: {\it Let $n=2$ and $H$ be as in (1.25). Assume $\det(P,P_\psi)\neq0$.

(i) Assume $z(s)=(X(s,\psi),P(s,\psi))\in\Lambda_+^E$ be a special point (hence $\nabla\rho(x(s))\neq0$). 
If $X_\psi=0$, then $\rank d\pi_x(z(s))=1$ as in Proposition 4.3 (i).
If  $X_\psi\neq0$, then $ac\neq0$ so that $\rank d\pi_x(z)=2$. Near $z(s)$, $\Lambda_+^E$ is given by  $t=t(x)$, $\psi=\psi(x)$, and
${\partial t\over\partial x}\neq 0$, ${\partial\psi\over\partial x}\neq0$. The constraint $\lambda=1$ takes the form
$${\partial\lambda\over\partial\psi}+\big|{\partial\psi\over\partial x}\big|^{-2}
\langle {\partial\lambda\over\partial x},{\partial\psi\over\partial x}\rangle=0\leqno(4.41)$$

(ii) Assume $z(s)=(X(s,\psi),P(s,\psi))\in\Lambda_+^E$ be a residual point (i.e. $\nabla\rho(x(s))=0$). 
If $X_\psi\neq0$, then $c\neq0$ and $\rank d\pi_x(z)=2$. }
\smallskip
\noindent {\it Proof}: As in Proposition 4.3, the relations $\partial_t\Phi=\partial_\lambda\Phi=0$
being given by $(t,\lambda)=\bigl(t(x,\psi),\lambda(x,\psi)\bigr)$ we use (4.24). 
Since $X_\psi\neq0$, $\langle P,X_\psi\rangle=0$ and $\det(P,P_\psi)>0$ we have $c\neq0$.
\smallskip 
\noindent {\bf (i)} By the same geometric argument we have $\dot P\neq0$ by (2.2), and since $\langle\nabla\rho,P\rangle=0$, 
the relation $a=\langle\dot P,X_\psi\rangle=0$
would contradict $X_\psi\neq0$. 
So by second line (4.24), $-mE\partial_\psi\lambda=\langle\dot P,X_\psi\rangle\neq0$, or $\partial_\psi\lambda\neq0$.    
Now we need $\lambda=\lambda(x,\psi)=1$~; since $\partial_\psi\lambda\neq0$, the implicit functions theorem shows that $\psi=\psi(x)$. 
Then we have
$$x=X(t,\psi) \ \Longleftrightarrow \ t=t(x), \psi=\psi(x) \ \hbox{on} \ 
 \partial_t\Phi(x,t,\psi,\lambda=1)=\partial_\lambda\Phi(x,t,\psi,\lambda=1)=0\leqno(4.42)$$
Differentiating $\lambda=\lambda(x,\psi)$ gives 
$${\partial\lambda\over\partial x}+{\partial\lambda\over\partial\psi}{\partial\psi\over\partial x}=0$$
and together with the first equation (4.24)
$$\langle P_\psi,\dot X\rangle {}^t\bigl({\partial\psi\over\partial x}\bigr)=
mE\dot P-\langle\dot P,\dot X\rangle P$$
Let us show that ${\partial\psi\over\partial x}\neq0$. Otherwise, we would have 
$\langle P,\dot X\rangle\dot P=\langle \dot P,\dot X\rangle P$, and since we know that 
$\langle\nabla \rho(X(s,\psi)),P(s,\psi)\rangle=0$, this would contradict the fact that
$\dot P$ is parallel to $\nabla\rho(X(s,\psi)$. Moreover ${\partial t\over\partial\psi}=0$, ${\partial t\over\partial x}={1\over mH}{}^t P\neq0$,
which readily gives (4.41).
To compute the rank of $\pi_x$ at a special point, we are left to compute second derivative of the critical value, namely
$-\partial^2_\psi\Psi=ac\neq0$, so we conclude as in Proposition 4.3 that $\rank \pi|_{\Lambda_+}$ is 2. 
\smallskip
\noindent {\bf (ii)} Thinking of a residual point as the limit of special points, (4.24) shows that 
${\partial\lambda\over\partial x}={\partial\lambda\over\partial\psi}=0$, and (4.24) reduces to $mE{\partial t\over\partial x}={}^t P, 
{\partial t\over\partial \psi}=0$.
Note that on $C_\Phi$, (4.17) gives $\det\Hess_{(t,\psi,\lambda)}\Phi=(mE)^2\langle P_\psi,X_\psi\rangle\neq0$ if $X_\psi\neq0$, 
so by the implicit functions theorem
$$(t,\psi,\lambda)=\bigl(t_2(x),\psi_2(x),\lambda_2(x)\bigr)$$ 
Let us check again that $\lambda_2(x)=1$~: 
Differentiating $\Phi'_{t,\psi,\lambda}(x,t,\psi,\lambda)=0$ with respect to $x,\lambda$ gives the triangular system
$$\eqalign{
&\langle P_\psi,\dot X\rangle{}^t\bigl({\partial t_2\over\partial x}\bigr)+\langle P_\psi,X_\psi\rangle
{}^t\bigl({\partial \psi_2\over\partial x}\bigr)=P_\psi\cr
&\langle \dot P,\dot X\rangle{}^t\bigl({\partial t_2\over\partial x}\bigr)+\langle \dot P,X_\psi\rangle
{}^t\bigl({\partial \psi_2\over\partial x}\bigr)+mH{}^t\bigl({\partial\lambda_2\over\partial x}\bigr)=\dot P\cr
&mH{}^t\bigl({\partial t_2\over\partial x}\bigr)=P
}$$  
Since $\dot P=0$, using $a=0$ (see (4.16)) this reduces to
$$c{}^t\bigl({\partial \psi_2\over\partial x}\bigr)=P_\psi, \quad  {\partial\lambda_2\over\partial x}=0, 
\quad mH{}^t\bigl({\partial t_2\over\partial x}\bigr)=P$$
and in particular $\lambda=\lambda_2(x)=1$. There are no ``$\theta$-parameters'' left and so $\rank d\pi_x(z(s)=2$. 
Then $\det(P,P_\psi)\neq0$ implies $\det\bigl({\partial \psi_2\over\partial x},{\partial t_2\over\partial x}\bigr)\neq0$. 
$\clubsuit$
\smallskip
Note that if $X_\psi=0$ at a focal point of $\Lambda_+^E$, then $P_\psi\neq0$ (otherwise this would violate property (3) of Proposition A.1).

From Lemma 4.7, the set of focal points which are also special points is ${\cal S}(\Lambda_t)\cap{\cal F}(\Lambda_t)=
\{\psi: \langle P_\psi,X_\psi\rangle=0\}$.
In Example 2.3, we find 
$\langle P_\psi,X_\psi\rangle={1\over2}\sinh (2f(t))>$0 for all $t>0$ and vanishes at $t=0$.

Now to find the canonical charts for the phase functions,  
we use a connectedness argument. Assume (1.35), and let $s$ be the supremum of $I_\psi$, we have $f(s,\psi)=0$. Since $G(\rho)(s,\psi)>0$, 
we have $f(t,\psi)>0$
for all $t>s$, so all points $z(t)=(X(t,\psi),P(t,\psi))$ for $t>s$ are ordinary points.  
So far we proved (except for the case $X_\psi=0$ which can be handled similarly by replacing $\Phi$ by its Legendre transformation):
\smallskip
\noindent {\bf Proposition 4.8}: {\it Let $H(x,p)={|p|^m\over\rho(x)}$, $n=2$. 
Then there exists globally a smooth (parametric, i.e. defined by a non-degenerate phase function) solution $\Phi$ of HJ equation $H(x,\partial_x\Phi)=E$.
Let $C_\Phi=\{(t,x,\psi), t>0, \psi\in{\bf R}: \partial_t\Phi=\partial_\lambda\Phi=\partial_\psi\Phi=0\}$. Then the embedding
$$\iota_\Phi:C_\Phi\to T^*{\bf R}^2, \ (t,x,\psi,\lambda)\mapsto(x,\partial_x\Phi(x,t,\psi,\lambda=1))$$
such that $\iota_\Phi(C_\Phi)\subset\Lambda_+$ consists in charts of rank 1 or 2
[the rank is never 0 since $p\neq0$ in the energy shell $H=E$]. Under the defocussing condition (1.35) these charts can intersect the 
set of special points ${\cal S}$ only along a line. }
\medskip
\noindent{\it Remark 4.9}: The canonical charts in $\Lambda_+$ where $\Phi=\Phi(x)$, i.e. of WKB type are of course of maximal rank 2, in particular
there is a WKB solution near a special point $z$ such that $\langle X_\psi,P_\psi\rangle\neq0$. 
\bigskip
\noindent {\bf 5. $f_h$ is supported microlocally on ``Bessel cylinder''}
\smallskip
We recall $\Lambda$ from (1.3). 
When $n=2$ this is the wave-front set of Bessel function $f_h(x;h)=J_0({|x|\over h})$.
In case of the vertical plane, $\Lambda_+^E$ was parametrized by $(\psi,t)$, $\psi\in{\bf R}^{n-1}$ parametrizing
$\partial\Lambda_+^E$. This is no longer the case 
for Bessel cylinder, 
since $\varphi$ and $\psi$ are not independent variables on the energy surfaces $H=E-\tau$.
Because the eikonal on $\Lambda_+^E$
takes the form $S=\varphi+mEt$ (see (2.10)) there is indeed a ``pairing''
between $t$ and $\varphi$, we shall call the ``twin variables''. The condition $m=1$ shows to be necessary.
Moreover, glancing points naturally occur on $\Lambda_0$, and we will essentially restrict to the simplest case where the Hamiltonian
is of the form (1.25) with $\rho$ radially symmetric. 
The set of glancing points then reduces to $\{x=0\}$, and the situation is quite similar to this of the vertical plane.  
\medskip
\noindent {\it 5.1 Non degeneracy conditions}
\smallskip
Let us check first the Lagrangian intersection. 
\smallskip
\noindent {\bf Definition 5.1}: {\it The point $z\in\Lambda_0$ is called {\it glancing} if
$v_H(z)\in T_z\Lambda_0$. We denote by ${\cal G}(\Lambda_0)$ the set of glancing points on $\Lambda_0$. }

So $(\Lambda_0,\Lambda_+^E)$ is an intersecting pair whenever ${\cal G}(\Lambda_0)=\emptyset$.
Hamiltonian flow preserves the set of glancing points, i.e. $\exp tv_H({\cal G}(\Lambda_0))={\cal G}(\Lambda_t)$ for all $t>0$.  
\smallskip
\noindent {\bf Proposition 5.2}: {\it Let $H\in C^\infty(T^*M)$ 
and $\Lambda_0$ be Bessel cylinder
$\Lambda_0=\{x=X(\varphi,\psi)=\varphi\omega(\psi), p=P(\varphi,\psi)=\omega(\psi), \varphi\in{\bf R}\}$.
We set ${\cal H}(\varphi,\psi)=H|_{\Lambda_0}$. 
With the notation above, $z=(x,p)\in\Lambda_0$ is a glancing iff
$$\nabla{\cal H}(\varphi,\psi)=0\leqno(5.1)$$ 
In particular, let $H$ be homogeneous of degree $m$ with respect to $p$, then $z=(x,p)\in\Lambda_0$ is a glancing point at energy $E$ iff
$$\partial_p H+\varphi\partial_x H=mE\omega(\psi),  
\quad \langle-\partial_x H,\omega(\psi)\rangle=0, \quad H(z)=E\leqno(5.2)$$
and $z$ is a special point, in the sense of Definition 1.4. }
\smallskip
\noindent {\it Proof}: We complete $\omega(\psi)$ in ${\bf S}^{n-1}$ into a (direct) orthonormal basis 
$\omega^\perp(\psi)=\bigl(\omega_1(\psi),\cdots,\omega_{n-1}(\psi)\bigr)$ of ${\bf R}^n$, and denote by 
$\omega^\perp(\psi)\delta\psi=\omega_1(\psi)\delta\psi_1+\cdots+\omega_{n-1}(\psi)\delta\psi_{n-1}$ a section of $T{\bf S}^{n-1}$, 
$\delta\psi_j\in{\bf R}$.
The tangent space $T_z\Lambda_0$ has the parametric equations
$$\delta X=\omega(\psi)\delta\varphi+\varphi\omega^\perp(\psi)\delta\psi, \quad \delta P=\omega^\perp(\psi)\delta\psi, 
\quad (\delta\varphi,\delta\psi)\in{\bf R}^n$$
so $v_H\in T_z\Lambda_0$ iff there exist $(\delta\varphi,\delta\psi)$ such that
$$\partial_pH=\omega(\psi)\delta\varphi+\varphi\omega^\perp(\psi)\delta\psi, \quad -\partial_xH=\omega^\perp(\psi)\delta\psi$$
Taking scalar products with $\omega(\psi),\omega^\perp(\psi)$, we get
$$\delta\varphi=\langle\partial_p H,P(\psi)\rangle, \ \delta\psi=\langle-\partial_xH,\omega^\perp(\psi)\rangle$$
Since $(\omega(\psi),\omega^\perp(\psi))$ form a basis of ${\bf R}^n$, relations
$$\partial_p H+\varphi\partial_x H=\langle\partial_p H,P(\psi)\rangle\omega(\psi),   
\quad \langle-\partial_x H,\omega(\psi)\rangle=0\leqno(5.3)$$
are necessary and sufficient
for $v_H\in T_z\Lambda_0$. 

On the other hand
$$\nabla{\cal H}(\varphi,\psi)=\bigl(\langle\partial_xH,\omega^\perp(\psi)\rangle, \varphi\langle\partial_xH,\omega^\perp(\psi)\rangle+
\langle\partial_pH,\omega^\perp(\psi)\rangle\bigr)$$
so (5.3) readily gives (5.1).

Now, if $H$ is positively homogeneous of degree $m$ with respect to $p$,  using Euler identity, we get
$\delta\varphi=\langle\partial_p H,P(\psi)\rangle=mH$,
and (5.2) holds iff
for $v_H\in T_z\Lambda_+^E\cap T_z\Lambda_0$ when $H=E$. $\clubsuit$
\medskip

\noindent {\it Example 5.1} (Helmholtz equation with constant coefficients): 
When $H=p^2$, all points are glancing at energy 1.
Start from Helmholtz equation
$(-h^2\Delta-1)f_h(x)=0$ in ${\bf R}^2$. A radially symmetric solution is given by
$$f_h(x)=(2\pi/h)^{1/2}J_0\bigl({|x| / h}\bigr)$$
and microlocalized on Bessel cylinder.
In turn, all points of $\Lambda_0$ are glancing for $-h^2\Delta$ at energy 1, and a radially symmetric solution of Helmholtz equation
$$(-h^2\Delta-1)u_h(x)=f_h(x)$$
is given by the semi-classical Lagrangian distribution
$$u_h(r) = -J_1\big({|x|\over h}\bigr) {|x|\over2h}$$

\medskip
\noindent {\it Example 5.2} (Helmholtz equation with variable coefficients):  
When $H={|p|^m\over\rho(x)}$, $z(0)$ is a glancing point iff 
$$ \hbox{either}: \ \varphi\neq0 \ \hbox{and} \ \nabla\rho=0, \ \hbox{or}: \ \varphi=0 \ \hbox{and} \ 
\langle\nabla\rho(0),\omega(\psi)\rangle=0\leqno(5.4)$$
Second condition means that if $z(0)=(0,\omega(\psi))$ is a special point. Assuming the defocussing condition
(1.35) it follows that if $z(0)$
is a glancing point, $z(t)$ will be glancing but never special at later $t>0$. 
\medskip

We say that $z\in\Lambda_0$ is a {\it non degenerate glancing point} if ${\cal H}$ has a non-degenerate critical point at $(\varphi,\psi)$.
It follows from Proposition 5.2 
that the set of corresponding energies $E$ is discrete. Such a non-transversality is called a {\it kiss} in [ElGr]. 
\medskip
\noindent {\it Example 5.3} Let $n=2$, $m=1$, $H(z)={|p|\over\rho(x)}$, with
$\rho(x)={1\over2}(1+(x-x_0)^2)$. If $x_0=\varphi\omega(\psi)\neq0$, we have
$\rho^4(x_0)\det \nabla^2 (H|_{\Lambda_0})=\varphi^2$, $\rho^2(x_0)\Tr \nabla^2 (H|_{\Lambda_0})=-(1+\varphi^2)$.
Critical energy is given by $E=H(z)={|p|\over\rho(x)}$, i.e. $E_0=1/\rho(x_0)$
\medskip
Let us first describe $\partial\Lambda_+(\tau)$ in case (1.25).
The intersection of $\Lambda$ with the energy surface $H=E-\tau$, is given by the implicit equation
$H(\varphi\omega(\psi),\omega(\psi))=E-\tau$, which usually defines a smooth $(n-1)$ dimensional isotropic submanifold.

Assume that $\Lambda_+(\tau)$ (or simply $\Lambda_+^E$) can be parametrized by $(t,\psi)$ as in Sect.4 when $\Lambda_0$ is the vertical plane. 
Recall from Proposition A.1 that $\langle\dot X,P_\psi\rangle=\langle\dot P,X_\psi\rangle$.
So taking the limit $t\to0_+$ readily implies 
$$\varphi\langle\nabla\rho(\varphi\omega(\psi),\omega^\perp(\psi)\rangle=0\leqno(5.6)$$
on $\partial\Lambda_+(\tau)$. This shows that $\rho$ is necessarily radial symmetric.

We ignore herafter glancing points. So let $z_0=\bigl(X(\varphi_0,\psi_0),P(\varphi_0,\psi_0)\bigr)\in V\subset\Lambda_0\cap\Sigma_E$ not a glancing point.
Let us summarize our discussion so far in the following
\medskip
\noindent {\bf Lemma 5.3}: 
{\it Let ${\cal H}=H_{\Lambda_0}$ as in Proposition 5.2, and $\alpha_0$ be one of the $n$ coordinates $\alpha=(\varphi,\psi)$ such that
$\partial_{\alpha_0}{\cal H}(\varphi_0,\psi_0)\neq0$. Then
$(\varphi,\psi)\mapsto\tau=\tau(\varphi,\psi)$ is a local submersion, i.e.
$\varphi=\varphi(\tau,\psi)$, or $\psi_1=\psi_1(\tau,\varphi,\psi')$, where $\psi=(\psi_1,\psi')$, possibly renumbering the coordinates. 
Away from glancing points, either condition is met in some local chart $\psi\in V_0$, or
$(\varphi,\psi')\in V_1$.

In the case (1.25) with a radially symmetric conformal metric, this means that
$\varphi\mapsto\tau=\tau(\varphi)$ is a local diffeomorphism.}

\smallskip
Next we discuss the properties of Hamiltonian flow issued from $\Lambda_0$ using eikonal coordinates as in Sect.4.1. 
Consider Hamiltonian $\tau+H(x,p)$ on $T^*({\bf R}^{n}\times{\bf R}_+)$ and recall from (1.15)
the Lagrangian manifold $\Lambda_+=\bigcup_{\tau}\Lambda_+(\tau)$ in the extended phase-space, i.e.
$$\eqalign{
\Lambda_+=\{&(x,p;t,\tau): \tau+H(x,p)=E, z(t)=\bigl(X(t,\varphi,\psi,\tau),P(t,\varphi,\psi,\tau)\bigr),\cr 
&z(0)\in \partial\Lambda_+(\tau), t\geq0\}\subset T^*(M\times{\bf R}_+)
}\leqno(5.10)$$
which is a Lagrangian embedding if we take $t,\tau$ small enough.
Let $\widetilde X={t\choose X(t,\varphi,\psi,\tau)}$, and $\widetilde P={\tau\choose P(t,\varphi,\psi,\tau)}$
in the extended phase-space.
The action on $\Lambda_+$ is of the form
$$\langle\widetilde P,d\widetilde X\rangle=\langle P(t,\varphi,\psi,\tau),dX(t,\varphi,\psi,\tau)\rangle+\tau\,dt\leqno(5.11)$$ 
Here $\varphi$ is still considered as a variable.

Let us carry to Bessel cylinder the construction of ${\cal M}$ as in Sect.4. 

With the notations of (4.1)-(4.2), let $\widetilde n=n+1$, $k=2$, $\widetilde\psi=\psi$, $\widetilde\phi=(t,\varphi)$.
We call $(t,\varphi)$ the {\it twin variables}.

Recall $S(t,\varphi,\psi,\tau)=\varphi+m(E-\tau)t+S_0$ from (2.10)
where we can set $S_0=0$. In the extended
phase-space, we set $\widetilde S(t,\varphi,\psi,\tau)=S(t,\varphi,\psi,\tau)+t\tau$,
one has also to differentiate with respect to $\tau$, so that $dS+\tau\,dt+t\,d\tau=\langle P,dX\rangle+\tau\,dt$, or
$$d\varphi+m(E-\tau)\,dt+(1-m)t\,d\tau=\langle P,\dot X\rangle\,dt+X_\varphi\,d\varphi+X_\psi\,d\psi+X_\tau\,d\tau$$
Since $\tau=\tau(\varphi,\psi)$, we have $d\tau=\tau_\varphi\,d\varphi+\tau_\psi\,d\psi$
so we get by identification (still with the notation $y_x={\partial y\over \partial x}$)
$$\eqalign{
&\langle P,X_\varphi\rangle+\tau_\varphi\langle P,X_\tau\rangle=1+(1-m)t\tau_\varphi\cr
&\langle P,X_\psi\rangle+\tau_\psi\langle P,X_\tau\rangle=(1-m)t\tau_\psi\cr
&\langle P,\dot X\rangle=mH
}\leqno(5.12)$$
We look for a ``left inverse'' of
$$\bigl(\dot{\widetilde X},\widetilde X_{\varphi}\bigr)=\pmatrix{1&0\cr \dot X&X_\varphi+X_\tau\tau_\varphi}$$ 
We try $\widetilde\Pi_1=\pmatrix{1&-mH\cr 0&P}$, which gives, using (5.12)
$$\widetilde\Pi_1^*\pmatrix{1&0\cr \dot X&X_\varphi+X_\tau\tau_\varphi}=\pmatrix{1&0\cr 0&1+(1-m)t\tau_\varphi}$$
so we choose 
$$\widetilde\Pi=\widetilde\Pi_1\pmatrix{1&0\cr 0&{1\over1-(m-1)t\tau_\varphi}}= 
\pmatrix{1&-\alpha mH\cr 0&\alpha P}, \ \alpha=\alpha(t,\varphi,\psi)=\bigl(1+(1-m)\tau_\varphi\bigr)^{-1} \leqno(5.13)$$
which is well defined if $1-(m-1)t\tau_\varphi>0$, which is granted for $m=1$. 
By [DoNaSh,Lemma 6] (which doesn't assume $\widetilde k$ to be the rank of $\widetilde\Lambda_+$), we know in particular that the equation 
$$\widetilde \Pi^*{0\choose x-X(t,\varphi,\psi,\tau)}=0 \Longleftrightarrow \langle P(t,\varphi,\psi,\tau),x-X(t,\varphi,\psi,\tau)\rangle=0\leqno(5.14)$$
has a unique solution $t=t_1(x,\psi,\tau)$, $\varphi=\varphi_1(x,\psi,\tau)$.

Thus as in [DoNaSh,Lemma 6] 
we can prove the following result, at least for small $t>0$. Recall Lemma 5.3. For short, we restrict to the case $\varphi=\varphi(\tau,\psi)$,
the case $\psi_1=\psi(\tau,\varphi,\psi')$ being similar. 
\smallskip
\noindent {\bf Proposition 5.4}: {\it Let us assume  $\varphi=\varphi(\tau,\psi)$. Then~:

(i) There is a open set $U\subset{\bf R}^{n+1}$, $(t,\varphi(\tau,\psi),\psi)\in U$, $\psi\in V_0\subset{\bf R}^{n-1}$ such that 
$$\kappa: U\mapsto{\bf R}^{n+1}_{t,x}\times V_0, \quad (t,\varphi(\tau,\psi),\psi)\mapsto
\bigl(\widetilde X\bigl(t,\varphi(\tau,\psi),\psi\bigr),\psi\bigr)\leqno(5.15)$$
is an embedding.

(ii) There is a neighbhd $V$ of $C=\kappa(U)$, such that 
$$\langle P(t,\varphi(\tau,\psi),\psi),x-X(t,\varphi(\tau,\psi),\psi)\rangle=0, \ (t,x,\psi)\in V\leqno(5.16)$$ 
has a unique solution $t=t_1(x,\psi)$, $\varphi=\varphi_1(x,\psi)$
satisfying the condition}
$$(t,x,\psi)\in C \ \Longrightarrow \ x=X(t_1(x,\psi),\varphi_1(x,\psi),\psi)\leqno(5.17)$$
\smallskip
Consider now
${\cal M}(t,\varphi,\psi,\tau)=\bigl(\widetilde \Pi;\widetilde P_{\widetilde\psi}-\widetilde P_{\widetilde\phi}\widetilde\Pi^*\widetilde X_{\widetilde\psi}\bigr)$
as in Sect.4.1. 
In [DoNaSh, Lemma 9], ${\cal M}$ is constructed  in a canonical chart of rank  $\widetilde k$ 
for the Lagrangian embedding (without boundary) $\iota:\Lambda\to T^*M$, 
so that ${\cal M}$ becomes an invertible matrix.  
Here we consider instead the expression of ${\cal M}$ as an Ansatz, and check that $\det{\cal M}$ defines the (inverse) density on $\Lambda_+$. 
With $\alpha$ as in (5.13) we have
$${\cal M}=\pmatrix{1&-mH\alpha&\tau_\psi-\alpha(m-1)\tau_\varphi\tau_\psi t\cr
0&\alpha P&P_\psi+\tau_\psi P_\tau+\alpha(m-1)\tau_\psi t(P_\varphi+\tau_\varphi P_\tau)}$$
so that 
$$\det {\cal M}=\alpha(t,\varphi,\psi)\det\bigl(P,P_\psi+\alpha\tau_\psi(P_\tau+(m-1)tP_\varphi)\bigr)\leqno(5.19)$$
Note that when $m=1$ (which implies that $\tau+H(x,p)$ is homogeneous of degree 1 as an Hamiltonian on $T^*(M\times{\bf R})$), 
$\det{\cal M}=\det(P, d_\psi P)$, where 
$$d_\psi P=P_\psi+\tau_\psi P_\tau\leqno(5.20)$$ 
On the other hand when $\tau_\psi=0$,
$\det {\cal M}=\alpha(t,\varphi,\psi)\det(P,P_\psi)$ is non zero for small $t$.

By the discussion above, this holds for Hamiltonian (1.25) when $\rho$ is radially symmetric, so    
for $m=1$ we recover the same (inverse) density on $\Lambda_+$ as in Sect.4. Moreover
$\det {\cal M}=\det(P,P_\psi)>0$ with the condition of Proposition 4.4 that carries
to this case (and in particular for small $t$). 

We shall recover the (inverse) density from the non-degenerate phase function.

Recall from Sect.2 that when $m=1$,
$\varphi$ is the eikonal action on $\Lambda_+$, which we complete
to $(\varphi,\psi,t)$ as eikonal coordinates. 

\medskip
\noindent {\it 5.2 Construction of the phase functions in the extended phase-space}
\smallskip
We begin to find a parametrization for $\Lambda_+$, then for $\Lambda_+^E$.
Among the $\theta$-variables, it is convenient to choose the eikonal.
\medskip
\noindent {\bf Proposition 5.5}: 
{\it Let $H(x,p)$ be positively
homogeneous of degree 1 with respect to $p$. We set $\theta = (\lambda, \varphi, \psi) \in {\bf R}^{n+1}$
Then
$$\Phi(x, t, \theta) = \varphi + \lambda \langle P(t,\varphi, \psi) , x - X(t,\varphi, \psi) \rangle\leqno(5.26)$$
is a generating family for
$\Lambda_+ \subset T^* {\bf R}^{n+1}_{x,t}$
at the points satisfying the inequality $\det(P,P_\psi)\neq0$, which holds at least for small $t$. Moreover $\Phi$ verifies
the initial condition $\Phi(x, 0, (1, \varphi, \psi))=\langle x,\omega(\psi)\rangle$. }
\smallskip
\noindent {\it Proof}:
Taking partial derivatives in (5.26) with respect to variables $\varphi,\psi,\lambda$, using 
with $\langle P,X_\psi\rangle=0$, $\langle P,X_\varphi\rangle=1$, we find
$$\eqalign{
&\partial_\psi\Phi=\lambda\langle P_\psi,x-X\rangle\cr
&\partial_\lambda\Phi=\langle P,x-X\rangle\cr
&\partial_\varphi\Phi=1-\lambda+\lambda\langle P_\varphi,x-X\rangle
}\leqno(5.27)$$
The critical point $x=X(t,\varphi,\psi)$ is uniquely determined (for small $t$)~: namely the 
determinant of the $n\times n$ system for the first 2 Eq. is given by $\det(P,P_\psi)$, and this is non zero at $t=0$.
So when $\lambda=1$, $(x,t)=(X(t,\varphi,\psi),t)$ belongs to the projection of the critical set $C_\Phi$. Taking differential on $C_\Phi$ gives
$$\eqalign{
&d\partial_\psi\Phi=\langle P_\psi,dx-dX\rangle\cr
&d\partial_\lambda\Phi=\langle P,dx-dX\rangle\cr
&d\partial_\varphi\Phi=-d\lambda+\langle P_\varphi,dx-dX\rangle
}$$
Identifying the coefficient of $d\lambda$ in the linear combination
$$\alpha d\partial_\psi\Phi+\beta d\partial_\lambda\Phi+\gamma d\partial_\varphi\Phi=0$$
we find $\gamma=0$, and as above the condition $\det(P,P_\psi)\neq0$ gives $\alpha=\beta=0$.

At last we check that for $\lambda=1$
$$\Phi|_{t=0}=\int P\,dX+\langle\omega(\psi),x\rangle-\varphi=\varphi+\langle\omega(\psi),x\rangle-\varphi=\langle\omega(\psi),x\rangle$$
so $\Phi|_{t=0}=\langle\omega(\psi),x\rangle$ satisfies the initial condition. $\clubsuit$. 
\medskip
Next we determine the (inverse) density on $\Lambda_+$, which results essentially from the same computations. Actually the non-vanishing of the
(inverse) density is equivalent to the linear independence of the differentials $d(\partial_{\theta_1}\Phi),\cdots,d(\partial_{\theta_N}\Phi)$. 
\smallskip
\noindent {\bf Proposition 5.6}: {\it Under hypotheses of Proposition 5.4, let $y=(t,\varphi,\psi)$ be the coordinates on $C_\Phi$.
Then the (inverse) density $d\mu=\iota^*(dy)$ on $\Lambda_+=\iota(C_\Phi)$ is the absolute value of
$F[\Phi,dy]|_{C_\Phi}$ defined (up to sign) as the quotient of two $(2n+2)$-forms}
$$F[\Phi,dy]|_{C_\Phi}=
{dt\wedge\,d\psi\wedge\,d\varphi\wedge\,d(\partial_\varphi\Phi)\wedge\,d(\partial_\lambda\Phi)\wedge\,d(\partial_\psi\Phi)\over
  dx\wedge\,dt\wedge\,d\psi\wedge\,d\varphi\wedge\,d\lambda}=
\pm\det(P,\partial_\psi P)\leqno(5.29)$$
\smallskip
\noindent {\it Proof}: We expand the numerator of (5.29) as $\omega_1+\omega_2$, where
$$\eqalign{
    &\omega_2=dt\wedge\,d\psi\wedge\,d\varphi\wedge\, \langle P_\varphi,dx-dX\rangle \wedge \langle P,dx-dX\rangle
      \wedge\, \langle P_\psi,dx-dX\rangle\cr
      &\omega_1=-dt\wedge\,d\psi\wedge\,d\varphi\wedge\,d\lambda \wedge \langle P,dx-dX\rangle
      \wedge\, \langle P_\psi,dx-dX\rangle}
$$
Write $dX=\dot X\,dt+X_\varphi\,d\varphi+X_\psi\,d\psi$. It is easy to see that $\omega_2$ is a sum of terms, each containing
twice one of the factors $dt$,$d\varphi$ or $d\psi$. Hence $\omega_2=0$.
In turn, for the same reason, $\omega_1$ reduces to
$$\omega_1=-dt\wedge\,d\psi\wedge\,d\varphi\wedge\,d\lambda\wedge \langle P,dx\rangle \wedge \langle P_\psi,dx\rangle$$
The last two products equal $\det(P,P_\psi)$. Simplifying the quotient on the RHS we readily get (5.29). $\clubsuit$
\medskip
Consider now $\Lambda_+^E$. 
Recall HJ equation for the phase function parametrizing $\Lambda_+^E$ as in (2.12)-(2.13), expressed in the $(t,x)$ variables as
$$\partial_t\Psi+H(x,\partial_x\Psi)=E, \  \Psi|_{t=0}=\langle x,\omega(\psi)\rangle\leqno(5.30)$$
We find as before an integral manifold $\Lambda_{\Psi,t}=\{p=\partial_x\Psi(t,x)\}\subset T^*{\bf R}^n$. As in Sect.4,
we look instead for a generating family, which is the 1-jet of $\Psi$ along $\Lambda_+^E$.

Let $H(x,p)$ be positively homogeneous of degree 1. Near a non glancing point $z_0\in\Lambda_0$,
and at the points satisfying the inequality $\det(P,P_\psi)\neq0$, which holds at least for small $t$, we could try
$$\widetilde\Phi^E(x,t,\varphi,\psi,\lambda)= \Phi(x,t,\varphi,\psi,\lambda)+Et\leqno(5.31)$$
with $(t,\varphi,\psi,\lambda)\in{\bf R}^{n+2}$ as new $\theta$-parameters.
and check this defines a generating family for $\Lambda_+^E$ near $z_0$. 
In general this is a difficult task, so
for simplicity, we assume that $\Lambda_+^E$ can be (locally) parametrized by $(t,\psi)$ which implies (locally) $\rho(x)$
radially symmetric. 

So let us take instead
$$\Phi^E(x,t,\varphi,\psi,\lambda)= Et + \lambda \langle P(t,\varphi, \psi) , x - X(t,\varphi, \psi) \rangle\leqno(5.32)$$
with $(t,\psi,\lambda)$ as new $\theta$-variables, and $\varphi$ as a parameter.
Compared to (5.26) this amounts to permute the ``twin variables'' $(t,\varphi)$,
which play a symmetric role in the expression $\varphi+Et$ of the eikonal on $\Lambda_+^E$. 

As in Proposition 5.5, 
taking partial derivatives in (5.32) with respect to variables $t,\psi,\lambda$, we find
$$\eqalign{
&\partial_\psi\Phi^E=\lambda\langle P_\psi,x-X\rangle\cr
&\partial_\lambda\Phi^E=\langle P,x-X\rangle\cr
&\partial_t\Phi^E=(1-\lambda)E+\lambda\langle\dot P,x-X\rangle\cr
}$$
This determines again the critical point $\lambda=1$, $x=X(t,\varphi,\psi)$ when $\det(P,P_\psi)\neq0$. So everything holds as if there
were no $\varphi$ parameter.

Concerning the (inverse) density on $\Lambda_+^E$ the same computations as in Proposition 5.6 give~:
\smallskip
\noindent {\bf Proposition 5.7}: {\it Let $H$ be as (1.25). Assume $\Lambda_+^E$ is (locally) parametrized by $(t,\psi)$ which implies that
$\rho$ is (locally) radially symmetric. Assume $\rho(0)\neq1/E$, so that 
$\varphi=\varphi(E)$ on $\partial\Lambda_+^E$, and the set of glancing points on $\Lambda_0$ projects onto $x=0$.
Let $\theta=(t,\lambda,\psi)$
and $y=(t,\psi)$ be local coordinates near $z_0$ on $C_{\Phi^E}$ with fixed $E$.
Then the (inverse) density $d\mu_+^E$ on $\Lambda_+^E$ is the absolute value of
the quotient of two $(2n+1)$-forms}
$$F[\Phi^E,dy]|_{C_{\Phi^E}}=
{dt\wedge\,d\psi\wedge\,d\dot\Phi^E\wedge\,d(\partial_\lambda\Phi^E)\wedge\,d(\partial_\psi\Phi^E)\over
  dx\wedge\,d\psi\wedge\,dt\wedge\,d\lambda}=
\pm E\det(P,\partial_\psi P)$$
\smallskip
\noindent {\it Remark}: The (inverse) density on $\Lambda_+^E$ has the same expression as this on $\Lambda_+$.
It may becomes singular for $t>0$ when $E$ approaches a critical value of ${\cal H}(\varphi,\psi)$. This reflects the fact
that $\det(P,P_\psi)$ (which is always equal to 1 on $\Lambda_0$) may vanish rapidly as $t$ takes positive values.

In the general case of an Hamiltonian positively homogeneous of degree 1, and away from the glancing points, we conjecture that we 
we need only $(n+1)-\theta$ -parameters to define $\Phi^E$, one of  them being a function of the ``twin variables'' $(t,\varphi)$
(possibly in some degenerate cases such as radial symmetric potentials, a function of $t$ alone),
depending on the chart in $\Lambda_0$ where $\varphi=\varphi(\tau,\psi)$, or $\psi_1=\psi_1(\tau,\varphi,\psi')$, as stated in Lemma 5.2.
This chart is in turn prescribed by the point $x$ where we are computing the solution, moving back along the trajectories up to $\Lambda_0$.

\medskip
\noindent {\it 5.3 Reduced parametrizations of $\Lambda_+^E$ in case of the ``conformal metric''}.
\smallskip

Let us find the critical point $(t,\psi)$ of $\Phi^E$, i.e.  
consider the 
system $\partial_{(t,\psi)}\Phi^E=0$. We have 
$$\partial^2_{t,\psi}\Phi^E=\pmatrix{{\partial^2\Phi\over\partial t^2}&{\partial^2\Phi\over\partial t\partial\psi}\cr 
{\partial^2\Phi\over\partial \psi\partial t}&{\partial^2\Phi\over\partial\psi^2}}=
\pmatrix{\langle \dot P,\dot X\rangle&\langle P_\psi,\dot X\rangle\cr\langle P_\psi,\dot X\rangle&\langle P_\psi,X_\psi\rangle}\leqno(5.33)$$
where $\langle P_\psi,\dot X\rangle=\langle X_\psi,\dot P\rangle$. Let 
$D(z(t))=\det \partial^2_{t,\psi}\Phi^E$.

For $t=0$, we have $\det D(z(0))=\varphi\langle -\partial_xH,\partial_pH\rangle-\langle\partial_pH,\omega(\psi)^\perp\rangle^2$,
so for Hamiltonian (1.25) this is non vanishing when $z(0)$ is an ordinary point (in the sense of Definition 1.4) where $\varphi\neq0$, and $D(z(0))=0$ if 
$z(0)$ is glancing. So when $D(z(t))\neq0$, implicit function theorem shows that $\partial_{t,\psi}\Phi^E=0$ is equivalent to
$t=t_0(x,\varphi),\psi=\psi_0(x,\varphi)$.

Take polar coordinates on $M$ of the form
$(r,\theta)$ such that $(r,\theta)=(\varphi,\psi)$ parametrize a point on $\Lambda_0$ near $\partial\Lambda_+^E$. We make the identification 
$x=(r,\theta)$.

For simplicity we restrict to an ordinary points. We show that near an ordinary point $z(0)\in\partial\Lambda_+^E$, $\rank d\pi_x(z(0))=2$. Namely we have:
\smallskip
\noindent {\bf Proposition 5.8}: {\it As in Lemma 5.3, let us assume  $\varphi=\varphi(\tau,\psi)$.
Assume also $D(z(0))\neq0$ i.e. $z(0)=\varphi\omega(\psi)$ is an ordinary point with $\varphi\neq0$.
Then for $x$ sufficiently close to $z(0)$,
the system
$\partial_{t,\psi}\Phi^E(x,t,\psi,\lambda=1)=0$ is equivalent to $t=t(x),\psi=\psi(x)$, and $\varphi=\varphi(x)$ when $x$ is sufficiently close to
$\varphi\omega(\psi)$. The critical value of $\Phi$ takes the form
$$\Psi(r,\theta;\varphi,\psi)=r\cos(\theta-\psi(x))+{\cal O}(|r-\varphi(x),\theta-\psi(x)|)=r+{\cal O}(|r-\varphi(x),\theta-\psi(x)|)\leqno(5.34)$$
In particular $\rank d\pi_x(z(0))=2$. }
\smallskip
\noindent {\it Proof}: 
Since $D(z(0))\neq0$, implicit function theorem shows that
$\partial_{t,\psi}\Phi^E_{t,\psi}=0$ 
(we omitted $\lambda=1$) has a unique solution
$t=t_0(x,\varphi),\psi=\psi_0(x,\varphi)$, i.e. $t=t_0(r,\theta,\varphi),\psi=\psi_0(r,\theta,\varphi)$
Differentiating $\partial_{t,\psi}\Phi^E=0$ along $\Lambda_+^E$ 
with respect to $r,\theta$ and $\varphi$ we find
$$\eqalign{
&\langle\dot P,\dot X\rangle {\partial t_0\over\partial r}+\langle\dot P,X_\psi\rangle{\partial\psi_0\over \partial r}=
\langle\dot P,\omega(\theta)\rangle\cr
&\langle P_\psi,\dot X\rangle {\partial t_0\over\partial r}+\langle P_\psi,X_\psi\rangle{\partial\psi_0\over\partial r}=
\langle P_\psi,\omega(\theta)\rangle\cr
&\langle\dot P,\dot X\rangle {1\over r}{\partial t_0\over\partial \theta}+\langle\dot P,X_\psi\rangle{1\over r}{\partial\psi_0\over \partial r}=
\langle\dot P,\omega^\perp(\theta)\rangle\cr
&\langle P_\psi,\dot X\rangle {1\over r}{\partial t_0\over\partial r}+\langle P_\psi,X_\psi\rangle{1\over r}{\partial\psi_0\over\partial r}=
\langle P_\psi,\omega^\perp(\theta)\rangle\cr
&\langle \dot P,\dot X\rangle{\partial t_0\over\partial\varphi}+\langle \dot P,X_\psi\rangle{\partial \psi_0\over\partial\varphi}=
-\langle\dot P,X_\varphi\rangle\cr
&\langle P_\psi,\dot X\rangle{\partial t_0\over\partial\varphi}+\langle P_\psi,X_\psi\rangle{\partial \psi_0\over\partial\varphi}
=-\langle P_\psi,X_\varphi\rangle\cr
}\leqno(5.35)$$
Using the relation $\langle \dot P,X_\psi\rangle=\langle P_\psi,\dot X\rangle$, the 3 sub-systems have determinant
$$D(z(t))=\det\pmatrix{\langle \dot P,\dot X\rangle&\langle P_\psi,\dot X\rangle\cr
\langle P_\psi,\dot X\rangle&\langle P_\psi,X_\psi\rangle}\neq0$$
So (5.37) has the unique solution
with the condition
$$t_0(\varphi,\psi,\varphi)=0, \quad \psi_0(\varphi,\psi,\varphi)=\psi\leqno(5.36)$$
On the other hand we know from Proposition 5.3 that on $C_{\Phi^E}$, i.e. when $x=X(t,\varphi,\psi)$, we have $\varphi=\varphi_1(x,\psi)$ and
$t=t_1(x,\psi)$. This gives~:
$$g(x,\psi)=\psi-\psi_0\bigl(x,\varphi_1(x,\psi)\bigr)=0$$
Compute $\partial_\psi g(x,\psi)=1-{\partial \psi_0\over\partial\varphi}{\partial\varphi_1\over\partial\psi}$ at $x=\varphi\omega(\psi)$.
Combining Proposition 5.3 and [DoNaSh,Lemma 7] (which still doesn't assume $\widetilde k=2$ to be the rank of $\widetilde\Lambda_+$),  we get 
$$\partial_\psi{t_1\choose\varphi_1}(x,\psi)=-\pmatrix{1&0\cr -H&{}^t P}{0\choose X_\psi}+
{\cal O}\bigl(\big|x-X\bigl(t_1(x,\psi),\varphi_1(x,\psi),\psi\bigr)\big|\bigr)$$
Since $\langle P,X_\psi\rangle=0$ , we get in particular 
$${\partial\varphi_1\over\partial\psi}(x,\psi)={\cal O}\bigl(\big|x-X\bigl(t_1(x,\psi),\varphi_1(x,\psi),\psi\bigr)\big|\bigr)=o(1)$$
Thus $\partial_\psi g(x,\psi)\neq0$ along $C_\Phi$ and implicit function theorem
gives $\psi=\psi_1(x)$. Sustituting into $t=t_1(x,\psi)=t_0(x,\varphi)$ we get also
$t=t_1(x,\psi_1(x))=t_0(x,\varphi_1(x,\psi_1(x))$. 
Substituting into (5.32) (where we have assumed $m=1$) gives the critical value (for $\lambda=1$) where all ``$\theta$-variables'' have been 
eliminated
$$\Psi(r,\theta)=Et_1(t,\theta)+\varphi_1(r,\theta)+\langle P(t_1,\varphi_1,\psi_1),r\omega(\theta)-X(t_1,\varphi_1,\psi_1)\rangle\leqno(5.37)$$
which we expand around $(r,\theta)=(\varphi,\psi)$ using (5.36). With 
$$t_1(t,\theta)=o(1), \ \varphi_1(t,\theta)=\varphi+o(1), \ P(t_1,\varphi_1,\psi_1)=\omega(\psi)+o(1), \ 
X(t_1,\varphi_1,\psi_1)=\varphi\omega(\psi)+o(1)$$
substituting into (5.37) we find $\Psi(r,\theta)=r\cos(\theta-\psi)+o(1)$,
where 
$$o(1)={\cal O}\bigl(\big|x-X\bigl(t_1(x,\psi),\varphi_1(x,\psi),\psi\bigr)\big|\bigr)={\cal O}(|r-\varphi,\theta-\psi|)$$
This proves (5.34) (no ``$\theta$-variables''). $\clubsuit$ 

\medskip
\noindent {\it 5.4 Example in the constant coefficient case}
\smallskip
To close this section we consider as in Remark 3.1, $H=-h^2\Delta$, $n=2$ and 
$$f_h(x)=J_0\bigl({|x|\over h}\bigr)
=(2\pi)^{-1}\int_{-\pi}^\pi e^{i\langle x,\omega(\psi)\rangle}\,d\psi=(2\pi)^{-1}\int_{-\pi}^\pi e^{i|x|\sin\psi}\,d\psi$$ 
The outgoing solution is given by the oscillating integral
$$u_h(x)=(4\pi^2h)^{-1/2}\int e^{i\Phi(x,y,\psi)/h}\bigl(\xi^2-E-i0\bigr)^{-1}\,d\psi\,dy\,d\xi$$
which we compute (formally) by stationary phase in $(y,\xi)$. The critical point is given by $y=x$, $\xi={x\over|x|}\sin\psi$, and 
$$u_h(x)=h^{1/2}\int_{-\pi}^\pi e^{i|x|\sin\psi/h}{1\over\sin^2\psi-E-i0}\,d\psi+{\cal O}(h^{3/2})\leqno(5.40)$$
and the leading term can be simply evaluated by contour integrals. 
The phase of course is the same as in Bessel function $J_0$. In fact $u_h(x)$ takes the exact value comuted in Example 5.1.

We can also consider more
general $f_h$ and stick in an amplitude of the form (see [DoMaNaTu])
$$A(x,\psi)={1\over2}\bigl(a(|x|,\psi)+a(-|x|,\psi)\bigr)+{\langle x,\omega(\psi)\rangle\over2|x|} \bigl(a(|x|,\psi)-a(-|x|,\psi)\bigr)$$
so that
$$f_h(x)=\bigl({i\over2\pi h}\bigr)^{1/2}\int_{-\pi}^\pi e^{i\langle x,\omega(\psi)\rangle}A(x,\psi)\,d\psi$$
As in (5.40) we get
$$u_h(x)=h^{1/2}\int_{-\pi}^\pi e^{i|x|\sin\psi/h}{A(x,\psi)\over\sin^2\psi-E-i0}\,d\psi+{\cal O}(h)\leqno(5.41)$$
When $A$ is independent of $\psi$, this can lead to significant simplifications. 
\bigskip
\noindent {\bf Appendix~:} {\it Lagrange immersions and global half-densities}
\medskip
Recall first some well-known properties of Lagrangian immersions (see e.g. [DoZh], [DoNaSh])~: 
\smallskip
\noindent {\bf Proposition A.1}: {\it Let $\iota:\Lambda\to T^*M$, 
be a Lagrangian immersion, parametrized on a canonical chart $U$ by 
$\varphi\mapsto z=\iota(\varphi)=(X(\varphi),P(\varphi))$. Introduce the 
Jacobian matrices $B(z)={\partial P\over\partial\varphi}$, $C(z)={\partial X\over\partial\varphi}$. Then:

(1) the matrix $(B(z), C(z))$ is of rank $n$. 

(2) the matrix ${}^tC(z)B(z)$ is symmetric.

(3) $C(z)\pm iB(z)$ is non degenerate.}
\smallskip
The symmetry of ${}^tC(z)B(z)$ expresses for instance in the situation of Sect.4 as the symmetry of Gram matrix
$$\pmatrix{\langle \dot X,\dot P\rangle&\langle \dot X,P_\psi\rangle\cr
\langle X_\psi,\dot P\rangle&\langle X_\psi,P_\psi\rangle}\leqno(A.1)$$

We consider the rank of projections $\pi_x:\Lambda_+\to M$. It is equal to the rank of $\pi_{x,t}:\widetilde\Lambda_+\to M\times{\bf R}_t$.  

In general we call {\it focal point} a point $z\in\Lambda$ where
$\pi_*:T\Lambda\to TM$ is singular, and {\it caustics} the projection ${\cal C}$ of the set
of focal points onto $M$. Assume $n=2$, and let $z$ be a focal point, so $C(z)$ cannot be of rank 2, and by property (1) above either
$B(z)$ is of rank 2 (and $C(z)$ has rank at most 1, since the projection $\pi:\Lambda_+\mapsto M$ is not a diffeomorphism at $z$)
or both $C(z)$ and $B(z)$ are of rank 1. 
\smallskip
Assume now $H(x,p)={|p|^m\over\rho(x)}$, $n=2$. Let $\Lambda=\Lambda_+$ be an integral manifold of $v_H$ in the energy shell $H(x,p)=E$,
and $U\subset \Lambda_+$ be a canonical chart parametrized by $\varphi=(t,\psi)$, i.e. 
$z=(X(t,\psi),P(t,\psi))$ verifies $\dot X=\partial_xH(X,P), \dot P=-\partial_pH(X,P)$, such that $U\to T^*M$ is an immersion
(not necessarily an embedding). Recall $\langle P,\dot X\rangle=mH$, $\langle P,X_\psi\rangle=0$. Actually 
$(X(t,\psi),P(t,\psi)$ may depend on additional parameters, as in Sect.5, but here only $t,\psi$ matter. Consider the quantity 
$\det(P,P_\psi)$. 

\medskip
\noindent {\bf Lemma A.2}: {\it Let $H(x,p)={|p|^m\over\rho(x)}$, $n=2$. 
Assume that at some point $z=\bigl( X(t,\psi),P(t,\psi)\bigr)$, we have $\det(P,P_\psi)=0$. 

1) If $\nabla\rho(X(t,\psi))\neq0$, then either $|P(t,\psi)|=1$, or $P_\psi(t,\psi)=0$.
In the latter case $C(z)=(\dot X,X_\psi)$ has rank 2, i.e. $\pi_x$ is regular at $z$.

2) If $\rho$ has a critical point at some $x_0=X(t,\psi)$,
then $\pi_x$ is regular at $z=(x_0,P(t,\psi))$ } 
\smallskip
\noindent {\it Proof}: We apply Proposition A.1 to the Lagrangian immersion
$\iota:U\to T^*M$, $\varphi=(t,\psi)$. 
Assume $\det(P,P_\psi)=0$ at some point $z=(X(t,\psi),P(t,\psi))$. Then either $P_\psi=0$, or 
by Hamilton equations (3.2), there is $\alpha\in{\bf R}$ such that $\dot X=\alpha P_\psi$. Let first $\nabla\rho\neq0$.

(i) Let $P_\psi=0$. In this case 
$X_\psi\neq0$, for otherwise this would contradict property (3) of Proposition A.1.
The symmetry of ${}^tC(z)B(z)$ shows that $\langle \dot P,X_\psi\rangle=0$, or
$\langle\nabla\rho(X),X_\psi\rangle=0$. Since $\langle P,X_\psi\rangle=0$, we find that $\dot X,P,\dot P,\nabla\rho(X)$ are parallel, 
and all orthogonal to $X_\psi$. In particular, $C(z)=(\dot X,X_\psi)$ is of rank 2. 

(ii) Let $P_\psi=0$. By property (2) of Proposition A.1, the matrix
$${}^tC(z)B(z)=\pmatrix{*&\dot X_1\partial_\psi P_1+\dot X_2\partial_\psi P_2\cr \dot P_1\partial_\psi X_1+\dot P_2\partial_\psi X_2&*}$$
has to be symmetric. If $\dot X=\alpha P_\psi$, this implies $\alpha=\langle \dot P,\partial_\psi X\rangle|P_\psi|^{-2}$, and
hence differentiating the dispersion relation (2.2)
$$\dot X= \langle \dot P,\partial_\psi X\rangle |P_\psi|^{-2} P_\psi={|P|^m\over\rho(X)^2}\langle \nabla\rho(X),X_\psi\rangle |P_\psi|^{-2} P_\psi
\leqno(A.5)$$ 

\noindent $\bullet$ Assume $B(z)$ is of rank 2, and $C(z)$ of rank at most 1.  
We know that $\langle P(t,\psi),\dot X(t,\psi)\rangle=m(E-\tau)\neq0$, $\dot X\neq0$ parallel to $P\neq0$, 
and $\langle P(t,\psi),\partial_\psi X(t,\psi)\rangle=0$. 

Assume $X_\psi\neq0$.
Since $C(z)$ is of rank at most 1, $X_\psi$ parallel to $\dot X$, which is itself parallel to $P$. So $P$ is parallel to $X_\psi$,
which contradicts $\langle P,X_\psi\rangle=0$. Hence $X_\psi=0$. By first equality (A.6) we have $\langle P,P_\psi\rangle=0$.
If $\det(P,P_\psi)=0$, then we would have $P$ both orthogonal to $X_\psi$, and parallel to $P_\psi\neq0$.
But $B(z)$ is of rank 2, which is a contradiction.

\noindent $\bullet$ So by property (3) of Proposition A.1
we must have $B(z)$ and $C(z)$ of rank 1. So either $\dot X=\lambda X_\psi$ and $\dot P=\mu P_\psi$ for some 
$\lambda,\mu\in{\bf R}$, or $X_\psi=0, \dot P=\mu P_\psi$, or $P_\psi=0, \dot X=\lambda X_\psi$, or $X_\psi=P_\psi=0$. 

Examine the first case:
Identifying the off-diagonal terms
of $^tB(z)C(z)$, which is symmetric by property (2) of Proposition A.1,
we find that either $\lambda=\mu$ or ${\partial X\over\partial\psi}
\perp {\partial P\over\partial\psi}$. 
But $\lambda\neq\mu$ since otherwise the complex
matrices $C(z)\pm iB(z)$ would be degenerate, which violates
property (3) of Proposition 3.2. So ${\partial X\over\partial\psi}$
and ${\partial X\over\partial t}$ are colinear, and orthogonal
to both ${\partial P\over\partial\psi}$ and 
${\partial P\over\partial t}$. Assume first $\lambda\mu\neq0$. 
Using again (A.5), we find that $P\perp P_\psi$, since $P_\psi\neq0$, we find that $\det(P,P_\psi)\neq0$.
Assume then $\lambda\mu=0$. If $\lambda=0, \mu\neq0$, we would have $\dot X=0$, which is impossible.
Let now $\lambda\neq0$, then $B+iC=(P_\psi+iX_\psi,i\lambda X_\psi)$ has rank 2, so has $(P_\psi,X_\psi)$.
So by (A.6) $P$, ${\partial X\over\partial\psi}$
and ${\partial X\over\partial t}$ are colinear, and orthogonal
to both ${\partial P\over\partial\psi}$ and 
$\nabla\rho$. As before, this implies $\det(P,P_\psi)\neq0$. 

Examine the second case $X_\psi=0, \dot P=\mu P_\psi$. Writing that ${}^tB(z)C(z)$ is symmetric, we find $\dot X\perp P_\psi$.
Since $B(z)+iC(z)$ has rank 2, $P_\psi\neq0$. Then $P$ and $\dot X$ are parallel, and both orthogonal to $\nabla\rho$, $P_\psi$ and $\dot P$. 
As before, we find $\det(P,P_\psi)\neq0$. 

The last two cases are similar. So $\det(P,P_\psi)\neq0$ at any focal point where $\nabla\rho\neq0$. .

(iii) Let at last, $\nabla\rho(x_0)=0$. The first situation above cannot hold, since this would imply $\dot X=0$ by (A.6), hence $P=0$. 
Hence $\dot P=P_\psi=0$ which implies again $C(z)$ of rank 2 by Proposition A.1. 
$\clubsuit$
\bigskip

\centerline{\bf References}
\medskip
\noindent [Ar] V.Arnold. {\bf 1}. 
On a characteristic class entering into conditions of quantization. Funct. Anal. Appl. 1, p.1-13, 1967.
{\bf 2}. M\'ethodes Math\'ematiques de la M\'ecanique Classique. Mir, Moscou, 1976.

\noindent [BaWe] S.Bates, A.Weinstein. Lectures on the geometry of quantization. Berkeley Math. Lect. Notes 88,
American Math. Soc. 1997.

\noindent [AnDoNaRo] A.Anikin, S.Dobrokhotov, V.Nazaikinskii, M.Rouleux. {\bf 1}. Maslov's canonical operator on a pair of 
Lagrangian manifolds and asymptotic solutions of stationary equations with localized right-hand sides.
Doklady Akad. Nauk, Vol. 76, No1, p.1-5, 2017.
{\bf 2}.
Semi-classical Green functions.
Proceedings ``Days of Diffraction 2018'', Saint-Petersburg, IEEE.
{\bf 3}. Lagrangian manifolds and the construction
of asymptotics for (pseudo)differential equations
with localized right-hand sides, Theor. Math. Phys. 214(1): p.1-23, 2023. 

\noindent [Bad] N.Baddour, Operational and convolution properties of two-dimensional Fourier transforms in polar coordinates,
J. Opt. Soc. Am. A, Vol.26, p.1767-1777, 2009.

\noindent [BeBa] M.V. Berry, N.L. Balazs. Non spreading wave packets. American J. Physics 47, p.264, 1979.

\noindent [BoRo] I.Bogaevskii, M.Rouleux. {\bf 1}. Lagrangian intersections and glancing points: typical transitions of phase in semiclassical
approximations. Proceedings Intern. Conference Days of Diffraction 2023, IEEE, p.12-18, 2023.
{\bf 2}. Typical transitions of phase in semi-classical approximations. Submitted.

\noindent [Bon] J.-F. Bony. Mesures limites pour l'\'equation de Helmholtz dans le cas non captif.
Ann. Fac. Sci. Toulouse Math. (6), 18(3): p.459-493, 2009

\noindent [Ca] F. Castella. The radiation condition at infinity for the high-frequency Helmholtz equation with source term: A wave-packet approach.
Journal of Functional Analysis 223(1), p.204-257, 2005.
    
\noindent [CdV] Y.Colin de Verdi\`ere. M\'ethodes semi-classiques et th\'eorie spectrale. 
https://www-fourier.ujf-grenoble.fr/~ycolver/ All-Articles/93b.pdf

\noindent [DoMaNa] S.Dobrokhotov, G.Makrakis, V.Nazaikinskii. {\bf 1}. Maslov's canonical operator, H\"ormander formula, 
and localization of Berry-Balazs solution in the theory of wave beams. Th. Math. Phys. 180(2), p.894-916, 2014.
{\bf 2}. Fourier integrals and a new representation of Maslov canonical operator near caustics. American Math. Soc. Transl. Vol. 233, P.95-115 2014 

\noindent [DoMaNaTu] S.Dobrokhotov, G.Makrakis, V.Nazaikinskii, T.Tudorovskii. New formulas for 
Maslov's canonical operator in a neighborhood of focal points and caustics in 2D semiclassical asymptotics. 
Th. Math. Phys. 177(2), p.1679-1605, 2013.

\noindent [DoMiNa] S.Dobrokhotov, D.Minenkov, V.Nazaikinskii. Representation of Bessel functions by the Maslov canonical operator.
Theoretical and Math. Physics, 208(2), p.1018-1037, 2021.

\noindent [DoMiRo] S.Dobrokhotov, D.Minenkov, M.Rouleux. The Maupertuis-Jacobi principle for Hamiltonians
$F(x,|p|)$ in 2-D stationary semiclassical problems. Math. Notes, Vol. 97, No.1, p.42-49, 2015.

\noindent [DoNa] S.Dobrokhotov, V.Nazaikinskii. {\bf 1}. Punctured Lagran- gian manifolds and asymptotic solutions of the
linear water-wave equations with localized initial solutions.
Math. Notes, 101, No.6, p.130-137, 2017. {\bf 2}. {\it Private communication}.

\noindent [DoNaSh] S.Dobrokhotov, V.Nazaikinskii, A.Shafarevich. New integral representations of Maslov canonical
operator in singular charts. Izv. Math. 81:2, p.286-328, 2017.

\noindent [DoNaTs] S.Dobrokhotov, V.Nazaikinskii, A.Tsvetkova.
Asymptotics of the localized Bessel beams and Lagrangian manifolds.
J. Communications Technology and Electronics, Vol.26 (6), p.625-638, 2023.

\noindent [DoRo] S.Dobrokhotov, M.Rouleux.
The semi-classical Maupertuis-Jacobi correspondence for quasi-periodic Hamiltonian flows with applications to linear 
water waves theory.
Asympt. Analysis, Vol.74 (1-2), p.33-73, 2011.

\noindent [Dui] J.J. Duistermaat. {\bf 1}. Oscillatory integrals, Lagrange immersions and unfolding of singularities. 
Comm. Pure Appl. Math. 27, p.207-281, 1974. {\bf 2}. Fourier Integral Operators. Birkh\"auser.

\noindent [Dui] J.J. Duistermaat, L.H\"ormander. Fourier Integral Operators II. Acta Mathematica. 128, p.183-269, 1972. 

\noindent [ElGr] Y.Eliashberg, M.Gromov. Lagrangian intersections theory. A finite dimensionnal approach. 
https://www.ihes.fr/~gromov/wp-content/uploads/2018/08/998.pdf

\noindent [EstHagHedLitt] Ilya Esterlis, Hal M.Haggard, Austin Hedeman, Robert G.Littlejohn. Maslov indices, Poisson brackets, and singular differential forms.
arXiv:1402.0786v1. 

\noindent [Fe] M.V. Fedoriuk (ed.) Partial Differential Equations V
Asymptotic Methods for Partial Differential Equations. The Encyclopedia of mathematical sciences 34. Springer, 199.

\noindent [ForHasHil] G.A.Ford, A.Hassell, L.Hillairet. Wave propagation on euclidean surfaces with conical singularities I: Geometric Diffraction.
arXiv:1505.01043v1

\noindent [GreUhl]  A.Greenleaf,  G.Uhlmann Composition of some singular Fourier integral operators and estimates for restricted X-ray transforms
Ann. Inst. Fourier 40(2), p.443-466, 1990.

\noindent [GriSj] A.Grigis, J.Sj\"ostrand. Microlocal analysis for differential operators. An Introduction. Cambridge Univ. Press, 1994.

\noindent [GuMel] V.Guillemin, R.Melrose. Poisson formula for manifolds with boundary. Adv. in Math. 32(3), p.204-232, 1979. 

\noindent [GuSt1] V.Guillemin, S.Sternberg. {\bf 1} Geometric Asymptotics.
American Math. Soc. Surveys, 14, Providence, Rhode Island, 1977. 
{\bf 2}. Semi-Classical Analysis. International Press, Boston, 2013.

\noindent [H\"o] L.H\"ormander. The Analysis of Linear Partial Differential Operators I-IV. Springer.

\noindent [Iv] V.Ivrii. Microlocal Analysis and Precise Spectral Asymptotics. Springer-Verlag, Berlin, 1998.

\noindent [JoMeRa] J.L. Joly, G.M\'etivier and J.Rauch. {\bf 1} Non linear oscillations beyond caustics.
Comm. Pure and Appl. Math., Vol. XLIX. p.443-527, 1996. {\bf 2}. Recent results in non-linear geometric optics, {\it in}
Hyperbolic Problems: Theory, Numerics, Applications, 
ISNM, Vol.130,  pp 723–736, Springer-Verlag, 2021. 

\noindent [Jos] M.S. Joshi. Geometric proofs of composition theorems for generalized Fourier Integral Operators. 
Portugaliae Mathematica, 56(2), p.1-26, 1999. 

\noindent [KaMa] M.V Karasev, V.P. Maslov. Asymptotic and geometric quantization. Russ. Math. Surve. 39, p.133-205, 1984.

\noindent [KlCa] A.Klak, F.Castella. Radiation condition at infinity for the high-frequency Helmholtz
equation: optimality of a non-refocusing criterion. Hokkaido Mathematical Journal, 2014, 43 (3),
pp.275-325.

\noindent [Ku] V.Kucherenko. Quasi-classical asymptotics of a point source function for the stationary Schr\"o- dinger equation.
Teor. i Mat. Fiz. Vol.1, No.3, p.384-406. Transl. Consultants Bureau, 1970. 

\noindent [Ler] J.Leray. Analyse lagrangienne 
et M\'ecanique Quantique. S\'eminaire EDP Coll\`ege de France, 1976-77.

\noindent [MRa] J E. Marsden, T.S. Ratiu. Introduction to Mechanics and Symmetry. Springer, 1994.

\noindent [Ma] V.P.Maslov. {\bf 1}. Th\'eorie des perturbations et m\'ethodes asymptotiques. Dunod, Paris, 1972.
{\bf 2}. Operational Methods. Moscow: Mir~Publ. 1976.

\noindent [MeSj] A. Melin, J.Sj\"ostrand. Fourier integral operators with complex valued phase functions, Springer Lect. Notes in Maths, No459, 
p.120-223.

 \noindent [MelUh] R.~B. Melrose, G.~A. Uhlmann,
Lagrangian intersection and the Cauchy problem,
Comm. Pure Appl. Math. 32 (4), p.483-519, 1979. 

\noindent [MiBeSa] J.~ C. Mi\~nano, P. Benitez, A. Santamaria. Hamilton-Jacobi equation in momentum space. OSA. Vol. 14, No. 20

\noindent [ReMiKaDo] K.Reinjders, D.Minenkov, M.Katznelson, S.Dobrokhotov. Electronic optics in graphene in the
semiclassical approximation. Annals of Physics 397 p.65-135, 2018. 

\noindent [PeSto] V.Petkov, L.Stoyanov. Geometry of Generalized Geodesic Flows and Inverse spectral problems, 2:nd edition, Wiley, 2024.

\noindent [PeVo]  V.Petkov, G.Vodev, Asymptotics of the number of the interior transmision eigenvalues, J. Spectral
Theory 7, p.1-31, 2017. 

\noindent [So] J.M.Souriau, Construction explicite de l'indice de Maslov. Applications, {\it in}:
Group theoretical methods in Physics, Nijmegen. Lect. Notes in Physics 50, Springer, 1976.

\noindent [SteSh] B.Sternin, V.Shatalov. On a method of solving equations with simple characteristics. Sb. Math. 44(1), p.23-59, 1983.

\noindent [UhlZh] G.Uhlmann, J.Zhai. Inverse problems for nonlinear hyperbolic equations.  Discrete and Continuous Dynamical Syst. 41(1),
p.455-469,  2021. 

\noindent [Vo] G.Vodev. {\bf 1} Transmission eigenvalue-free regions. Commun. Math. Phys. 336, p. 1141-1166, 2015.
{\bf 2}. Parabolic transmission eigenvalue-free regions in the degenerate isotropic case. Asymptot. Anal. 106, no. 3-4, p.147-168, 2018. 
{\bf 3} Improved parametrix in the glancing region for the interior Dirichlet-to-Neumann map,  Comm. Partial Differential Equations 44, no. 5, p.367-396,
2019. {\bf 4}. Semiclassical parametrix for the Maxwell equation and applications to the electromagnetic transmission eigenvalues,
Res. Math. Sci. 8 (3), Paper 35, 2021.

\noindent [Ta] M.Taylor. Finsler structures and wave propagation, {\it in}: V.Isakov (ed.), Sobolev spaces in Mathematics III.
Int. Math. Series, Springer, 2009.

\bye